\documentclass[english,final,a4paper,11pt,reqno]{amsart}
\usepackage[hmargin=3cm,vmargin=3.5cm]{geometry}
\usepackage{amssymb,graphicx,color,enumitem}
\usepackage[plainpages=false,hypertexnames=false,colorlinks=true,bookmarksopen=true,bookmarksopenlevel=1]{hyperref}
\usepackage{amsmath,epsfig,psfrag,subfigure,tikz}

\renewcommand{\ldots}{\dotsc}

\def\ds{\displaystyle}

\newcommand{\R}{\mathbb{R}}
\newcommand{\C}{\mathbb{C}}

\newcommand{\dd}[2]{{\displaystyle \frac{d #1}{d #2}}}

\numberwithin{equation}{section}

\begin{document}

\title[periodic normal forms for codim 2 bifurcations of limit cycles]{Numerical periodic normalization for codim 2 bifurcations of
limit cycles with center manifold of dimension higher than $3$}

\author[V. De Witte]{Virginie De Witte}
\address{V. De Witte\\Department of Applied Mathematics and Computer Science\\
Ghent University\\ Krijgslaan 281-S9\\ B-9000 Gent\\ Belgium}
\email{Virginie.DeWitte@UGent.be}

\author[W. Govaerts]{Willy Govaerts}
\address{W. Govaerts\\\\Department of Applied Mathematics and Computer Science\\
Ghent University\\ Krijgslaan 281-S9\\ B-9000 Gent\\ Belgium}
\email{Willy.Govaerts@UGent.be}

\author[Yu.A. Kuznetsov]{Yuri A. Kuznetsov}
\address{Yu.A. Kuznetsov\\ Mathematical Insitute\\ Utrecht University\\ The Netherlands \and Department of Applied Mathematics\\ University of Twente\\ The Netherlands}
\email{I.A.Kouznetsov@uu.nl}

\author[H.~G.~E. Meijer]{Hil Meijer}
\address{H.~G.~E. Meijer\\MIRA Inst Biomed Technol \& Tech Med\\ Twente University\\ 7500 AE Enschede\\ The Netherlands}
\email{meijerhge@math.utwente.nl}


\date{\today}

\begin{abstract}
Explicit computational formulas for coefficients of the periodic normal forms of the 
three most complex codim $2$ bifurcations of limit cycles with dimension of the center 
manifold equal to $4$ or to $5$ in generic autonomous ODEs are derived. The resulting 
formulas are independent of the dimension of the phase space and involve solutions of 
certain boundary-value problems as well as multilinear functions from the Taylor expansion 
of the ODE right-hand side near the cycle. The formulas allow one to distinguish between 
the complicated bifurcation scenarios which can happen near these codim $2$ bifurcations 
of limit cycles, where $3$-tori and $4$-tori can be present. We apply our techniques to 
the study of a known laser model, a novel model from population biology, and one for 
mechanical vibrations; these models exhibit Limit Point--Neimark-Sacker, Period-Doubling--Neimark-Sacker 
and double Neimark-Sacker bifurcations. Lyapunov exponents are computed to numerically 
confirm the results of the normal form analysis, in particular with respect to the existence 
of stable invariant tori of various dimensions and chaos.
\end{abstract}

\maketitle

\section{Introduction}
\label{sec1}
Consider a smooth system of ODEs 
\begin{equation} \label{ODE}
\dot{x}=f(x,p), ~~~x\in\mathbb{R}^{n},
\end{equation}
smoothly depending on a parameter vector $p\in\mathbb{R}^{m}$. Typically, the dynamics of 
such systems show qualitative transitions, i.e. bifurcations, upon variation of a parameter. It is hard to use 
simulations to characterize such transitions correctly and efficiently. Numerical continuation 
software such as {\sc auto} \cite{AUTO97} or {\sc matcont} \cite{sac2003,MATCONT,NewMATCONT} may be used 
to track bifurcations from a stable equilibrium to a periodic oscillation by a Hopf bifurcation 
and even the appearance of (un)stable invariant tori with multi-frequency oscillations by a 
secondary Hopf, or Neimark-Sacker bifurcation. Bifurcations of these invariant tori $\mathbb{T}^{m\geq2}$ 
into other tori or chaos, however, are out of reach of the standard numerical analysis. 

One possibility to study bifurcations of tori -- if they are stable -- is to compute Lyapunov exponents. 
The dimension of the torus for a given parameter value then equals the number of exponents equal to zero. 
Varying one parameter one can observe that exponents become zero and this indicates a bifurcation. The 
exact nature of the bifurcation is however obscured from this analysis and should be elucidated with 
additional means. Yet, in many cases, bifurcations of tori first emerge from codim 2 bifurcations of 
limit cycles. Specifically, these codim 2 bifurcations are points in the parameter plane where one 
Neimark-Sacker bifurcation curve intersects a Limit Point of cycles, a Period-Doubling or another 
Neimark-Sacker bifurcation curve. The intersections produce LPNS, PDNS, or NSNS bifurcations, 
respectively. This paper focuses on these bifurcations, occuring in generic systems (\ref{ODE}) when $m \geq 2$
and $n$ is sufficiently large. The bifurcations are well understood theoretically with Poincar\'e maps 
and the corresponding normal forms \cite{Ar:83,Io:79,GuHo:83,Ku:2004,KuMe:2006,GoGhKuMe:07,BSV:2010}. 
The results of the analyis of the normal form for these codim 2 bifurcations can be used to verify 
nondegeneracy conditions and classify the bifurcation structure. Hence, we need an algorithm for 
the numerical computation of the coefficients of each critical normal form to enable this analysis. 

There is a straightforward approach to obtain the critical normal forms of the codim 2 bifurcations 
of the limit cycle. In the Poincar\'e map, the limit cycle is a fixed point and one can use techniques 
developed for maps to obtain the critical normal form \cite{KuMe:2006,GoGhKuMe:07}. However, in this case 
partial derivatives of the map up to order $k$, most often $k=3$, sometimes $k=5$, are needed. This 
may be done using software such as {\sc capd} \cite{capd} or {\sc tides} \cite{tides}. These packages 
can compute the solution and the derivatives of the solution with respect to the initial condition with 
arbitrary precision using Taylor series. Alternatively one could integrate the variational equations 
\cite{Simo} or use automatic differentiation \cite{Guckenheimer:2000, KuMe:2005} to obtain the derivatives
of the Poincar\'e map. All these methods, however, have two drawbacks that make them less (time) efficient. 
First, these are shooting methods that are slower when the system is very sensitive to perturbations. 
Second, the full Poincar\'e map is computed while only certain expressions are needed for the normalization. 
There is an alternative technique that is more suitable in the context of numerical continuation of 
periodic orbits using collocation as the whole periodic orbit is available. It uses periodic normalization 
\cite{Io:88,IoAd:92} and has been applied to codim 1 bifurcations of limit cycles and implemented in {\sc matcont} 
\cite{KuDoGoDh:05}. Recently, we have extended this algorithm to codim 2 bifurcations of limit cycles with 
center manifold dimension at most 3 \cite{drdwgk}. Here we consider the three remaining cases, 
LPNS, PDNS, and NSNS, that are characterized by a center manifold of the critical cycle of
dimension 4 or 5. These three cases always involve a -- possibly unstable -- two-dimensional torus 
$\mathbb{T}^{2}$.

We have implemented our algorithm in the numerical continuation toolbox {\sc matcont} which automatically 
invokes the algorithm whenever the corresponding bifurcation is detected. Hence, any user is able to use 
it and take advantage of the automated normal form analysis. Here we document precisely what our algorithm 
does. First, its aim is to compute coefficients of a periodic critical normal form. We present these normal 
forms in Section \ref{Section:Normalforms} using (contrary to \cite{KuDoGoDh:05,drdwgk}) the original
Iooss \cite{Io:88} representation. Remark that these normal forms are closely related to 
the normal forms for the Zero-Hopf and Hopf-Hopf bifurcations of equilibria. We discuss the 
correspondence and the interpretation of the bifurcation diagrams of the generic unfoldings 
for the LPNS, PDNS, and NSNS bifurcations. Next, we present the formulas to compute 
the critical normal form coefficients in Section \ref{Section:Theory}. Here we also comment on the 
implementation which is similar to \cite{drdwgk}. Finally in Section \ref{Section:Examples}, 
we consider several examples that involve tori bifurcations: a laser model, a model from population biology, 
and one for mechanical vibrations. In these models we find and analyze the three codim 2 bifurcations 
that we focus on. We compute the critical normal form coefficients using our algorithm to predict the 
bifurcation diagram near each of these codim 2 points. Next we corroborate the predictions using 
Lyapunov exponents. In fact, we argue that the classification from the critical normal form guides the 
correct interpretation of the Lyapunov exponents. 

\section{Normal forms on the center manifold and their bifurcations}
\label{Section:Normalforms}
Write (\ref{ODE}) at the critical parameter values as
\begin{equation}\label{eq:P.1}
\dot{u}=F(u)
\end{equation}
and suppose that there is a limit cycle $\Gamma$ corresponding to a
periodic solution $u_0(t)=u_0(t+T)$, where $T>0$ is its (minimal)
period. Expand $F(u_0(t)+v)$ into the Taylor series
\begin{equation} \label{eq:MULT}
\begin{array}{rcl}
F(u_0(t)+v) &=& F(u_0(t)) + \\
&&A(t)v + {\displaystyle \frac{1}{2}B(t;v,v) + \frac{1}{3!}C(t;v,v,v)} + \\
&&{\displaystyle\frac{1}{4!}D(t;v,v,v,v) +
\frac{1}{5!}E(t;v,v,v,v,v) + O(\|v\|^6)},
\end{array}
\end{equation}
where $A(t)=F_u(u_0(t))$ and
$$
B(t;v_1,v_2)=F_{uu}(u_0(t))[v_1,v_2],~~
C(t;v_1,v_2,v_3)=F_{uuu}(u_0(t))[v_1,v_2,v_3],
$$
etc. The matrix $A$ and the multilinear forms $B,C,D,$ and $E$ are periodic in $t$
with period $T$ but this dependence will often not be indicated
explicitly.
\par
Consider the initial-value problem for the fundamental matrix
solution $Y(t)$, namely,
\begin{equation*} \label{P.2}
\frac{dY}{dt}=A(t)Y,\ \ Y(0)=I_n,
\end{equation*}
%f
where $I_n$ is the $n\times n$ identity matrix. The eigenvalues of
the monodromy matrix $M=Y(T)$ are called ({\em Floquet}) {\em
multipliers} of the limit cycle. The multipliers with $|\mu|=1$ are
called {\em critical}. There is always a ``trivial" critical
multiplier $\mu_n=1$. We denote the total number of critical
multipliers by $n_c$ and assume that the limit cycle is
non-hyperbolic, i.e. $n_c > 1$. In this case, there exists an
invariant $n_c$-dimensional {\em critical center manifold}
$W^c(\Gamma) \subset {\mathbb R}^n$ near $\Gamma$\footnote{This
manifold should not be confused with the $(n_c-1)$-dimensional center
manifold of the corresponding Poincar\'{e} map.}.

\subsection{Critical normal forms}
It is well known \cite{Ar:83,Ku:2004} that in generic
two-parameter systems (\ref{ODE}) only eleven codim 2 local bifurcations
of limit cycles occur. To describe the normal forms of (\ref{eq:P.1}) on the
critical center manifold $W^c(\Gamma)$ for these codim 2 cases, we
parameterize  $W^c(\Gamma)$ near $\Gamma$ by $(n_c-1)$ transverse coordinates and 
$\tau \in [0,kT]$ for $k \in \{1,2,3,4\}$, depending on the bifurcation. The $8$ cases 
where $n_c\leq 3$ were treated in \cite{drdwgk}. Based on \cite{Io:88} we show 
in Appendix \ref{Appendix:1} that the restriction of (\ref{eq:P.1}) to the
corresponding critical center manifold $W^c(\Gamma)$ with $n_c=4$ or $n_c=5$ will 
take one of the following {\it Iooss normal forms}. 

\subsubsection{LPNS}
The {\em Limit Point -- Neimark-Sacker} bifurcation occurs when the
trivial critical multiplier $\mu_n=1$ corresponds to a two-dimensional Jordan block and there are only
two more critical simple multipliers $\mu_{1,2}=e^{\pm i\theta}$
with $\theta\neq \frac{2 \pi}{j}$, for $j=1,2,3,4$. The four-dimensional Iooss normal form at the LPNS 
bifurcation is derived in Appendix \ref{Appendix:1_lpns} and can be written as
\begin{equation} \label{eq:NF-LPNS}
\begin{cases}
 \ds\dd{\tau}{t}=1-\xi_1+\alpha_{200} \xi_1^2+\alpha_{011}\left|\xi_2\right|^2+\alpha_{300} \xi_1^3+\alpha_{111} \xi_1\left|\xi_2\right|^2+\ldots, \smallskip \\
 \ds\dd{\xi_1}{\tau}= a_{200} \xi_1^2+a_{011}\left|\xi_2\right|^2+a_{300} \xi_1^3+a_{111} \xi_1\left|\xi_2\right|^2+\ldots, \smallskip \\
 \ds\dd{\xi_2}{\tau}= i\omega \xi_2+b_{110} \xi_1\xi_2+b_{210} \xi_1^2 \xi_2+b_{021} \xi_2\left|\xi_2\right|^2 +\ldots, \smallskip\\
\end{cases}
\end{equation}
where $\tau \in [0,T]$, $\omega = \theta/T$, $\xi_1$ is a real coordinate and $\xi_2$ is a complex
coordinate on $W^c(\Gamma)$ transverse to $\Gamma$, $\alpha_{ijk},a_{ijk}\in \R, b_{ijk} \in \C$, and 
the dots denote the $O(\|\xi^4\|)$-terms, which are $T$-periodic in $\tau$. The equations (\ref{eq:NF-LPNS}) 
implicitly describe motions on the $4$-dimensional invariant manifold $W^c(\Gamma)$ with one cyclic 
coordinate $\tau$.

\subsubsection{PDNS}
The {\em Period-Doubling -- Neimark-Sacker} bifurcation occurs when the
trivial critical multiplier $\mu_n=1$ is simple and there are only
three more critical simple multipliers, namely $-1$ and $\mu_{1,2}=e^{\pm i\theta}$
with $\theta\neq \frac{2 \pi}{j}$, for $j=1,2,3,4$. The four-dimensional Iooss normal form at the PDNS 
bifurcation is derived in Appendix \ref{Appendix:1_pdns} and can be written as
\begin{equation} \label{eq:NF-PDNS}
\begin{cases}
 \ds\dd{\tau}{t}=1+\alpha_{200} \xi_1^2+\alpha_{011}\left|\xi_2\right|^2+\alpha_{400} \xi_1^4+\alpha_{022}\left|\xi_2\right|^4+\alpha_{211} \xi_1^2\left|\xi_2\right|^2+\ldots, \smallskip \\
 \ds\dd{\xi_1}{\tau}= a_{300} \xi_1^3+a_{111}\xi_1\left|\xi_2\right|^2+a_{500} \xi_1^5+a_{122} \xi_1\left|\xi_2\right|^4+a_{311} \xi_1^3\left|\xi_2\right|^2+\ldots, \smallskip \\
 \ds\dd{\xi_2}{\tau}= i\omega \xi_2+b_{210} \xi_1^2\xi_2+b_{021} \xi_2\left|\xi_2\right|^2+b_{410} \xi_1^4\xi_2+b_{221} \xi_1^2\xi_2\left|\xi_2\right|^2+b_{032} \xi_2\left|\xi_2\right|^4+\ldots, \smallskip\\
\end{cases}
\end{equation}
where $\tau \in [0,2T]$, $\omega = \theta/T$, $\xi_1$ is a real coordinate and $\xi_2$ is a complex
coordinate on $W^c(\Gamma)$ transverse to $\Gamma$, $\alpha_{ijk},a_{ijk}\in \R, b_{ijk} \in \C$, and 
the dots denote the $O(\|\xi^6\|)$-terms, which are $2T$-periodic in $\tau$. The equations 
(\ref{eq:NF-PDNS}) implicitly describe motions on the $4$-dimensional invariant manifold $W^c(\Gamma)$ that is 
doubly covered by the selected coordinates.

\subsubsection{NSNS}
The {\em double Neimark-Sacker} bifurcation occurs when the
trivial critical multiplier $\mu_n=1$ is simple and there are only
four more critical simple multipliers $\mu_{1,4}=e^{\pm i\theta_1}$ and $\mu_{2,3}=e^{\pm i\theta_2}$
with $\theta_{1,2}\neq \frac{2 \pi}{j}$, for $j=1,2,3,4,5,6$ and $l\theta_1 \neq j 
\theta_2$ for $l,j\in \mathbb{Z}$ with $l+j \leq 4$ (see \cite{GoGhKuMe:07}). 
The five-dimensional periodic normal form at the  NSNS bifurcation is derived 
in Appendix \ref{Appendix:1_nsns} and can be written as
\begin{equation} \label{eq:NF-NSNS}
\left\{\begin{array}{rcl} \ds\dd{\tau}{t}&=&1+\alpha_{1100}\left|\xi_1\right|^2+\alpha_{0011}\left|\xi_2\right|^2\\
&&+\alpha_{2200}\left|\xi_1\right|^4+\alpha_{0022} \left|\xi_2\right|^4+\alpha_{1111} \left|\xi_1\right|^2\left|\xi_2\right|^2+\ldots, \\
 \ds\dd{\xi_1}{\tau}&=&i\omega_1 \xi_1 + a_{2100} \xi_1\left|\xi_1\right|^2+a_{1011}\xi_1\left|\xi_2\right|^2\\
 && + a_{3200}\xi_1\left|\xi_1\right|^4+a_{1022}\xi_1\left|\xi_2\right|^4+a_{2111}\xi_1\left|\xi_1\right|^2\left|\xi_2\right|^2+\ldots, \\
 \ds\dd{\xi_2}{\tau}&=& i\omega_2 \xi_2 + b_{0021} \xi_2\left|\xi_2\right|^2+b_{1110}\xi_2\left|\xi_1\right|^2\\
 && +b_{0032}\xi_2\left|\xi_2\right|^4+b_{2210}\xi_2\left|\xi_1\right|^4 +b_{1121}\xi_2\left|\xi_1\right|^2\left|\xi_2\right|^2+\ldots, \\
\end{array}
\right.
\end{equation}
where $\tau \in [0,T]$, $\omega_{1,2} = \theta_{1,2}/T$, $\xi_1$ and $\xi_2$ are complex
coordinates on $W^c(\Gamma)$ transverse to $\Gamma$, $\alpha_{ijkl}\in \R,a_{ijkl},b_{ijkl} \in \C$, 
and the dots denote the $O(\|\xi^6\|)$-terms, which are $T$-periodic in $\tau$. The equations 
(\ref{eq:NF-NSNS}) implicitly describe motions on a $5$-dimensional manifold with one cyclic 
coordinate $\tau$.

\subsection{Generic unfoldings of the critical normal forms}
\label{Appendix:2}
Here we describe how the coefficients of the critical normal forms can be used to predict
bifurcations of the phase portraits near the critical limit cycles for nearby parameter values.
We introduce certain quantities -- computable in terms of these coeffcients -- that are reported in 
the {\sc matcont} output and used to distinguish between various bifurcation scenarios in examples
in Section \ref{Section:Examples}. 

In generic two-parameter systems (\ref{ODE}) the considered bifurcations occur at isolated parameter values.
By translating the origin of the parameter plane to one of such points, we can consider an {\em unfolding}
of the corresponding bifurcation and study its canonical local bifurcation diagram for nearby parameter values.
It is well known that the critical center manifold $W^c(\Gamma)$ can be smoothly continued w.r.t. $p$ 
in a neighborhood of the bifurcation point, so that the restriction of (\ref{ODE}) to this manifold can 
be studied. Choosing appropriate coordinates $(\xi,\tau)$ on this parameter-dependent invariant manifold, 
one can transform the restricted system into a parameter-dependent normal form in which 
$\frac{d\xi}{d\tau}$ has a $\tau$-independent principle part and higher-order
terms which are $kT$-periodic in $\tau$ with $k=1$ for LPNS and NSNS and $k=2$ for PDNS.
Below we describe bifurcations of these principle parts, i.e., the truncated parameter-dependent
autonomous normal forms. Since the dynamics is determined by the $\xi$-equations, we first focus on 
their bifurcations and then interpret appearing bifurcation diagrams for the original system (\ref{ODE}). 
The new unfolding parameters will be denoted by $(\beta_1,\beta_2)$.

\subsubsection{LPNS}\label{Appendix:2LPNS}
\begin{figure}[htb]
\centering
 \subfigure[$s=1, \theta>0$]{
 \includegraphics[width=.48\textwidth]{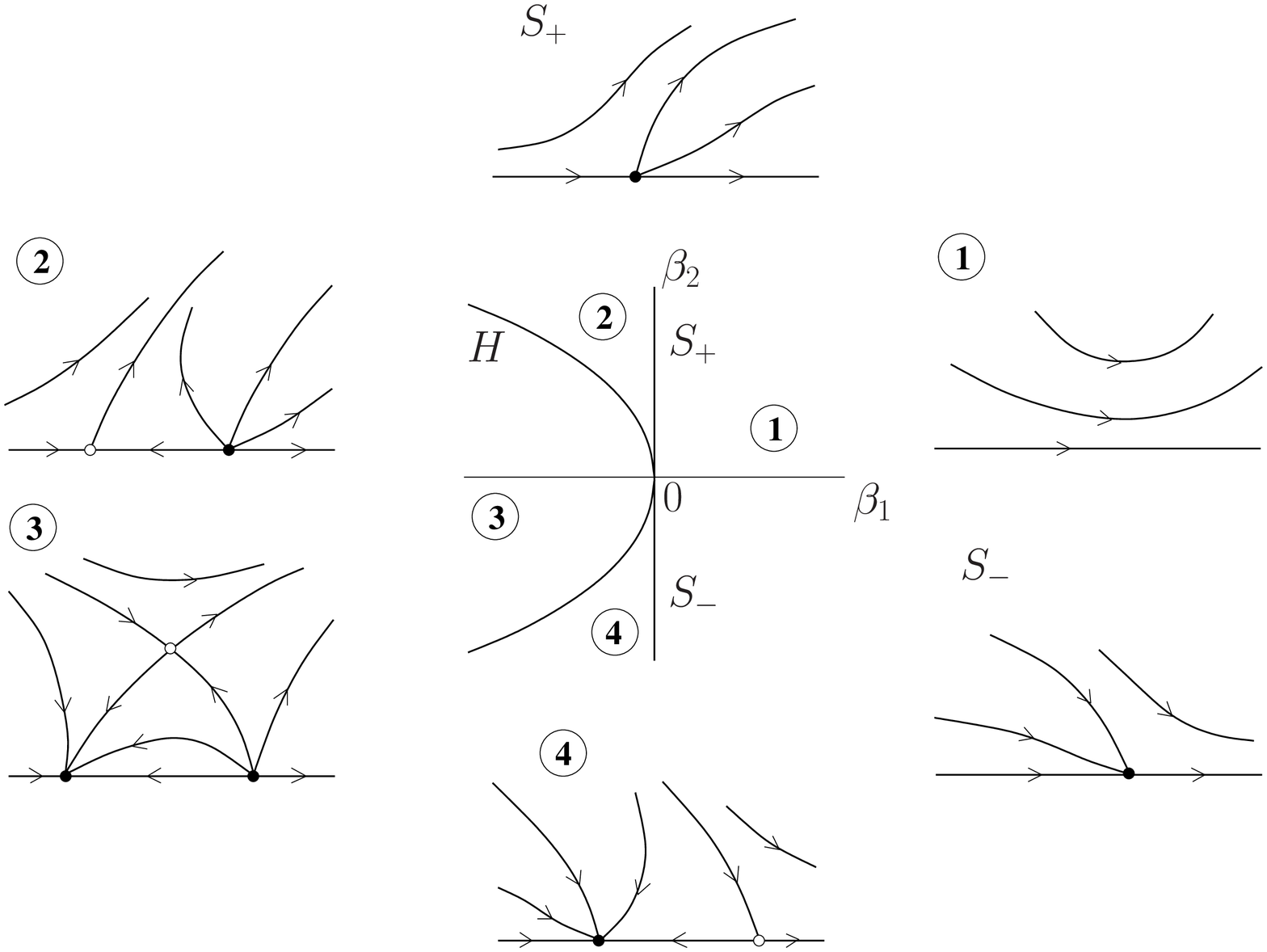}\label{fig:NF_LPNS_1}}\subfigure[$s=-1, \theta<0$]{
 \includegraphics[width=.48\textwidth]{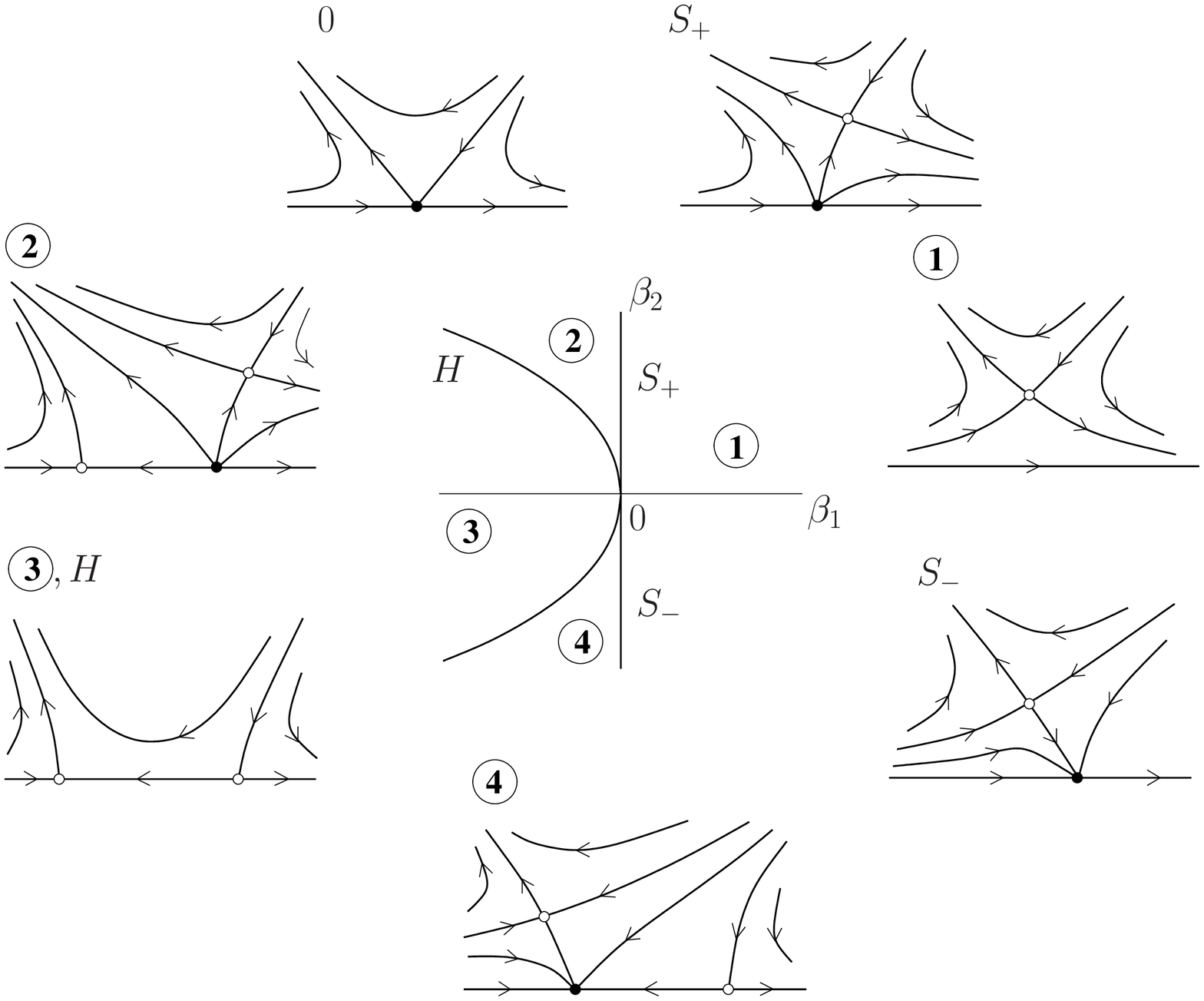}}
 \subfigure[$s=1, \theta<0$]{
 \includegraphics[width=.48\textwidth]{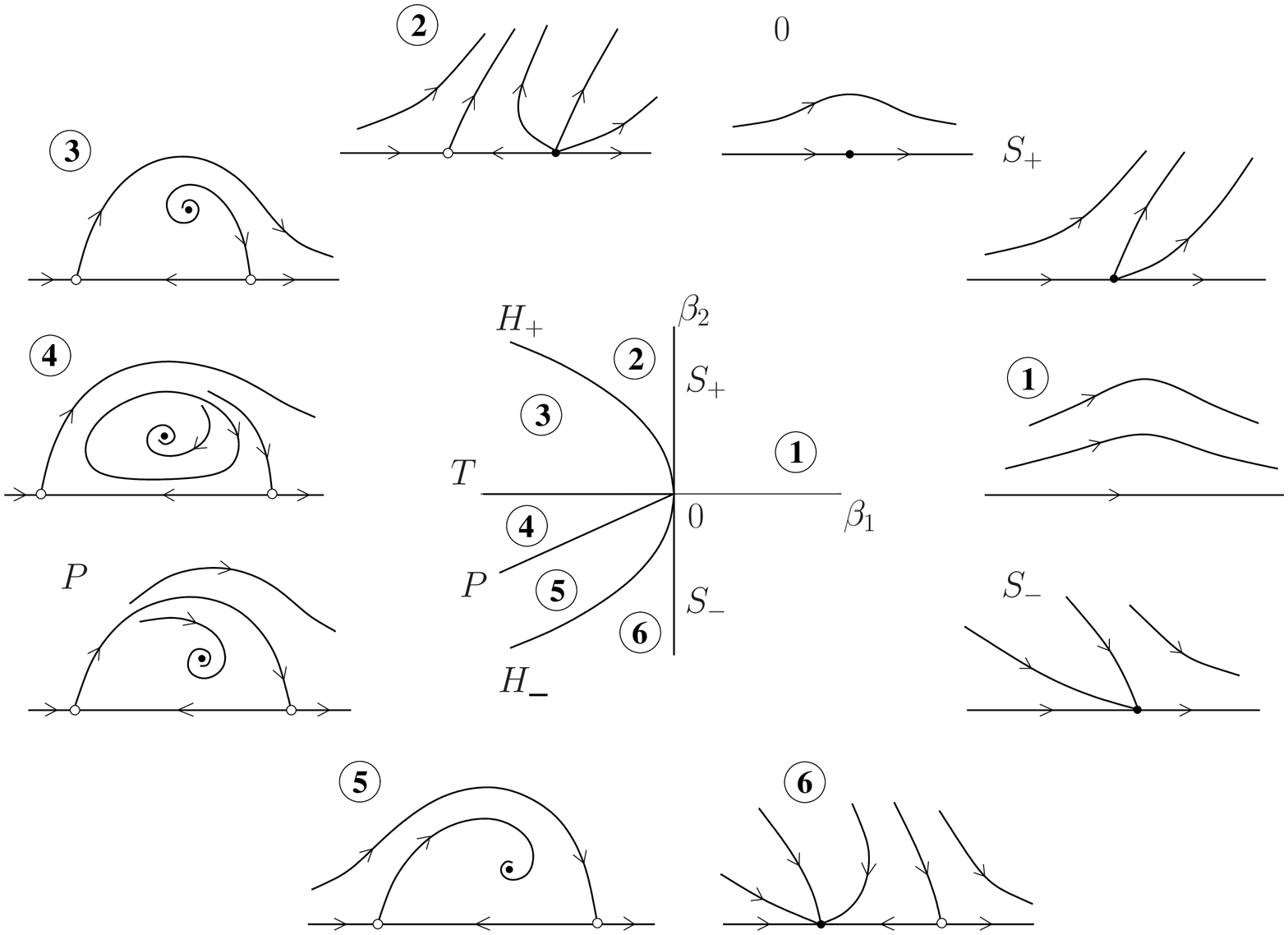}}\subfigure[$s=-1, \theta>0$]{
 \includegraphics[width=.48\textwidth]{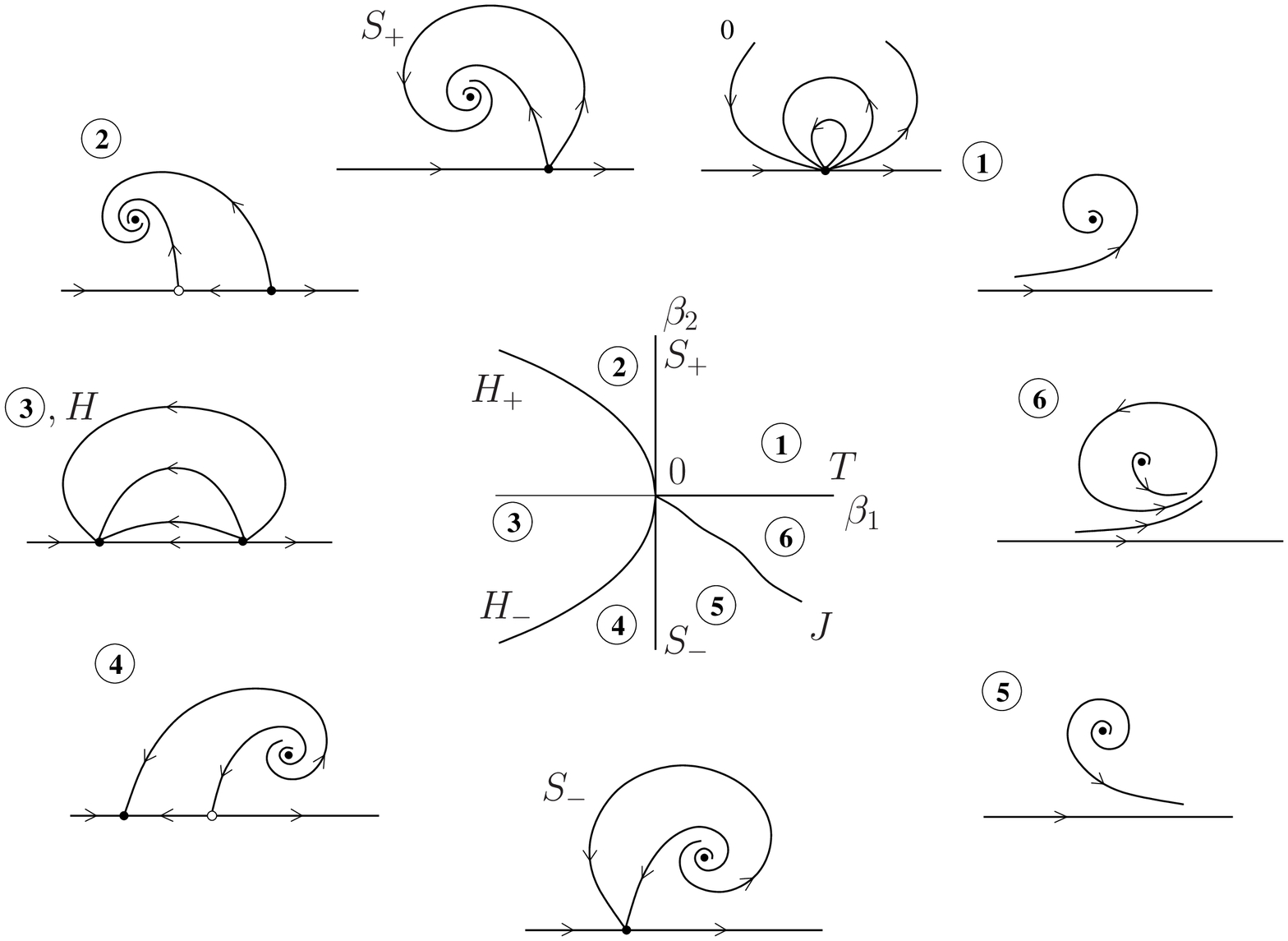}}
 \caption{Bifurcation diagrams of the truncated amplitude system (\ref{lpnsamplitude}) for the LPNS bifurcation.} 
\label{fig:NF_LPNS}
\end{figure}

Generically, a two-parameter unfolding of (\ref{ODE}) near this bifurcation restricted to the center manifold
is smoothly orbitally equivalent (with possible time reversal) to a system in which the equations for 
the transverse coordinates have the form
\begin{equation}
\label{lpnsnormalform}
\left\{\begin{array}{rcl}
 {\displaystyle \frac{d\xi}{d\tau}} &=& \beta_1 + \xi^2+ s \left|\zeta\right|^2 +
 O(\|(\xi,\zeta,\bar{\zeta})\|^4), \vspace{0.2cm}\\ 
 {\displaystyle \frac{d\zeta}{d\tau}} &=&  (\beta_2+i\omega_1)\zeta + 
 (\theta +i\vartheta)\xi \zeta+\xi^2\zeta + O(\|(\xi,\zeta,\bar{\zeta})\|^4),
 \end{array}\right.
\end{equation}
where the $O$-terms are still $T$-periodic in $\tau$. This system is similar to the
normal form for the Zero-Hopf bifurcation of equilibria (cf. Theorem 8.6 on page 338 in \cite{Ku:2004}). 
In Figure \ref{fig:NF_LPNS} the four possible bifurcation diagrams of the amplitude system 
for (\ref{lpnsnormalform}) without the $O$-terms,
\begin{equation}
\label{lpnsamplitude}
\left\{\begin{array}{rcl}
\dot{\xi}&=& \beta_1 + \xi^2 + s \rho^2,\\
\dot{\rho} &=& \rho (\beta_2 + \theta \xi + \xi^2),
\end{array}
\right.
\end{equation}
are reported depending on the sign of the normal form coefficients $s$ and $\theta$ \cite{Ku:2004}. 
Here and in what follows a dot means the derivative w.r.t. $\tau$.
\par
Let us now discuss the interpretation of the phase portraits in the $(\xi,\rho)$-plane of the truncated 
amplitude system in the context of the bifurcating limit cycle. The fixed points or limit cycles have additional dimensions from 
the phases of the periodic orbit itself plus the phases ignored in the reduction to the amplitude 
system. We note that in the amplitude system the vertical direction always corresponds to a 
Neimark-Sacker bifurcation, but that the horizontal component of the phase space has a different 
meaning. For LPNS, equilibria on the horizontal axis correspond to limit cycles. Equilibria
off the horizontal axis correspond to invariant 2D tori $\mathbb{T}^{2}$ and the periodic 
orbit which exists if $s\theta<0$ corresponds to an invariant 3D torus $\mathbb{T}^{3}$. 
\par
The critical values of $s$ and $\theta$ can be expressed in terms of the coeffcients of 
(\ref{eq:NF-LPNS}) as
$$
s=\mbox{sign }(a_{200}a_{011}),~~\theta = \frac{\Re(b_{110})}{a_{200}}.
$$ 
These values determine the bifurcation scenario. For $s\theta<0$, a $3$-torus appears in the unfolding via a Neimark-Sacker
bifurcation. The stability of 
this torus is determined by the third order terms in (\ref{eq:NF-LPNS}). Indeed, the sign of the 
corresponding first Lyapunov coefficient for the Hopf bifurcation in (\ref{lpnsamplitude}) is opposite to 
that of $\theta$ but the `time' in  (\ref{lpnsnormalform}) is rescaled with factor
$$
E = \Re\left(b_{210}+b_{110}\left( \frac{\Re(b_{021})}{a_{011}} - \frac{3a_{300}}{2a_{200}}+ 
\frac{a_{111}}{2a_{011}}\right)- \frac{b_{021}a_{200}}{a_{011}}\right).
$$
(see page 337 in \cite{Ku:2004}). If $E \cdot l_1<0$, an unstable $3$-torus appears, 
if $E \cdot l_1>0$, the $3$-torus is stable. The output given by {\sc matcont} is $(s,\theta,E)$ 
\footnote{Remark that $E=NaN$ is reported when terms up to only the second order are computed.}.

Note that Figure \ref{fig:NF_LPNS} presents  bifurcations of 
the truncated system (\ref{lpnsnormalform}) that only {\em approximates} the full normalized unfolding. 
In particular, the orbit structure on the invariant tori can differ from that
for the approximating system due to phase locking. Moreover,  the destruction of $\mathbb{T}^{3}$ via a 
heteroclinic bifurcation in case (c) of Figure \ref{fig:NF_LPNS} becomes a complicated sequence of global bifurcations involving 
stable and unstable invariant sets of cycles and tori. All these bifurcations, however, occur in the
exponentially-small parameter wedge near the heteroclinic bifurcation curve $P$. For detailed discussions
of the effects of the truncation, also in the two other cases, we refer to \cite{KuMe:2006,BSV:2010} and references therein.

\subsubsection{PDNS} \label{Appendix:2PDNS}
Generically, a two-parameter unfolding of (\ref{ODE}) near this bifurcation restricted to the center manifold
is smoothly orbitally equivalent to a system in which the equations for 
the transverse coordinates have the form
\begin{equation}
\left\{ 
\begin{array}{rcl}
\dot v_1 &=&\mu_1v_1+P_{11}v_1^3+P_{12}v_1\left|v_2\right|^2+S_1v_1\left|v_2\right|^4 + O(\|(v_1,v_2,\bar{v}_2)\|^6), \smallskip \\ 
\dot v_2 &=&(\mu_2+i\omega_2)v_2+P_{21}v_1^2v_2+P_{22}v_2\left|v_2\right|^2+S_2v_1^4v_2+
iR_2v_2\left|v_2\right|^4 + O(\|(v_1,v_2,\bar{v}_2)\|^6),
\end{array}
\right.
\label{pdnsnormalform}
\end{equation}
where the $O$-terms are still $T$-periodic in $\tau$. This system is similar to one of the normal forms
for the Hopf-Hopf bifurcations of equilibria (cf. Lemma 8.14 on page 354 in \cite{Ku:2004}).

The amplitude system for (\ref{pdnsnormalform}) without the $O$-terms is
\begin{equation}
\label{pdnsamplitude}
\left\{\begin{array}{rcl}
\dot{r}_1 &=& r_1(\mu_1 + p_{11}r_1^2 + p_{12}r_2^2 + s_1 r_2^4),\\
\dot{r}_2 &=& r_1(\mu_2 + p_{21}r_1^2 + p_{22}r_2^2 + s_2 r_1^4),
\end{array}
\right.
\end{equation}
where
$$
p_{11} = P_{11},~  p_{12} = P_{12},~  p_{21} = \Re(P_{21}),~  p_{22} = \Re(P_{22}),~  s_1 = S_{1},~  s_2 = \Re(S_{2}).
$$
The values of $p_{jk}$ and $s_j$, for $j,k=1,2$, and the quantities
\begin{eqnarray*}
\theta = \frac{p_{12}}{p_{22}}, \delta = \frac{p_{21}}{p_{11}}, 
\Theta = \frac{s_1}{p_{22}^2}, \Delta = \frac{s_2}{p_{11}^2}
\label{bifanalysispdns1}
\end{eqnarray*}
indicate in which bifurcation scenario we are (see Section 8.6.2 in \cite{Ku:2004}). 
\par
\begin{figure}[htb]
\centering \subfigure[]{
 \includegraphics[width=.44\textwidth]{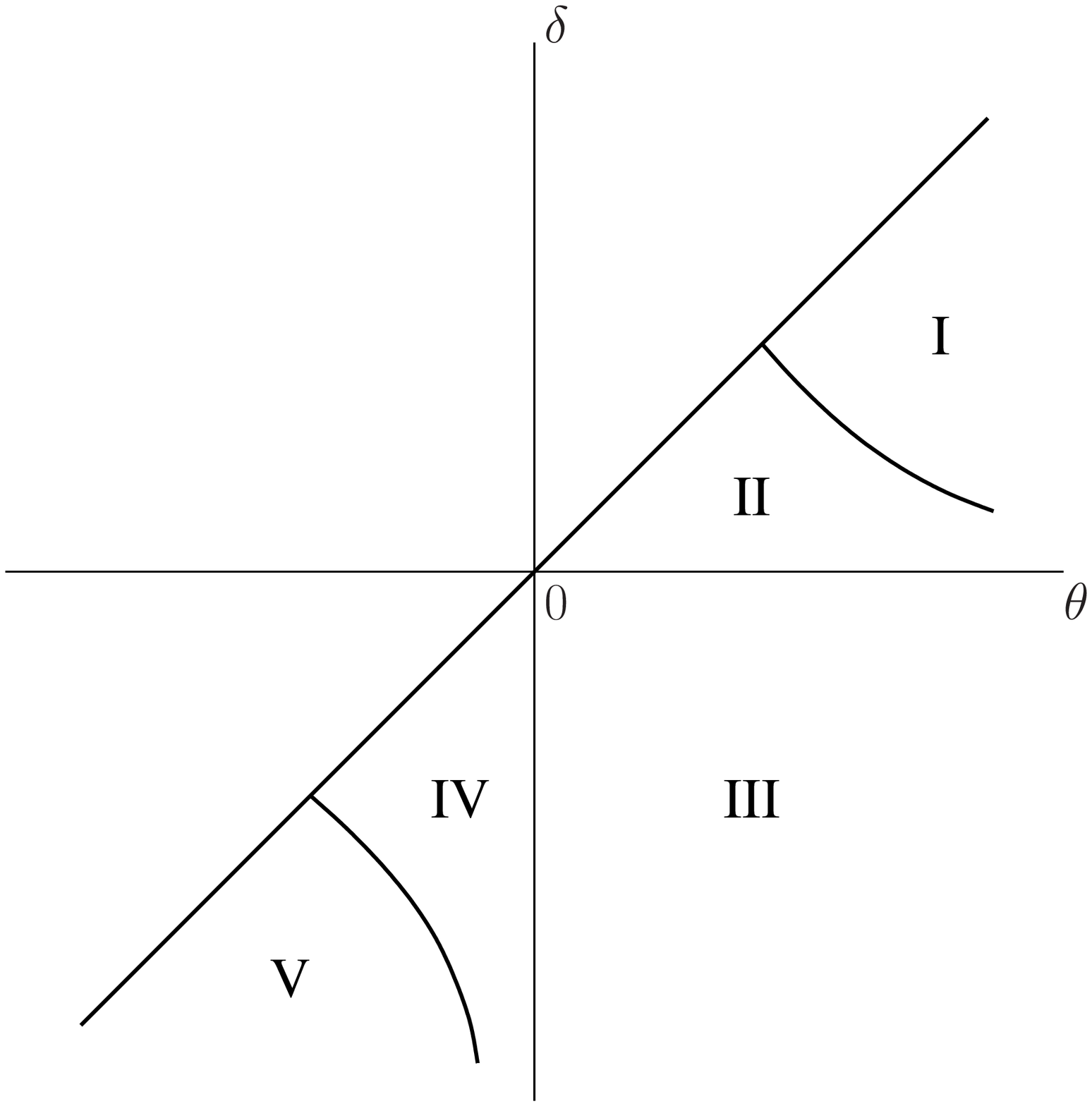}
 \label{fig:NF_PDNSoverviewleft}}
 \subfigure[]{
 \includegraphics[width=.44\textwidth]{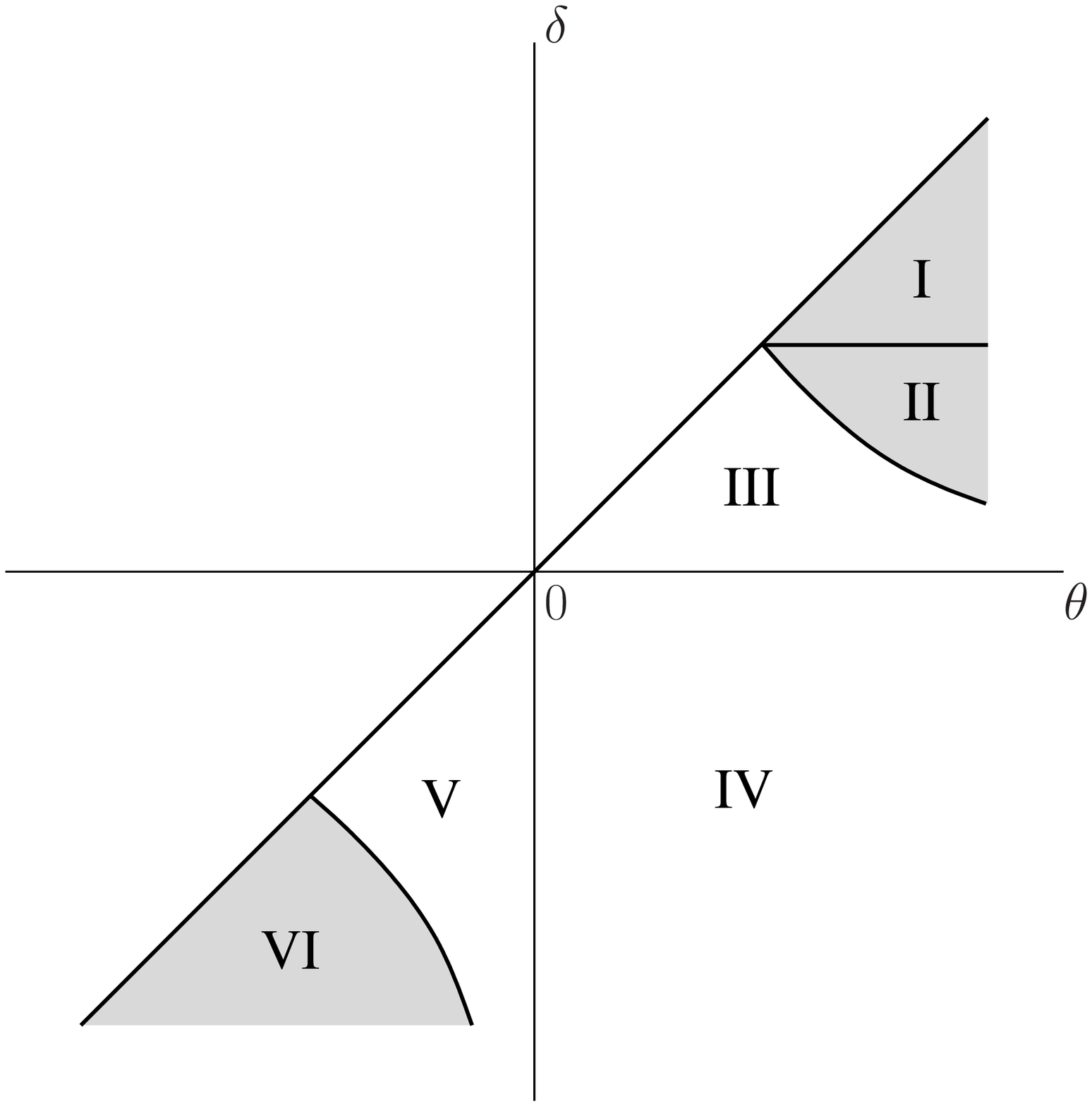}}
 \caption{(a) the five subregions in the $(\theta,\delta)$-plane in the ``simple" case; 
 (b) the six subregions in the $(\theta,\delta)$-plane in the ``difficult" case.} 
 \label{fig:NF_PDNSoverview}
\end{figure}
\begin{figure}[htb]
\centering \subfigure[]{
 \includegraphics[width=.44\textwidth]{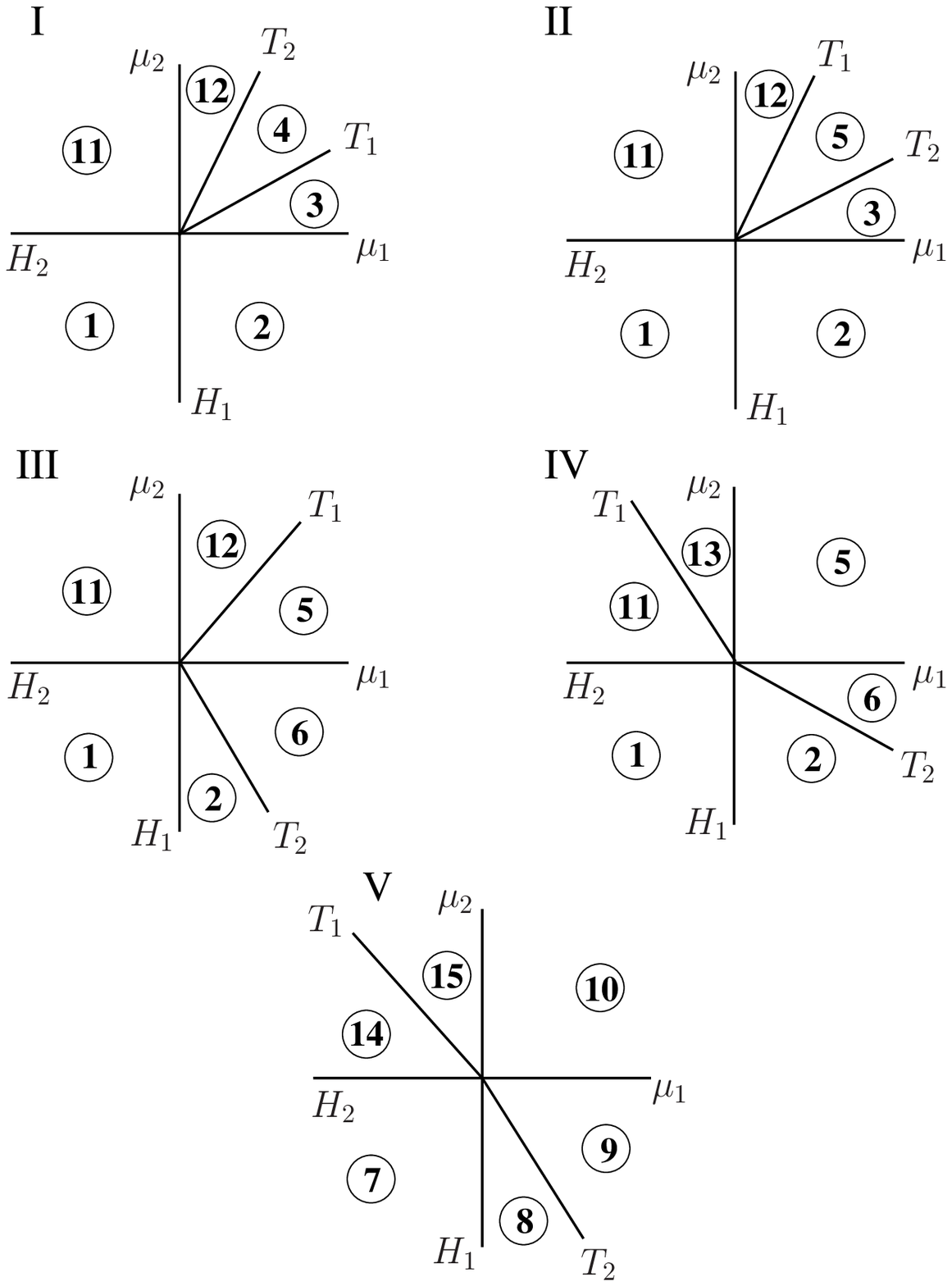}
 \label{fig:NF_PDNSsimple_1}} 
 \subfigure[]{
 \includegraphics[width=.44\textwidth]{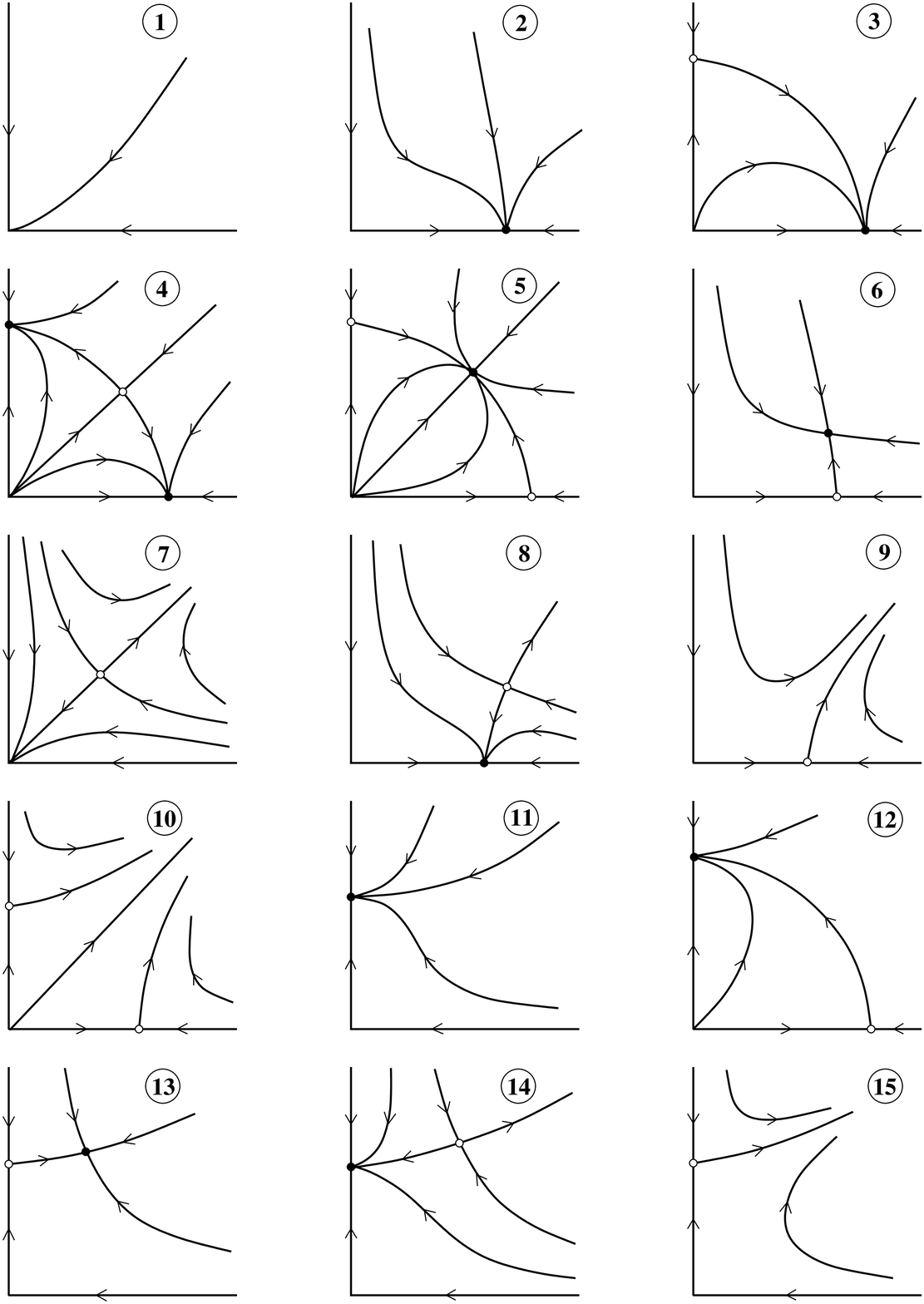}
 \label{fig:NF_PDNSsimple_2}}
 \caption{Bifurcation diagrams of the amplitude system (\ref{pdnsamplitude}) for the PDNS and NSNS bifurcations: 
 (a) parametric portraits in the ``simple" case; (b) phase portraits in the ``simple" case.} 
\label{fig:NF_PDNSsimple}
\end{figure}

In the ``simple" case where $p_{11}p_{22}>0$, there are five topologically different bifurcation diagrams
of the truncated amplitude system (\ref{pdnsamplitude}), corresponding to the following cases: 
\begin{itemize}
	\item[I.] $\theta>0, \delta>0, \theta\delta>1$
	\item[II.] $\theta>0, \delta>0, \theta\delta<1$
	\item[III.] $\theta>0, \delta<0$
	\item[IV.] $\theta<0, \delta<0, \theta\delta<1$
	\item[V.] $\theta<0, \delta<0, \theta\delta>1$	
\end{itemize}
If $\delta >\theta$, reverse the role of $\theta$ and $\delta$. Each case corresponds with a 
region in the $(\theta,\delta)$-plane, see Figure \ref{fig:NF_PDNSoverview} (a). The 
parametric portraits belonging to the different regions can be seen in Figure \ref{fig:NF_PDNSsimple} (a), 
with corresponding phase portraits in the $(r_1,r_2)$-plane in Figure \ref{fig:NF_PDNSsimple} (b). 
The phase portraits are only shown for the case when $p_{11}<0$ and $p_{22}<0$. The case 
$p_{11}>0$ and $p_{22}>0$ can be reduced to the considered one by reversing time.

\begin{figure}[htb]
\centering \subfigure[]{
 \includegraphics[width=.48\textwidth]{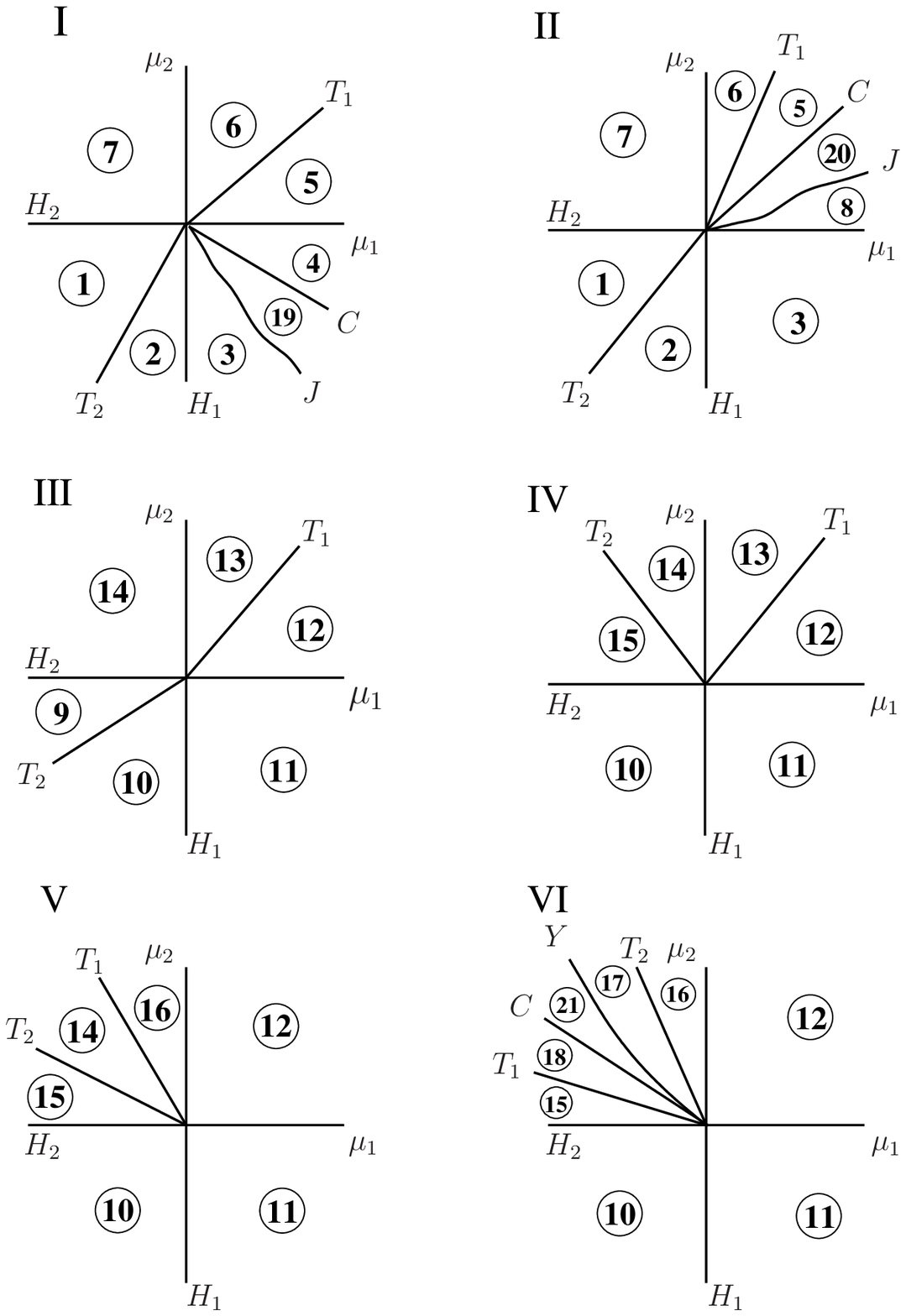}}
 \subfigure[]{
 \includegraphics[width=.48\textwidth]{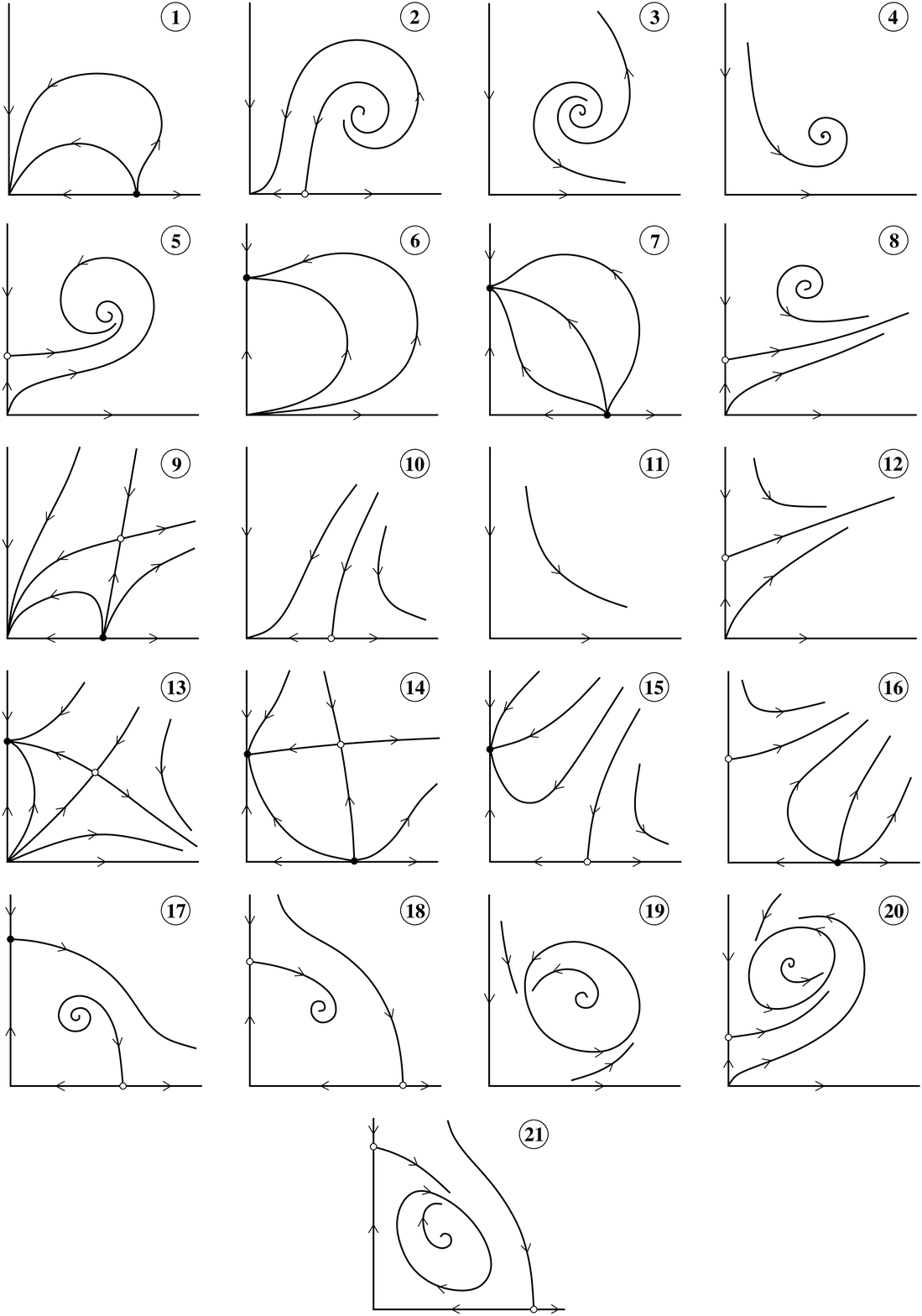}}
 \caption{Bifurcation diagrams of the amplitude system (\ref{pdnsamplitude}) for the PDNS and NSNS bifurcations: 
 (a) parametric portraits in the ``difficult" case; 
 (b) phase portraits in the ``difficult" case.} 
 \label{fig:NF_PDNSdifficult}
\end{figure}

In the ``difficult" case where $p_{11}p_{22}<0$ however, there are six essentially different bifurcation 
diagrams: 
\begin{itemize}
	\item[I.] $\theta>1, \delta>1$
	\item[II.] $\theta>1, \delta<1, \theta\delta>1$
	\item[III.] $\theta>0, \delta>0, \theta\delta<1$
	\item[IV.] $\theta>0, \delta<0$
	\item[V.] $\theta<0, \delta<0, \theta\delta<1$	
	\item[VI.] $\theta<0, \delta<0, \theta\delta>1$	
\end{itemize}
The regions in the $(\theta,\delta)$-plane are shown in Figure \ref{fig:NF_PDNSoverview} (b). 
The related parametric portraits and phase portraits of (\ref{pdnsamplitude}) 
are given in Figure \ref{fig:NF_PDNSdifficult}.
Only the case $p_{11}>0$ and $p_{22}<0$ is presented, to which the opposite one can be easily reduced.

We note that Section 8.6.2 in \cite{Ku:2004} for the ``difficult" case contains
a few errors in the figures and in the asymptotic expression for the heteroclinic 
bifurcation curve\footnote{Unfortunately, there is also a minor misprint in our earlier ``correction" 
for the heteroclinic curve given in \cite{KuMe:2006}.}. Therefore, for completeness, we provide the figures
and correct asymptotics in Appendix \ref{Appendix:6}.

The critical values of $P_{jk}$ and $S_j$ can be expressed in terms of the coeffcients of 
(\ref{eq:NF-PDNS}) as
$$
P_{11}= a_{2100},~ P_{12}=a_{1011},~ \Re(P_{21})= \Re(b_{1110}),~ \Re(P_{22}) = \Re(b_{0021}),
$$
and
\begin{eqnarray*}
S_1 &=& a_{1022} + a_{1011}\left( \frac{\Re(b_{1121})}{\Re(b_{1110})}-2 \frac{\Re(b_{0032})}{\Re(b_{0021})} - 
\frac{a_{3200}\Re(b_{0021})}{a_{2100}\Re(b_{1110})} \right),\\
\Re(S_2) &=&\Re(b_{2210}) + \Re(b_{1110})\left( \frac{a_{2111}}{a_{1011}}-2 \frac{a_{3200}}{a_{2100}} - 
\frac{a_{2100}\Re(b_{0032})}{a_{1011}\Re(b_{0021})} \right).
\end{eqnarray*}
(see page 356 in \cite{Ku:2004}). 

The fifth-order terms in (\ref{eq:NF-PDNS}) determine the stability of 
the tori in the ``difficult" cases. In fact, the sign of the first Lyapunov coefficient for 
the Neimark-Sacker bifurcation is given by
\begin{equation}
\label{l_1}
\mbox{sign } l_1 = -\mbox{sign}\left(\theta(\theta(\theta-1)\Delta+\delta(\delta-1)\Theta) \right).
\end{equation}
The output of {\sc matcont} is $(p_{11},p_{22},\theta,\delta,\mbox{sign } l_1)$\footnote{Remark 
that $\mbox{sign } l_1=NaN$ is reported when terms up to only the third order are computed.}.

For PDNS we have an interpretation analogous to LPNS, but the invariant sets may be ``doubled''. 
The origin always corresponds the original limit cycle. Other fixed points on the horizontal axis represent 
the period-doubled limit 
cycles, while a fixed point on the vertical axis corresponds to a $\mathbb{T}^{2}$. 
Fixed points off the coordinate 
axes correspond to doubled tori $\mathbb{T}^{2}$ and periodic orbits correspond to $\mathbb{T}^{3}$.
As in the LPNS case, Figures \ref{fig:NF_PDNSsimple} and \ref{fig:NF_PDNSdifficult} present  
bifurcations of the truncated amplitude system that only {\em approximates} the 
full normalized unfolding. In particular, one has to be carefull with `torus doubling', which is in fact a complicated
quasiperiodic bifurcation \cite{Los:1989,ViBrSi:2011}.

\subsubsection{NSNS}\label{Appendix:2NSNS}
Generically, a two-parameter unfolding of (\ref{ODE}) near this bifurcation restricted to the center manifold
is smoothly orbitally equivalent to a system in which the equations for 
the transverse coordinates have the form
\begin{equation} 
\begin{cases}
\dot v_1 =(\mu_1+i\omega_1)v_1+P_{11}v_1\left|v_1\right|^2+P_{12}v_1\left|v_2\right|^2
+iR_1v_1\left|v_1\right|^4+S_1v_1\left|v_2\right|^4 + O(\|(v,\bar{v})\|^6), \\ 
\dot v_2 =(\mu_2+i\omega_2)v_2+P_{21}v_2\left|v_1\right|^2+P_{22}v_2\left|v_2\right|^2+
S_2v_2\left|v_1\right|^4+iR_2v_2\left|v_2\right|^4 + O(\|(v,\bar{v})\|^6),
\end{cases}
\label{nsnsnormalform}
\end{equation}
where the $O$-terms are $T$-periodic in $\tau$. Neglecting this periodicity, system (\ref{nsnsnormalform}) is the
normal form for the Hopf-Hopf bifurcation of equilibria (cf. Lemma 8.14 on page 354 in \cite{Ku:2004}).

The truncated amplitude system for (\ref{nsnsnormalform}) is the same as (\ref{pdnsamplitude}), where now
\begin{eqnarray*}
&&p_{11}=\Re(P_{11})= \Re(a_{2100}), ~p_{12}=\Re(P_{12})= \Re(a_{1011}),\\
&&p_{21}=\Re(P_{21})= \Re(b_{1110}), ~p_{22}=\Re(P_{22}) = \Re(b_{0021}),
\end{eqnarray*}
and 
\begin{eqnarray*}
&& s_1= \Re(S_1) = \Re(a_{1022}) + 
\Re(a_{1011})\left( \frac{\Re(b_{1121})}{\Re(b_{1110})}-2 
\frac{\Re(b_{0032})}{\Re(b_{0021})} - \frac{\Re(a_{3200})\Re(b_{0021})}{\Re(a_{2100})\Re(b_{1110})} \right),\\
&& s_2 = \Re(S_2) =\Re(b_{2210}) + 
\Re(b_{1110})\left( \frac{\Re(a_{2111})}{\Re(a_{1011})}-2 
\frac{\Re(a_{3200})}{\Re(a_{2100})} - \frac{\Re(a_{2100})\Re(b_{0032})}{\Re(a_{1011})\Re(b_{0021})} \right).
\end{eqnarray*}
The output of {\sc matcont} is $(p_{11},p_{22},\theta,\delta,\mbox{sign } l_1)$\footnote{Remark 
that $\mbox{sign } l_1=NaN$ is reported when terms up to only the third order are computed.}.

Although the phase portraits of the truncated amplitude system are the same as for PDNS, 
their interpretation is slightly different, since they `live' in the $(|v_1|,|v_2|)$-plane. Here, on both axes the fixed points 
correspond to invariant 2D tori $\mathbb{T}^{2}$ 
for the original system. Fixed points off the coordinate axes and limit cycles correspond to $\mathbb{T}^{3}$ 
and $\mathbb{T}^{4}$, respectively. The usual remark on the approximate nature of the bifurcation diagrams
applies here as well.

\section{Computation of critical coefficients}
\label{Section:Theory}
As was mentioned in the previous section, the stability of the extra torus appearing in the ``difficult" cases is determined by third 
order terms for the LPNS bifurcation and fifth order terms for the PDNS and NSNS bifurcations. 
In the ``simple" cases, second order derivatives are sufficient to determine the behaviour in the LPNS 
bifurcations and third order derivatives are sufficient in the PDNS and NSNS bifurcations. 
Therefore, we restrict our computations in this section to second order terms in the LPNS case and up 
to and including third order terms in the PDNS and NSNS cases. The expressions of the third order 
coefficients for { LPNS} and fourth and fifth order coefficients for PDNS and NSNS are 
given in Appendix \ref{Appendix:5}. Remark that for efficiency reasons these higher order coefficients 
are not computed in {\sc matcont}, unless explicitly requested by the user.

\subsection{LPNS}\label{Section:LPNS}

The four-dimensional critical center manifold $W^c(\Gamma)$ at
the LPNS bifurcation can be parametrized locally by
$(\xi_1,\xi_2,\tau) \in \R \times \C \times [0,T]$ as 
\begin{equation}
u=u_0(\tau)+\xi_1 v_1(\tau) + \xi_2 v_2(\tau) +\bar \xi_2 \bar v_2(\tau)+ H(\xi_1,\xi_2,\tau), 
\label{eq:CM_LPNS}
\end{equation}
where $H$ satisfies $H(\xi_1,\xi_2,T)=H(\xi_1,\xi_2,0)$ and has the Taylor expansion 
\begin{eqnarray}
 H(\xi_1,\xi_2,\tau)= \sum_{\substack{2\leq i+j+k\leq 3}} 
 \frac{1}{i!j!k!}h_{ijk}(\tau)\xi_1^i\xi_2^j\bar{\xi}_2^k + O(\|\xi\|^4) , \label{eq:H_LPNS}
\end{eqnarray}
where the eigenfunctions $v_1$ and $v_2$ are defined by 
\begin{eqnarray} \label{eq:EigenFunc_LPNS}
&&\left\{\begin{array}{rcl}
\dot{v}_1-A(\tau)v_1 - F(u_0) & = & 0,\ \tau \in [0,T], \\
v_1(T)-v_1(0) & = & 0,\\
\int_{0}^{T} {\langle v_1,F(u_0)\rangle d\tau} & = & 0,\\
\end{array}
\right. 
\end{eqnarray}
and
\begin{eqnarray}\label{eq:EigenFunc_LPNS_2}
&&\left\{\begin{array}{rcl}
\dot{v}_2-A(\tau)v_2+i\omega v_2 & = & 0,\ \tau \in [0,T], \\
v_2(T) - v_2(0) & = & 0,\\
\int_{0}^{T} {\langle v_2,v_2\rangle d\tau} - 1  & = & 0.\\
\end{array}
\right.
\end{eqnarray}

The functions $v_1$ and $v_2$ exist because of Lemma~2 of \cite{Io:88}. 
The functions $h_{ijk}$ will be found by solving appropriate
BVPs, assuming that \eqref{eq:P.1} restricted to $W^c(\Gamma)$ has
the normal form \eqref{eq:NF-LPNS}. 

The coefficients of the normal
form  arise from the solvability conditions for the BVPs as
integrals of scalar products over the interval $[0,T]$.
Specifically, those scalar products involve among other things the quadratic and
cubic terms of \eqref{eq:MULT} near the periodic solution $u_0$,
the generalized eigenfunction $v_1$ and eigenfunction $v_2$, and the
adjoint eigenfunctions $\varphi^*$, $v_1^*$ and $v_2^*$ as solutions
of the problems
\begin{eqnarray}
&&\left\{\begin{array}{rcl}
 \dot{\varphi}^*+A^{\rm T}(\tau)\varphi^* & = & 0,\ \tau \in [0,T], \\
 \varphi^*(T)-\varphi^*(0) & = & 0, \label{eq:AdjEigenFunc_LPNS}\\
 \int_{0}^{T} {\langle \varphi^*,v_1 \rangle d\tau} -1 & = & 0,
 \end{array} \right.
 \end{eqnarray}
 \begin{eqnarray}
 &&\left\{\begin{array}{rcl}
 \dot{v}_1^*+A^{\rm T}(\tau) v_1^*+\varphi^* & = & 0,\ \tau \in [0,T], \\
 v_1^*(T)-v_1^*(0) & = & 0, \label{eq:AdjEigenFunc_LPNS_2}\\
 \int_{0}^{T} {\langle v_1^*,v_1 \rangle d\tau} & = & 0,
 \end{array} \right. 
 \end{eqnarray}
and
 \begin{eqnarray}
 &&\left\{\begin{array}{rcl}
 \dot{v}_2^*+A^{\rm T}(\tau) v_2^*+i\omega v_2^* & = & 0,\ \tau \in [0,T], \\
 v_2^*(T)-v_2^*(0) & = & 0, \label{eq:AdjEigenFunc_LPNS_3}\\
 \int_{0}^{T} {\langle v_2^*,v_2\rangle d\tau} -1 & = & 0.
 \end{array} \right.
\end{eqnarray}

In what follows we will make use of the orthogonality condition
\begin{eqnarray}
\int_{0}^{T} {\langle \varphi^*,F(u_0) \rangle d\tau} = 0,
\label{orthLPNS}
\end{eqnarray}
and the normalization condition
\begin{eqnarray}
\int_{0}^{T} {\langle v_1^*,F(u_0) \rangle d\tau} = 1,
\label{normLPNS}
\end{eqnarray}
which can be easily obtained from \eqref{eq:EigenFunc_LPNS}, \eqref{eq:AdjEigenFunc_LPNS} 
and \eqref{eq:AdjEigenFunc_LPNS_2}.

To derive the normal form coefficients we write down the
homological equation and compare term by term. We therefore substitute (\ref{eq:CM_LPNS}) 
into (\ref{eq:P.1}), using (\ref{eq:MULT}), (\ref{eq:NF-LPNS}) and (\ref{eq:H_LPNS}). 
By collecting the constant and linear terms we get the identities
$$\dot{u}_0=F(u_0), \qquad \dot v_1-F(u_0)=A(\tau) v_1, \qquad \dot{v}_2+i\omega v_2=A(\tau) v_2,$$
and the complex conjugate of the last equation.

By collecting the $\xi_1^2$-terms we find an equation for $h_{200}$
\begin{equation}\label{eq:h200_LPNS}
 \dot h_{200}-A(\tau) h_{200}=B(\tau;v_1,v_1)-2 a_{200} v_1-2 \alpha_{200} \dot u_0+2 \dot v_1,
\end{equation}
to be solved in the space of functions satisfying
$h_{200}(T)=h_{200}(0)$. In this space, the differential operator
$\frac{d}{d\tau}-A(\tau)$ is singular and its null-space is
spanned by $\dot u_0$. The Fredholm solvability
condition
\[
 \int_0^{T} \langle \varphi^*,B(\tau;v_1,v_1)-2 a_{200} v_1-2 \alpha_{200} \dot u_0+
 2 \dot v_1 \rangle\; d \tau =0
\]
allows one to calculate the coefficient $a_{200}$ in (\ref{eq:NF-LPNS}) due to the required normalization 
in \eqref{eq:AdjEigenFunc_LPNS}, i.e.
\begin{equation} \label{eq:a200_LPNS}
 \boxed{a_{200}=\frac{1}{2}\int_0^{T} \langle \varphi^*,B(\tau;v_1,v_1)+2 A(\tau) v_1 \rangle\; d \tau,}
\end{equation}
taking \eqref{eq:EigenFunc_LPNS} and \eqref{orthLPNS} into account. With $a_{200}$ defined in this way, 
let $h_{200}$ be a solution of (\ref{eq:h200_LPNS}) in the space of functions satisfying $h_{200}(0)=h_{200}(T)$. Notice that if
$h_{200}$ is a solution of \eqref{eq:h200_LPNS}, then also $h_{200}+\varepsilon_1 F(u_0)$ satisfies
\eqref{eq:h200_LPNS}, since $F(u_0)$ is in the kernel of
the operator $\frac{d}{d\tau}-A(\tau)$. In order to obtain a unique solution (without a component along  the null
eigenspace) we impose the following orthogonality condition which determines the value of $\varepsilon_1$
$$ \int_0^{T} \langle v_1^*,h_{200}\rangle \; d\tau =0,$$
since \eqref{normLPNS} holds. Thus $h_{200}$ is the unique solution of the BVP
\begin{equation}  \label{eq:BVP_h200_LPNS}
\left\{\begin{array}{rcl}
\dot{h}_{200}-A(\tau) h_{200} - B(\tau;v_1,v_1) -2A(\tau)v_1+2 a_{200} v_1+2 \alpha_{200} 
\dot u_0-2\dot u_0 & = & 0,\ \tau \in [0,T], \\
h_{200}(T)- h_{200}(0) & = & 0,\\
\int_0^{T} \langle v_1^*,h_{200}\rangle \; d\tau & = & 0.
\end{array}
\right.
\end{equation}

By collecting the $\xi_2^2$-terms (or $\bar \xi_2^2$-terms) we find an equation for $h_{020}$
\begin{equation*}\label{eq:h020_LPNS}
 \dot h_{020}-A(\tau) h_{020}+2i\omega h_{020}=B(\tau;v_2,v_2),
\end{equation*}
(or its complex conjugate). This equation has a unique solution $h_{020}$ satisfying
$h_{020}(T)=h_{020}(0)$, since due to the spectral assumptions $e^{2i\omega T}$ is not a 
multiplier of the critical cycle. Thus, $h_{020}$ can be found by solving
\begin{equation}  \label{eq:BVP_h020_LPNS}
\left\{\begin{array}{rcl}
\dot h_{020}-A(\tau) h_{020}+2i\omega h_{020}-B(\tau;v_2,v_2) & = & 0,\ \tau \in [0,T], \\
h_{020}(T)- h_{020}(0) & = & 0.
\end{array}
\right.
\end{equation}

By collecting the $\xi_1\xi_2$-terms we obtain an equation for $h_{110}$
\begin{equation*}\label{eq:h110_LPNS}
\dot h_{110}-A(\tau) h_{110} + i\omega h_{110} = B(\tau;v_1,v_2)-b_{110}v_2 +\dot v_2+i\omega v_2,
\end{equation*}
to be solved in the space of functions satisfying $h_{110}(T)=h_{110}(0)$. In this space, the 
differential operator
$\frac{d}{d\tau}-A(\tau)+i\omega$ is singular, since $e^{i\omega T}$ is a critical multiplier. 
So we can impose the following Fredholm solvability condition
$$\int_0^{T} \langle v_2^*, B(\tau;v_1,v_2)-b_{110}v_2 +\dot v_2+i\omega v_2\rangle \; d\tau = 0,$$
which due to the normalization condition in \eqref{eq:AdjEigenFunc_LPNS_3} determines the value 
of the normal form coefficient $b_{110}$, yielding
\begin{equation} \label{eq:b110_LPNS}
 \boxed{b_{110}=\int_0^{T} \langle v_2^*,B(\tau;v_1,v_2)+ A(\tau) v_2 \rangle\; d \tau.}
\end{equation}
The nullspace belonging to the operator $\frac{d}{d\tau}-A(\tau)+i\omega$ is one-dimensional and 
spanned by $v_2$. To determine $h_{110}$ uniquely, we need to impose an orthogonality condition 
with a vector whose inproduct with $v_2$ is non-zero. $v_2^*$ can be choosen because of the 
normalisation condition in \eqref{eq:AdjEigenFunc_LPNS_3}. In fact, $h_{110}$ only appears in 
the normal form coefficient $b_{210}$ (see Appendix \ref{Appendix:5_lpns}), and a different normalization 
of $h_{110}$ does not influence the value of that normal form coefficient. Therefore, we obtain 
$h_{110}$ as the unique solution of the BVP
\begin{equation}  \label{eq:BVP_h110_LPNS}
\left\{\begin{array}{rcl}
\dot h_{110}-A(\tau) h_{110} + i\omega h_{110} - B(\tau;v_1,v_2)+b_{110}v_2 -A(\tau) v_2& = & 0,\ 
\tau \in [0,T], \\
h_{110}(T)- h_{110}(0) & = & 0,\\
\int_0^{T} \langle v_2^*,h_{110}\rangle \; d\tau & = & 0.
\end{array}
\right.
\end{equation}

By collecting the $ \left|\xi_2\right|^2$-terms we obtain a singular equation for $h_{011}$, namely
\begin{equation*}\label{eq:h011_LPNS}
\dot h_{011}-A(\tau) h_{011} = B(\tau;v_2,\bar v_2)-a_{011}v_1-\alpha_{011}\dot u_0,
\end{equation*}
to be solved in the space of functions satisfying $h_{011}(T)=h_{011}(0)$. The non-trivial kernel 
of the operator
$\frac{d}{d\tau}-A(\tau)$ is spanned by $\dot u_0$. So, the following Fredholm solvability condition 
is involved
$$\int_0^{T} \langle \varphi^*,B(\tau;v_2,\bar v_2)-a_{011}v_1-\alpha_{011}\dot u_0\rangle \; d\tau = 0,$$
which gives us the expression for the normal form coefficient $a_{011}$, i.e.
\begin{equation} \label{eq:a011_LPNS}
 \boxed{a_{011}=\int_0^{T} \langle \varphi^*,B(\tau;v_2,\bar v_2) \rangle\; d \tau.}
\end{equation}
We impose the orthogonality condition with the adjoint generalized eigenfunction $v_1^*$ to obtain $h_{011}$ 
as the unique solution of 
\begin{equation}  \label{eq:BVP_h011_LPNS}
\left\{\begin{array}{rcl}
\dot h_{011}-A(\tau) h_{011} - B(\tau;v_2,\bar v_2)+a_{011}v_1+\alpha_{011}\dot u_0 & = & 0,\ \tau \in [0,T], \\
h_{011}(T)- h_{011}(0) & = & 0,\\
\int_0^{T} \langle v_1^*,h_{011}\rangle \; d\tau & = & 0.
\end{array}
\right.
\end{equation}

We remark that the values of $\alpha_{200}$ and $\alpha_{011}$ are not determined by the homological equation. 
We therefore put them equal to zero.

Third order coefficients are only needed to determine the stability of the torus, if this torus exists. 
For completeness, we have listed these terms in Appendix \ref{Appendix:5}.

\subsection{PDNS} \label{Section:PDNS}
The four-dimensional critical center manifold $W^c(\Gamma)$ at
the PDNS bifurcation can be parametrized locally by
$(\xi_1,\xi_2,\tau) \in \R \times \C \times [0,2T]$ as 
\begin{equation}
u=u_0(\tau)+\xi_1 v_1(\tau) + \xi_2 v_2(\tau) +\bar \xi_2 \bar v_2(\tau)+ H(\xi_1,\xi_2,\tau), \label{eq:CM_PDNS}
\end{equation}
where $H$ satisfies $H(\xi_1,\xi_2,2T)=H(\xi_1,\xi_2,0)$ and has the Taylor expansion 
\begin{equation}
 H(\xi_1,\xi_2,\tau)= \sum_{2\leq i+j+k\leq 5} \frac{1}{i!j!k!}h_{ijk}(\tau)
 \xi_1^i\xi_2^j\bar{\xi}_2^k
+ O(\|\xi\|^6) , 
\label{eq:H_PDNS}
\end{equation}
while the eigenfunctions $v_1$ and $v_2$ are defined by 
\begin{eqnarray} \label{eq:EigenFunc_PDNS}
&&\left\{\begin{array}{rcl}
\dot{v}_1-A(\tau)v_1 & = & 0,\ \tau \in [0,T], \\
v_1(T)+v_1(0) & = & 0,\\
\int_{0}^{T} {\langle v_1,v_1\rangle d\tau} -1 & = & 0,\\
\end{array}
\right. 
\end{eqnarray}
with $v_1(\tau+T)=-v_1(\tau)$ for $\tau \in [0,T]$ and
\begin{eqnarray}\label{eq:EigenFunc_PDNS_2}
&&\left\{\begin{array}{rcl}
\dot{v}_2-A(\tau)v_2+i\omega v_2 & = & 0,\ \tau \in [0,T], \\
v_2(T) - v_2(0) & = & 0,\\
\int_{0}^{T} {\langle v_2,v_2\rangle d\tau} - 1  & = & 0.\\
\end{array}
\right.
\end{eqnarray}

The functions $v_1$ and $v_2$ exist because of Lemma~5 of \cite{Io:88}. 
The functions $h_{ijk}$ can be found by solving appropriate
BVPs, assuming that \eqref{eq:P.1} restricted to $W^c(\Gamma)$ has
the normal form \eqref{eq:NF-PDNS}. Moreover, $u(\tau,\xi_1,\xi_2,\bar\xi_2) = 
u(\tau+T,-\xi_1,\xi_2,\bar\xi_2)$ so that
\begin{equation}
h_{ijk}(\tau)= (-1)^ih_{ijk}(\tau+T),
\label{propertyhijk}
\end{equation}
for $\tau \in [0,T]$. Therefore, we can restrict our computations to the interval $[0,T]$ instead of $[0,2T]$.

The coefficients of the normal
form  arise from the solvability conditions for the BVPs as
integrals of scalar products over the interval $[0,T]$.
Specifically, those scalar products involve among other things the quadratic and
cubic terms of \eqref{eq:MULT} near the periodic solution $u_0$,
$v_1$, $v_2$, and the
adjoint eigenfunctions $\varphi^*$, $v_1^*$ and $v_2^*$ as solutions
of the problems
\begin{eqnarray}
&&\left\{\begin{array}{rcl}
 \dot{\varphi}^*+A^{\rm T}(\tau)\varphi^* & = & 0,\ \tau \in [0,T], \\
 \varphi^*(T)-\varphi^*(0) & = & 0, \label{eq:AdjEigenFunc_PDNS}\\
 \int_{0}^{T} {\langle \varphi^*,F(u_0) \rangle d\tau} -1 & = & 0,
 \end{array} \right.
 \end{eqnarray}
 \begin{eqnarray}
 &&\left\{\begin{array}{rcl}
 \dot{v}_1^*+A^{\rm T}(\tau) v_1^* & = & 0,\ \tau \in [0,T], \\
 v_1^*(T)+v_1^*(0) & = & 0, \label{eq:AdjEigenFunc_PDNS_2}\\
 \int_{0}^{T} {\langle v_1^*,v_1 \rangle d\tau} -1& = & 0,
 \end{array} \right. 
 \end{eqnarray}
and
 \begin{eqnarray}
 &&\left\{\begin{array}{rcl}
 \dot{v}_2^*+A^{\rm T}(\tau) v_2^*+i\omega v_2^* & = & 0,\ \tau \in [0,T], \\
 v_2^*(T)-v_2^*(0) & = & 0, \label{eq:AdjEigenFunc_PDNS_3}\\
 \int_{0}^{T} {\langle v_2^*,v_2\rangle d\tau} -1 & = & 0.
 \end{array} \right.
\end{eqnarray}

By collecting the constant and linear terms we get the identities
$$
 \dot{u}_0=F(u_0), \qquad \dot v_1=A(\tau) v_1, \qquad \dot{v}_2+i\omega v_2=A(\tau) v_2,
$$
and the complex conjugate of the last equation, which merely reflect the definition of $u_0$ and (\ref{eq:EigenFunc_PDNS}), (\ref{eq:EigenFunc_PDNS_2}).

By collecting the $\xi_1^2$-terms we find an equation for $h_{200}$
\begin{equation}\label{eq:h200_PDNS}
 \dot h_{200}-A(\tau) h_{200}=B(\tau;v_1,v_1)-2 \alpha_{200} \dot u_0,
\end{equation}
to be solved in the space of functions satisfying
$h_{200}(T)=h_{200}(0)$. In this space, the differential operator
$\frac{d}{d\tau}-A(\tau)$ is singular and its null-space is
spanned by $\dot u_0$. The Fredholm solvability
condition
\[
 \int_0^{T} \langle \varphi^*,B(\tau;v_1,v_1)-2 \alpha_{200} \dot u_0 \rangle\; d \tau =0
\]
gives us the possibility to calculate $\alpha_{200}$ in (\ref{eq:NF-LPNS}) by the required normalization in \eqref{eq:AdjEigenFunc_PDNS}, i.e.
\begin{equation} \label{eq:alpha200_PDNS}
 \alpha_{200}=\frac{1}{2}\int_0^{T} \langle \varphi^*,B(\tau;v_1,v_1) \rangle\; d \tau.
\end{equation}
As before, $h_{200}$ is determined up to the addition of a multiple of $\dot u_0$, since $h_{200}+\varepsilon_1 F(u_0)$ is a solution of \eqref{eq:h200_PDNS} for every value of $\varepsilon_1$. We fix the value of $h_{200}$ by demanding the orthogonality with the adjoint eigenfunction corresponding with multiplier $1$, i.e.
$$ \int_0^{T} \langle \varphi^*,h_{200}\rangle \; d\tau =0.$$
We obtain $h_{200}$ then as the unique solution of the BVP
\begin{equation}  \label{eq:BVP_h200_PDNS}
\left\{\begin{array}{rcl}
\dot{h}_{200}-A(\tau) h_{200} -B(\tau;v_1,v_1)+2 \alpha_{200} \dot u_0& = & 0,\ \tau \in [0,T], \\
h_{200}(T)- h_{200}(0) & = & 0,\\
\int_0^{T} \langle \varphi^*,h_{200}\rangle \; d\tau & = & 0.
\end{array}
\right.
\end{equation}

By collecting the $\xi_2^2$-terms (or $\bar \xi_2^2$-terms) we find the differential equation for $h_{020}$
\begin{equation*}\label{eq:h020_PDNS}
 \dot h_{020}-A(\tau) h_{020}+2i\omega h_{020}=B(\tau;v_2,v_2),
\end{equation*}
or its complex conjugate. Since $e^{2i\omega T}$ is not a critical multiplier, no Fredholm solvability condition has to be satisified. $h_{020}$ can thus simply be found by solving
\begin{equation}  \label{eq:BVP_h020_PDNS}
\left\{\begin{array}{rcl}
\dot h_{020}-A(\tau) h_{020}+2i\omega h_{020}-B(\tau;v_2,v_2) & = & 0,\ \tau \in [0,T], \\
h_{020}(T)- h_{020}(0) & = & 0.
\end{array}
\right.
\end{equation}

The equation found by comparing the $\xi_1\xi_2$-terms is given by
\begin{equation*}\label{eq:h110_PDNS}
\dot h_{110}-A(\tau) h_{110} + i\omega h_{110} = B(\tau;v_1,v_2).
\end{equation*}
From \eqref{propertyhijk} it follows that $h_{110}$ is anti-periodic. Now, since $-e^{i\omega T}$ is not a multiplier of the critical cycle, no solvability condition has to be satisfied. Therefore, we can immediately obtain $h_{110}$ from 
\begin{equation}  \label{eq:BVP_h110_PDNS}
\left\{\begin{array}{rcl}
\dot h_{110}-A(\tau) h_{110} + i\omega h_{110} - B(\tau;v_1,v_2)& = & 0,\ \tau \in [0,T], \\
h_{110}(T)+ h_{110}(0) & = & 0,\\
\end{array}
\right.
\end{equation}

The $\left|\xi_2\right|^2$-terms lead to a singular equation for $h_{011}$, namely
\begin{equation*}\label{eq:h011_PDNS}
\dot h_{011}-A(\tau) h_{011} = B(\tau;v_2,\bar v_2)-\alpha_{011}\dot u_0,
\end{equation*}
to be solved in the space of $T$-periodic functions. The non-trivial kernel of the operator
$\frac{d}{d\tau}-A(\tau)$ is spanned by $\dot u_0$. So, the Fredholm solvability condition with the corresponding $T$-periodic adjoint eigenfunction is involved, i.e.
$$\int_0^{T} \langle \varphi^*,B(\tau;v_2,\bar v_2)-\alpha_{011}\dot u_0\rangle \; d\tau = 0,$$
from which the expression for the normal form coefficient $\alpha_{011}$ can be derived
\begin{equation} \label{eq:alpha011_PDNS}
 \alpha_{011}=\int_0^{T} \langle \varphi^*,B(\tau;v_2,\bar v_2) \rangle\; d \tau.
\end{equation}
Now, we still need to uniquely determine the multiple of $F(u_0)$ which can be added to the function $h_{011}$, and will therefore impose the orthogonality condition with $\varphi^*$ to obtain $h_{011}$ as the unique solution of 
\begin{equation}  \label{eq:BVP_h011_PDNS}
\left\{\begin{array}{rcl}
\dot h_{011}-A(\tau) h_{011} - B(\tau;v_2,\bar v_2)+\alpha_{011}\dot u_0 & = & 0,\ \tau \in [0,T], \\
h_{011}(T)- h_{011}(0) & = & 0,\\
\int_0^{T} \langle \varphi^*,h_{011}\rangle \; d\tau & = & 0.
\end{array}
\right.
\end{equation}
We have now examined all order two terms, and continue with the order three terms. 

Collecting the $\xi_1^3$-terms gives an equation for $h_{300}$ and will give us the possibility to compute the normal form coefficient $a_{300}$ in \eqref{eq:NF-PDNS}. The differential equation
\begin{equation*}\label{eq:h300_PDNS}
 \dot h_{300}-A(\tau) h_{300}=C(\tau;v_1,v_1,v_1)+3B(\tau;v_1,h_{200})-6 \alpha_{200} \dot v_1-6a_{300}v_1
\end{equation*}
has to be solved in the space of functions satisfying $h_{300}(T)=-h_{300}(0)$. The non-trivial anti-periodic kernel of the operator
$\frac{d}{d\tau}-A(\tau)$ is spanned by $v_1$. So, the Fredholm solvability condition with the anti-periodic adjoint eigenfunction $v_1^*$ is involved, i.e.
$$\int_0^{T} \langle v_1^*,C(\tau;v_1,v_1,v_1)+3B(\tau;v_1,h_{200})-6 \alpha_{200} \dot v_1-6a_{300}v_1\rangle \; d\tau = 0$$
and thus
\begin{equation} \label{eq:a300_PDNS}
 \boxed{a_{300}= \frac{1}{6} \int_0^{T} \langle v_1^*,C(\tau;v_1,v_1,v_1)+3B(\tau;v_1,h_{200})-6 \alpha_{200} A(\tau)v_1 \rangle\; d \tau,}
\end{equation}
due to the normalization condition from \eqref{eq:AdjEigenFunc_PDNS_2}. 
The usual orthogonality condition with the adjoint eigenfunction $v_1^*$ is imposed to obtain $h_{300}$ as the unique solution of 
\begin{equation}  \label{eq:BVP_h300_PDNS}
\left\{\begin{array}{rcl}
\dot h_{300}-A(\tau) h_{300}-C(\tau;v_1,v_1,v_1)-3B(\tau;v_1,h_{200})\\
+6 \alpha_{200} A(\tau) v_1+6a_{300}v_1 & = & 0,\ \tau \in [0,T], \\
h_{300}(T)+ h_{300}(0) & = & 0,\\
\int_0^{T} \langle v_1^*,h_{300}\rangle \; d\tau & = & 0.
\end{array}
\right.
\end{equation}

The $\xi_2^3$ (or $\bar \xi_2^3$)-terms from the homological equation give an equation for $h_{030}$
\begin{equation*}\label{eq:h030_PDNS}
 \dot h_{030}-A(\tau) h_{030}+3i\omega h_{030}=C(\tau;v_2,v_2,v_2)+3B(\tau;v_2,h_{020}),
\end{equation*}
or its complex conjugate. This equation has a unique solution $h_{030}$ satisfying
$h_{030}(T)=h_{030}(0)$, since due to the spectral assumptions $e^{3i\omega T}$ is not a multiplier of the critical cycle. Thus, $h_{030}$ can be found by solving
\begin{equation}  \label{eq:BVP_h030_PDNS}
\left\{\begin{array}{rcl}
\dot h_{030}-A(\tau) h_{030}+3i\omega h_{030}-C(\tau;v_2,v_2,v_2)-3B(\tau;v_2,h_{020}) & = & 0,\ \tau \in [0,T], \\
h_{030}(T)- h_{030}(0) & = & 0.
\end{array}
\right.
\end{equation}

By collecting the $\xi_1^2\xi_2$-terms we find an equation for $h_{210}$
\begin{equation}\label{eq:h210_PDNS}
\begin{array}{rcl}
 \dot h_{210}-A(\tau) h_{210}+i\omega h_{210}&=&C(\tau;v_1,v_1,v_2)+B(\tau;v_2,h_{200})+2B(\tau;v_1,h_{110})\\
 &-&2 \alpha_{200} \dot v_2-2b_{210}v_2-2i\omega \alpha_{200}v_2,
\end{array}
\end{equation}
to be solved in the space of $T$-periodic functions. The non-trivial kernel of the operator
$\frac{d}{d\tau}-A(\tau)+i\omega$ is spanned by the complex eigenfunction $v_2$. So, the following Fredholm solvability condition has to be imposed
\begin{eqnarray*}
\int_0^{T} \langle v_2^*,C(\tau;v_1,v_1,v_2)+B(\tau;v_2,h_{200})+2B(\tau;v_1,h_{110})\\ 
-2 \alpha_{200} \dot v_2-2b_{210}v_2-2i\omega \alpha_{200}v_2\rangle \; d\tau &=& 0.
\end{eqnarray*}
From this, the expression for the normal form coefficient $b_{210}$ can be derived, namely
\begin{equation} \label{eq:b210_PDNS}
 \boxed{b_{210}=\frac{1}{2}\!\int_0^{T}\! \langle v_2^*,C(\tau;v_1,v_1,v_2)+B(\tau;v_2,h_{200})+2B(\tau;v_1,h_{110})-2 \alpha_{200} A(\tau) v_2 \rangle d \tau,}
\end{equation}
taking the normalization from \eqref{eq:AdjEigenFunc_PDNS_3} into account. Now, $h_{210}$ is defined by \eqref{eq:h210_PDNS} up to a multiple of $v_2$, therefore we impose the orthogonality condition with the adjoint eigenfunction $v_2^*$ to obtain $h_{210}$ as the unique solution of 
\begin{equation}  \label{eq:BVP_h210_PDNS}
\left\{\begin{array}{rcl}
\dot h_{210}-A(\tau) h_{210}+i\omega h_{210}-C(\tau;v_1,v_1,v_2)-B(\tau;v_2,h_{200})\\
-2B(\tau;v_1,h_{110})+2 \alpha_{200} A(\tau) v_2+2b_{210}v_2 & = & 0,\ \tau \in [0,T], \\
h_{210}(T)- h_{210}(0) & = & 0,\\
\int_0^{T} \langle v_2^*,h_{210}\rangle \; d\tau & = & 0.
\end{array}
\right.
\end{equation}

Since $\xi_1\xi_2^2$ is not a term in the normal form \eqref{eq:NF-PDNS}, we will find a non-singular equation for $h_{120}$ when collecting the $\xi_1\xi_2^2$-terms from the homological equation. Moreover, because of property \eqref{propertyhijk} $h_{120}$ is anti-periodic and thus 
\begin{equation}  \label{eq:BVP_h120_PDNS}
\left\{\begin{array}{rcl}
\dot h_{120}-A(\tau) h_{120}+2i\omega h_{120}-C(\tau;v_1,v_2,v_2)\\
-B(\tau;v_1,h_{020})-2B(\tau;v_2,h_{110}) & = & 0,\ \tau \in [0,T], \\
h_{120}(T)+ h_{120}(0) & = & 0.
\end{array}
\right.
\end{equation}

The two remaining third order terms are the $\xi_2\left|\xi_2\right|^2$-terms and the 
$\xi_1\left|\xi_2\right|^2$-terms, which both give a singular equation, namely
\begin{eqnarray*}\label{eq:h021_PDNS}
 \dot h_{021}-A(\tau) h_{021}+i\omega h_{021}&=&C(\tau;v_2,v_2,\bar v_2)+B(\tau;\bar v_2,h_{020})+
 2B(\tau;v_2,h_{011})\\
 &-&2 \alpha_{011} \dot v_2-2b_{021}v_2-2i\omega\alpha_{011}v_2
\end{eqnarray*}
and 
\begin{eqnarray*}\label{eq:h111_PDNS}
 \dot h_{111}-A(\tau) h_{111}&=&C(\tau;v_1,v_2,\bar v_2)+B(\tau;v_1,h_{011})+B(\tau;v_2,h_{101})+
 B(\tau;\bar v_2,h_{110})\\
 &-& \alpha_{011} \dot v_1-a_{111}v_1.
\end{eqnarray*}
The first function is $T$-periodic, the second one is anti-periodic. Both involve a Fredholm solvability condition, which leads to the computation of the two remaining unknown third order normal form coefficients of \eqref{eq:NF-PDNS}, i.e.
\begin{equation} \label{eq:b021_PDNS}
 \boxed{b_{021}=\frac{1}{2} \int_0^{T} \langle v_2^*,C(\tau;v_2,v_2,\bar v_2)+B(\tau;\bar v_2,h_{020})+2B(\tau;v_2,h_{011})-2 \alpha_{011} A(\tau) v_2 \rangle\; d \tau}
\end{equation}
and 
\begin{equation} \label{eq:a111_PDNS}
 \boxed{a_{111}= \int_0^{T} \langle v_1^*,C(\tau;v_1,v_2,\bar v_2)+B(\tau;v_1,h_{011})+2\Re(B(\tau;v_2,h_{101}))- \alpha_{011}A(\tau)v_1\rangle\; d \tau.}
\end{equation}
Since we need both $h_{021}$ and $h_{111}$ for the computation of higher order normal form coefficients, we also write 
down their BVPs
\begin{equation}  \label{eq:BVP_h021_PDNS}
\left\{\begin{array}{rcl}
\dot h_{021}-A(\tau) h_{021}+i\omega h_{021}-C(\tau;v_2,v_2,\bar v_2)-B(\tau;\bar v_2,h_{020})\\
-2B(\tau;v_2,h_{011})+2 \alpha_{011} A(\tau) v_2+2b_{021}v_2 & = & 0,\ \tau \in [0,T], \\
h_{021}(T)- h_{021}(0) & = & 0,\\
\int_0^{T} \langle v_2^*,h_{021}\rangle \; d\tau & = & 0
\end{array}
\right.
\end{equation}
and 
\begin{equation}  \label{eq:BVP_h111_PDNS}
\left\{\begin{array}{rcl}
\dot h_{111}-A(\tau) h_{111}-C(\tau;v_1,v_2,\bar v_2)-B(\tau;v_1,h_{011})\\
-2\Re(B(\tau;v_2,h_{101}))+\alpha_{011} A(\tau) v_1+a_{111}v_1 & = & 0,\ \tau \in [0,T], \\
h_{111}(T)+ h_{111}(0) & = & 0,\\
\int_0^{T} \langle v_1^*,h_{111}\rangle \; d\tau & = & 0.
\end{array}
\right.
\end{equation}

The stability of a possibly existing torus depends on the fourth and fifth order coefficients, which we have listed in Appendix \ref{Appendix:5}.

\subsection{NSNS}\label{Section:NSNS}
The five-dimensional critical center manifold $W^c(\Gamma)$ at
the NSNS bifurcation can be parametrized locally by
$(\xi,\tau) \in \C^2 \times [0,T]$ as 
\begin{equation}
u=u_0(\tau)+\xi_1 v_1(\tau)+\bar \xi_1 \bar v_1(\tau) + \xi_2 v_2(\tau) +\bar \xi_2 \bar v_2(\tau)+ 
H(\xi,\tau), \label{eq:CM_NSNS}
\end{equation}
where $H$ satisfies $H(\xi,T)=H(\xi,0)$ and has the Taylor
expansion %%
\begin{equation}
 H(\xi,\tau)= \sum_{2\leq i+j+k+l\leq 5} \frac{1}{i!j!k!l!}h_{ijkl}(\tau)
 \xi_1^i\bar{\xi}_1^j\xi_2^k\bar{\xi}_2^l + O(\|\xi\|^6) , 
\label{eq:H_NSNS}
\end{equation}
where the complex eigenfunctions $v_1$ and $v_2$ are given by 
\begin{eqnarray} \label{eq:EigenFunc_NSNS}
&&\left\{\begin{array}{rcl}
\dot{v}_1-A(\tau)v_1 +i\omega_1 v_1 & = & 0,\ \tau \in [0,T], \\
v_1(T)-v_1(0) & = & 0,\\
\int_{0}^{T} {\langle v_1,v_1\rangle d\tau}-1 & = & 0,\\
\end{array}
\right. 
\end{eqnarray}
and
\begin{eqnarray}\label{eq:EigenFunc_NSNS_2}
&&\left\{\begin{array}{rcl}
\dot{v}_2-A(\tau)v_2+i\omega_2 v_2 & = & 0,\ \tau \in [0,T], \\
v_2(T) - v_2(0) & = & 0,\\
\int_{0}^{T} {\langle v_2,v_2\rangle d\tau} - 1  & = & 0.\\
\end{array}
\right.
\end{eqnarray}

The functions $v_1$ and $v_2$ exist because of Lemma~2 of \cite{Io:88}. %%
The functions $h_{ijkl}$ will be found by solving appropriate
BVPs, assuming that \eqref{eq:P.1} restricted to $W^c(\Gamma)$ has
the normal form \eqref{eq:NF-NSNS}. 

The coefficients of the normal
form  arise from the solvability conditions for the BVPs as
integrals of scalar products over the interval $[0,T]$.
Specifically, those scalar products involve among other things the quadratic and
cubic terms of \eqref{eq:MULT} near the periodic solution $u_0$,
the eigenfunctions $v_1$ and $v_2$, and the
adjoint eigenfunctions $\varphi^*$, $v_1^*$ and $v_2^*$ as solution
of the problems
\begin{eqnarray}
&&\left\{\begin{array}{rcl}
 \dot{\varphi}^*+A^{\rm T}(\tau)\varphi^* & = & 0,\ \tau \in [0,T], \\
 \varphi^*(T)-\varphi^*(0) & = & 0, \label{eq:AdjEigenFunc_NSNS}\\
 \int_{0}^{T} {\langle \varphi^*,F(u_0) \rangle d\tau} -1 & = & 0,
 \end{array} \right.
 \end{eqnarray}
 \begin{eqnarray}
 &&\left\{\begin{array}{rcl}
 \dot{v}_1^*+A^{\rm T}(\tau) v_1^*+i\omega_1v_1^* & = & 0,\ \tau \in [0,T], \\
 v_1^*(T)-v_1^*(0) & = & 0, \label{eq:AdjEigenFunc_NSNS_2}\\
 \int_{0}^{T} {\langle v_1^*,v_1 \rangle d\tau}-1 & = & 0,
 \end{array} \right. 
 \end{eqnarray}
and
 \begin{eqnarray}
 &&\left\{\begin{array}{rcl}
 \dot{v}_2^*+A^{\rm T}(\tau) v_2^*+i\omega_2 v_2^* & = & 0,\ \tau \in [0,T], \\
 v_2^*(T)-v_2^*(0) & = & 0, \label{eq:AdjEigenFunc_NSNS_3}\\
 \int_{0}^{T} {\langle v_2^*,v_2\rangle d\tau} -1 & = & 0.
 \end{array} \right.
\end{eqnarray}

By collecting the constant and linear terms we get the identities
\begin{eqnarray}
 \dot{u}_0=F(u_0), \qquad \dot v_1+ i\omega_1 v_1=A(\tau) v_1, \qquad \dot{v}_2+i\omega_2 v_2=A(\tau) v_2,
\label{basiceqnsns}
\end{eqnarray}
and the complex conjugates of the last two equations. (\ref{basiceqnsns}) merely reflects the definition of $u_0$ and the first equations in \eqref{eq:EigenFunc_NSNS}, \eqref{eq:EigenFunc_NSNS_2}.

By collecting the $\xi_1^2$ (or $\bar \xi_1^2$-terms)-terms we find an equation for $h_{2000}$
\begin{equation*}\label{eq:h2000_NSNS}
 \dot h_{2000}-A(\tau) h_{2000}+2i\omega_1h_{2000}=B(\tau;v_1,v_1),
\end{equation*}
(or its complex conjugate). This equation has a unique solution $h_{2000}$ satisfying
$h_{2000}(T)=h_{2000}(0)$, since due to the spectral assumptions $e^{2i\omega_1 T}$ is not a multiplier of the critical cycle. Thus, $h_{2000}$ can be found by solving
\begin{equation}  \label{eq:BVP_h2000_NSNS}
\left\{\begin{array}{rcl}
\dot h_{2000}-A(\tau) h_{2000}+2i\omega_1h_{2000}-B(\tau;v_1,v_1) & = & 0,\ \tau \in [0,T], \\
h_{2000}(T)- h_{2000}(0) & = & 0.
\end{array}
\right.
\end{equation}
The function $h_{0200}$ is just the complex conjugate of the function $h_{2000}$. Analogously, by comparing the $\xi_2^2$-terms, we find that $h_{0020}$ is the unique solution of 
\begin{equation}  \label{eq:BVP_h0020_NSNS}
\left\{\begin{array}{rcl}
\dot h_{0020}-A(\tau) h_{0020}+2i\omega_2h_{0020}-B(\tau;v_2,v_2) & = & 0,\ \tau \in [0,T], \\
h_{0020}(T)- h_{0020}(0) & = & 0.
\end{array}
\right.
\end{equation}

By collecting the $\left|\xi_1\right|^2$-terms we obtain a singular equation, as expected since this term is present in the normal form \eqref{eq:NF-NSNS}, namely
\begin{equation*}\label{eq:h1100_NSNS}
 \dot h_{1100}-A(\tau) h_{1100}=B(\tau;v_1,\bar v_1)-\alpha_{1100}\dot u_0,
\end{equation*}
to be solved in the space of functions satisfying
$h_{1100}(T)=h_{1100}(0)$. Since the null-space is
spanned by $\dot u_0$, the Fredholm solvability
condition
\[
 \int_0^{T} \langle \varphi^*,B(\tau;v_1,\bar v_1)-\alpha_{1100}\dot u_0 \rangle\; d \tau =0
\]
gives us the possibility to calculate parameter $\alpha_{1100}$ due to the normalization condition in \eqref{eq:AdjEigenFunc_NSNS}, i.e.
\begin{equation} \label{eq:alpha1100_NSNS}
 \alpha_{1100}=\int_0^{T} \langle \varphi^*,B(\tau;v_1,\bar v_1)\rangle\; d \tau.
\end{equation}
Function $h_{1100}$ is now determined up to the addition of a multiple of $\dot u_0$. As always, we will add an orthogonality condition, in this case with the adjoint eigenfunction coresponding with multiplier $1$. Therefore, with the value of $\alpha_{1100}$ from \eqref{eq:alpha1100_NSNS} we obtain $h_{1100}$ as the unique solution of the BVP
\begin{equation}  \label{eq:BVP_h1100_NSNS}
\left\{\begin{array}{rcl}
\dot{h}_{1100}-A(\tau) h_{1100} - B(\tau;v_1,\bar v_1) + \alpha_{1100} \dot u_0 & = & 0,\ \tau \in [0,T], \\
h_{1100}(T)- h_{1100}(0) & = & 0,\\
\int_0^{T} \langle \varphi^*,h_{1100}\rangle \; d\tau & = & 0.
\end{array}
\right.
\end{equation}

Analogously, function $h_{0011}$ can be obtained by solving
\begin{equation}  \label{eq:BVP_h0011_NSNS}
\left\{\begin{array}{rcl}
\dot{h}_{0011}-A(\tau) h_{0011} - B(\tau;v_2,\bar v_2) + \alpha_{0011} \dot u_0 & = & 0,\ \tau \in [0,T], \\
h_{0011}(T)- h_{0011}(0) & = & 0,\\
\int_0^{T} \langle \varphi^*,h_{0011}\rangle \; d\tau & = & 0,
\end{array}
\right.
\end{equation}
with 
\begin{equation} \label{eq:alpha0011_NSNS}
 \alpha_{0011}=\int_0^{T} \langle \varphi^*,B(\tau;v_2,\bar v_2)\rangle\; d \tau.
\end{equation}

By collecting the $\xi_1\xi_2$-terms we find the following differential equation for $h_{1010}$
\begin{equation*}\label{eq:h1010_NSNS}
 \dot h_{1010}-A(\tau) h_{1010}+i\omega_1 h_{1010}+i\omega_2 h_{1010}=B(\tau;v_1,v_2).
\end{equation*}
This equation has a unique solution $h_{1010}$ satisfying
$h_{1010}(T)=h_{1010}(0)$, since due to the spectral assumptions $e^{i(\omega_1+\omega_2) T}$ is not a multiplier of the critical cycle. Thus, $h_{1010}$ can be found by solving
\begin{equation}  \label{eq:BVP_h1010_NSNS}
\left\{\begin{array}{rcl}
\dot h_{1010}-A(\tau) h_{1010}+i\omega_1 h_{1010}+i\omega_2 h_{1010}-B(\tau;v_1,v_2) & = & 0,\ \tau \in [0,T], \\
h_{1010}(T)- h_{1010}(0) & = & 0.
\end{array}\right.
\end{equation}
We note that $h_{0101} = \overline{h_{1010}}$.

The last second order derivative coming from looking at the $\xi_1\bar \xi_2$-terms results in a non-singular differential equation, such that 
\begin{equation}  \label{eq:BVP_h1001_NSNS}
\left\{\begin{array}{rcl}
\dot h_{1001}-A(\tau) h_{1001}+i\omega_1 h_{1001}-i\omega_2 h_{1001}-B(\tau;v_1,\bar v_2) & = & 0,\ \tau \in [0,T], \\
h_{1001}(T)- h_{1001}(0) & = & 0.
\end{array}\right.
\end{equation}

We now come to the third order terms. From the $\xi_1^3$ and $\xi_2^3$-terms we immediately get the BVPs for $h_{3000}$ and $h_{0030}$, namely
\begin{equation}  \label{eq:BVP_h3000_NSNS}
\left\{\begin{array}{rcl}
 \dot h_{3000}-A(\tau) h_{3000}+3i\omega_1 h_{3000}-C(\tau;v_1,v_1,v_1)\\
-3B(\tau;v_1,h_{2000}) & = & 0,\ \tau \in [0,T], \\
h_{3000}(T)- h_{3000}(0) & = & 0
\end{array}\right.
\end{equation}
and 
\begin{equation}  \label{eq:BVP_h0030_NSNS}
\left\{\begin{array}{rcl}
 \dot h_{0030}-A(\tau) h_{0030}+3i\omega_2 h_{0030}-C(\tau;v_2,v_2,v_2)\\
-3B(\tau;v_2,h_{0020}) & = & 0,\ \tau \in [0,T], \\
h_{0030}(T)- h_{0030}(0) & = & 0.
\end{array}\right.
\end{equation}

Since the $\xi_1\left|\xi_1\right|^2$-term is present in the normal form for the double Neimark-Sacker bifurcation, a Fredholm solvability condition is involved coming from the following differential equation for $h_{2100}$
\begin{eqnarray*}
 \dot h_{2100}-A(\tau) h_{2100}+i\omega_1 h_{2100}&=&C(\tau;v_1,v_1,\bar v_1)+2B(\tau;v_1,h_{1100})+B(\tau;\bar v_1,h_{2000})\nonumber\\
 &&-2a_{2100}v_1-2i\omega_1\alpha_{1100}v_1-2\alpha_{1100} \dot v_1.\nonumber\\\label{eq:h2100_NSNS}
\end{eqnarray*}
The differential operator
$\frac{d}{d\tau}-A(\tau)+i\omega_1$ is singular with its null-space
spanned by $v_1$, so we get
condition
\begin{eqnarray*}
 \int_0^{T} \langle v_1^*,C(\tau;v_1,v_1,\bar v_1)+2B(\tau;v_1,h_{1100})+B(\tau;\bar v_1,h_{2000})-2a_{2100}v_1 && \\
 -2\alpha_{1100} \dot v_1 -2i\omega_1\alpha_{1100}v_1\rangle\; d \tau &=& 0.
 \end{eqnarray*}
Taking the normalization condition from \eqref{eq:AdjEigenFunc_NSNS_2} and the differential equation from \eqref{eq:EigenFunc_NSNS} into account, we get
\begin{equation} \label{eq:a2100_NSNS}
 \boxed{a_{2100}=\frac{1}{2} \int_0^{T} \langle v_1^*,C(\tau;v_1,v_1,\bar v_1)+2B(\tau;v_1,h_{1100})+B(\tau;\bar v_1,h_{2000})-2\alpha_{1100} A(\tau) v_1\rangle\; d \tau.}
\end{equation}
Therefore, we can compute $h_{2100}$ as the unique solution of the BVP
\begin{equation}  \label{eq:BVP_h2100_NSNS}
\left\{\begin{array}{rcl}
\dot h_{2100}-A(\tau) h_{2100}+i\omega_1 h_{2100}-C(\tau;v_1,v_1,\bar v_1)\\
-2B(\tau;v_1,h_{1100})-B(\tau;\bar v_1,h_{2000})+2a_{2100}v_1+2\alpha_{1100}A(\tau) v_1& = & 0,\ \tau \in [0,T], \\
h_{2100}(T)- h_{2100}(0) & = & 0,\\
\int_0^{T} \langle v_1^*,h_{2100}\rangle \; d\tau & = & 0.
\end{array}
\right.
\end{equation}

We can now immediately list the following four BVPs
\begin{equation}  \label{eq:BVP_h2010_NSNS}
\left\{\begin{array}{rcl}
 \dot h_{2010}-A(\tau) h_{2010}+2i\omega_1 h_{2010}+i\omega_2 h_{2010}-C(\tau;v_1,v_1,v_2)&&\\
 -B(\tau;v_2,h_{2000})-2B(\tau;v_1,h_{1010})& = & 0,\ \tau \in [0,T], \\
h_{2010}(T)- h_{2010}(0) & = & 0,
\end{array}\right.
\end{equation}
\begin{equation}  \label{eq:BVP_h2001_NSNS}
\left\{\begin{array}{rcl}
 \dot h_{2001}-A(\tau) h_{2001}+2i\omega_1 h_{2001}-i\omega_2 h_{2001}-C(\tau;v_1,v_1,\bar v_2)\\
-B(\tau;\bar v_2,h_{2000})-2B(\tau;v_1,h_{1001})& = & 0,\ \tau \in [0,T], \\
h_{2001}(T)- h_{2001}(0) & = & 0,
\end{array}\right.
\end{equation}
\begin{equation}  \label{eq:BVP_h1020_NSNS}
\left\{\begin{array}{rcl}
 \dot h_{1020}-A(\tau) h_{1020}+i\omega_1 h_{1020}+2i\omega_2 h_{1020}\\
-C(\tau;v_1,v_2,v_2)-B(\tau;v_1,h_{0020})-2B(\tau;v_2,h_{1010})& = & 0,\ \tau \in [0,T], \\
h_{1020}(T)- h_{1020}(0) & = & 0,
\end{array}\right.
\end{equation}
and
\begin{equation}  \label{eq:BVP_h0120_NSNS}
\left\{\begin{array}{rcl}
 \dot h_{0120}-A(\tau) h_{0120}-i\omega_1 h_{0120}+2i\omega_2 h_{0120}\\
-C(\tau;\bar v_1,v_2,v_2)-B(\tau;\bar v_1,h_{0020})-2B(\tau;v_2,h_{0110})& = & 0,\ \tau \in [0,T], \\
h_{0120}(T)- h_{0120}(0) & = & 0.
\end{array}\right.
\end{equation}

The $\xi_2\left|\xi_2\right|^2$-terms from the homological equation make it possible to compute $b_{0021}$. Indeed, the differential equation
\begin{eqnarray*}
 \dot h_{0021}-A(\tau) h_{0021}+i\omega_2 h_{0021}&=&C(\tau;v_2,v_2,\bar v_2)+B(\tau;\bar v_2,h_{0020})+2B(\tau;v_2,h_{0011})\nonumber\\
&& -2b_{0021}v_2-2\alpha_{0011}\dot v_2-2i\omega_2\alpha_{0011}v_2 \nonumber\\\label{eq:h0021_NSNS}
\end{eqnarray*}
results in a solvability
condition with $v_2^*$, i.e.
\begin{eqnarray*}
 \int_0^{T} \langle v_2^*,C(\tau;v_2,v_2,\bar v_2)+B(\tau;\bar v_2,h_{0020})+2B(\tau;v_2,h_{0011})\\
 -2b_{0021}v_2-2\alpha_{0011}\dot v_2-2i\omega_2\alpha_{0011}v_2 \rangle\; d \tau &=&0.
 \end{eqnarray*}
Therefore, considering the normalization condition from \eqref{eq:AdjEigenFunc_NSNS_3} and the differential equation from \eqref{eq:EigenFunc_NSNS_2}, we can calculate parameter $b_{0021}$
\begin{equation} \label{eq:b0021_NSNS}
 \boxed{b_{0021}=\frac{1}{2} \int_0^{T} \langle v_2^*,C(\tau;v_2,v_2,\bar v_2)+B(\tau;\bar v_2,h_{0020})+2B(\tau;v_2,h_{0011})\\
 -2\alpha_{0011}A(\tau) v_2 \rangle\; d \tau,}
\end{equation}
with $h_{0021}$ as the unique solution of the BVP
\begin{equation}  \label{eq:BVP_h0021_NSNS}
\left\{\begin{array}{rcl}
\dot h_{0021}-A(\tau) h_{0021}+i\omega_2 h_{0021}-C(\tau;v_2,v_2,\bar v_2)\\
-B(\tau;\bar v_2,h_{0020})-2B(\tau;v_2,h_{0011})\\
+2b_{0021}v_2+2\alpha_{0011}A(\tau) v_2& = & 0,\ \tau \in [0,T], \\
h_{0021}(T)- h_{0021}(0) & = & 0,\\
\int_0^{T} \langle v_2^*,h_{0021}\rangle \; d\tau & = & 0.
\end{array}
\right.
\end{equation}

The last two third order terms which we have to examine give us both the formula for a normal form coefficient. The first one, obtained from the $\left|\xi_1\right|^2\xi_2$-terms, gives us the BVP
\begin{equation}  \label{eq:BVP_h1110_NSNS}
\left\{\begin{array}{rcl}
\dot h_{1110}-A(\tau) h_{1110}+i\omega_2 h_{1110}-C(\tau;v_1,\bar v_1,v_2)-B(\tau;v_1,h_{0110})&&\\
-B(\tau;\bar v_1,h_{1010})-B(\tau;v_2,h_{1100})\\
+b_{1110}v_2+\alpha_{1100}A(\tau) v_2& = & 0,\ \tau \in [0,T], \\
h_{1110}(T)- h_{1110}(0) & = & 0,\\
\int_0^{T} \langle v_2^*,h_{1110}\rangle \; d\tau & = & 0,
\end{array}
\right.
\end{equation}
where from the solvability condition it follows that 
\begin{equation}
 \boxed{b_{1110}=\int_0^{T} \langle v_2^*,C(\tau;v_1,\bar v_1,v_2)+B(\tau;v_1,h_{0110})+B(\tau;\bar v_1,h_{1010})+B(\tau;v_2,h_{1100}) \nonumber \\
-\alpha_{1100}A(\tau)v_2  \rangle\; d \tau.} \label{eq:b1110_NSNS}
\end{equation}
Analogously, we obtain the BVP
\begin{equation}  \label{eq:BVP_h1011_NSNS}
\left\{\begin{array}{rcl}
\dot h_{1011}-A(\tau) h_{1011}+i\omega_1 h_{1011}-C(\tau;v_1,v_2,\bar v_2)-B(\tau;v_1,h_{0011})\\
-B(\tau;v_2,h_{1001})-B(\tau;\bar v_2,h_{1010})+a_{1011}v_1+\alpha_{0011}A(\tau) v_1& = & 0,\ \tau \in [0,T], \\
h_{1011}(T)- h_{1011}(0) & = & 0,\\
\int_0^{T} \langle v_1^*,h_{1011}\rangle \; d\tau & = & 0,
\end{array}
\right.
\end{equation}
with 
\begin{equation} 
 \boxed{a_{1011}=\int_0^{T} \langle v_1^*,C(\tau;v_1,v_2,\bar v_2)+B(\tau;v_1,h_{0011})+B(\tau;v_2,h_{1001})+B(\tau;\bar v_2,h_{1010}) \nonumber \\
-\alpha_{0011}A(\tau)v_1  \rangle\; d \tau.} \label{eq:a1011_NSNS}
\end{equation}
We still need the coefficients $b_{1101}$ and $a_{0111}$ which are determined by 
\begin{eqnarray}
 b_{1101}=\int_0^{T} \langle \bar v_2^*,C(\tau;v_1,\bar v_1,\bar v_2)+B(\tau;v_1,h_{0101})+B(\tau;\bar v_1,h_{1001})+B(\tau;\bar v_2,h_{1100}) \nonumber \\
-\alpha_{1100}A(\tau)\bar v_2  \rangle\; d \tau && 
\end{eqnarray}
and
\begin{eqnarray}
 a_{0111}=\int_0^{T} \langle \bar v_1^*,C(\tau;\bar v_1,v_2,\bar v_2)+B(\tau;\bar v_1,h_{0011})+B(\tau;v_2,h_{0101})+B(\tau;\bar v_2,h_{0110}) \nonumber\\
-\alpha_{0011}A(\tau)\bar v_1  \rangle\; d \tau. &&
\end{eqnarray}

As before, the higher order terms which determine the stability of the torus can be found in Appendix \ref{Appendix:5}.

\subsection{Implementation}\label{Section:Implementation}
Numerical implementation of the formulas derived in the previous section requires the evaluation of 
integrals of scalar functions over $[0,T]$ and the solution of nonsingular linear BVPs with integral 
constraints. Such tasks can be carried out within the standard continuation software such as {\sc auto} 
\cite{AUTO97}, {\sc content} \cite{CONTENT}, and {\sc matcont} \cite{MATCONT}. In these software packages, 
periodic solutions to (\ref{ODE}) are computed with the method of {\em orthogonal collocation} with 
piecewise polynomials applied to properly formulated BVPs \cite{BoSw:73,BVPbook:95}. 

We have implemented our algorithms in {\sc matcont} analogously to the eight cases with $n_c\leq3$. For 
further details we refer to \cite{drdwgk} where this is extensively discussed.

\section{Examples}
\label{Section:Examples}

\subsection{Laser model}

\begin{figure}[htbp]
\centering
 \footnotesize
 \psfrag{Omp1}[][c]{$\Omega_p$}
 \psfrag{Dcav1}[][c]{$\Delta_{cav}$}
  \includegraphics[width=.95\textwidth]{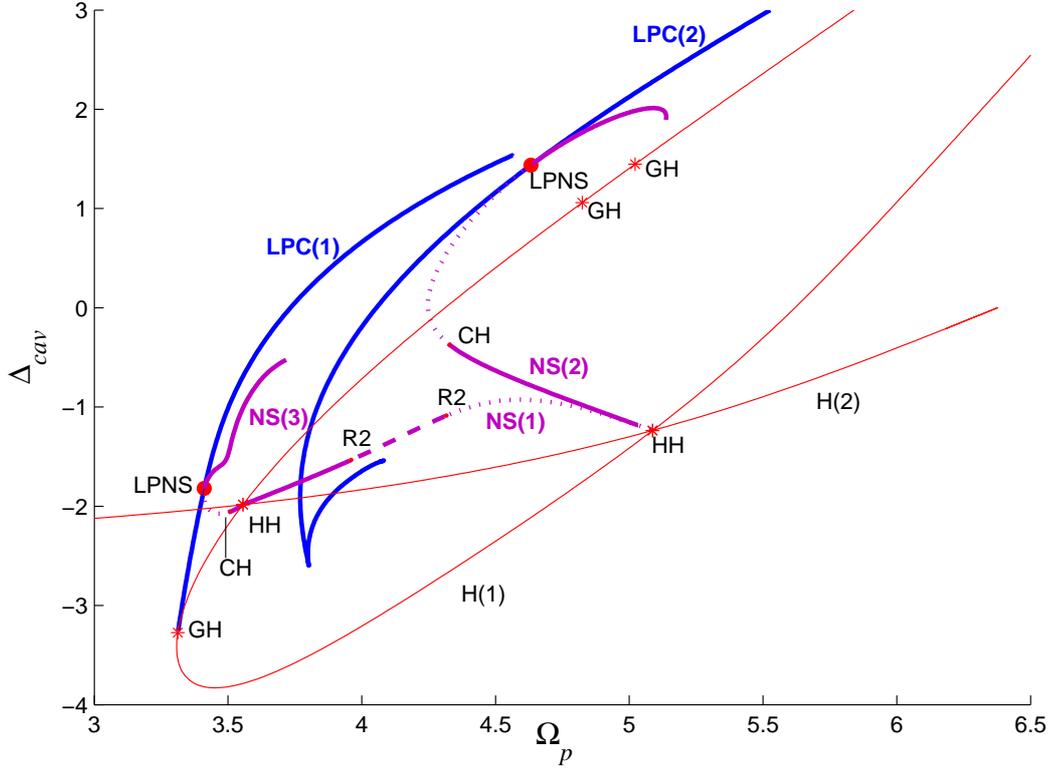}
 \caption{Bifurcation diagram of \eqref{eq:lasermodel}. The thin red curves are Hopf curves. In blue are limit point of cycles bifurcations and in magenta Neimark-Sacker bifurcations. Solid/dotted curves correspond to
 supercritical/subcritical bifurcations. The dashed curves are curves of neutral saddles.} \label{fig:lasermodel}
\end{figure}

In \cite{WieChow:06} a single-mode inversionless laser with a three-level phaser was studied and 
shown to operate in various modes. These modes are ``off" (non-lasing), continuous waves, periodic, 
quasi-periodic and chaotic lasing. The model is a $9$-dimensional system given by $3$ real and $3$ 
complex equations:
\begin{equation} \label{eq:lasermodel}
\begin{cases}
\dot \Omega_l = - \frac{\gamma_{cav}}{2}\Omega_l-g \Im(\sigma_{ab})\\
\dot \rho_{aa} = R_a- \frac{i}{2}(\Omega_l(\sigma_{ab}-\sigma_{ab}^*)+\Omega_p(\sigma_{ac}-\sigma_{ac}^*))\\
\dot \rho_{bb} = R_b + \frac{i}{2} \Omega_l (\sigma_{ab}-\sigma_{ab}^*)\\
\dot \sigma_{ab} = -(\gamma_1+i\Delta_l)\sigma_{ab}- \frac{i}{2}(\Omega_l(\rho_{aa}-\rho_{bb})-\Omega_p\sigma_{cb})\\
\dot \sigma_{ac} = -(\gamma_2+i\Delta_p)\sigma_{ac}- \frac{i}{2}(\Omega_p(2\rho_{aa}+\rho_{bb}-1)-\Omega_l\sigma_{cb}^*)\\
\dot \sigma_{cb} = -(\gamma_3+i(\Delta_l-\Delta_p))\sigma_{cb}- \frac{i}{2}(\Omega_l\sigma_{ac}^*-\Omega_p\sigma_{ab}),
\end{cases}
\end{equation}
with $R_a=-0.505 \rho_{aa}-0.405\rho_{bb}+0.45, R_b=0.0495 \rho_{aa}-0.0505\rho_{bb}+0.0055$ and 
$\Delta_l=\Delta_{cav}+g\Re(\sigma_{ab})\Omega_l$. The fixed parameters are $\gamma_1=0.275, 
\gamma_2=0.25525, \gamma_3=0.25025, \gamma_{cav}=0.03, g=100, \Delta_p=0$. The parameters 
$\Omega_p$ and $\Delta_{cav}$ are varied. The bifurcation diagram of \eqref{eq:lasermodel} is computed
in \cite{KuMeGoSa:2008} and is reproduced in Figure \ref{fig:lasermodel}. 

\subsubsection{The LPNS points}
\begin{figure}[htbp]
\centering
 \footnotesize
% \psfrag{Omegap}[][c]{$\Omega_p$}
% \psfrag{Deltacav}[][c]{$\Delta_{cav}$}
  \includegraphics[width=.42\textwidth]{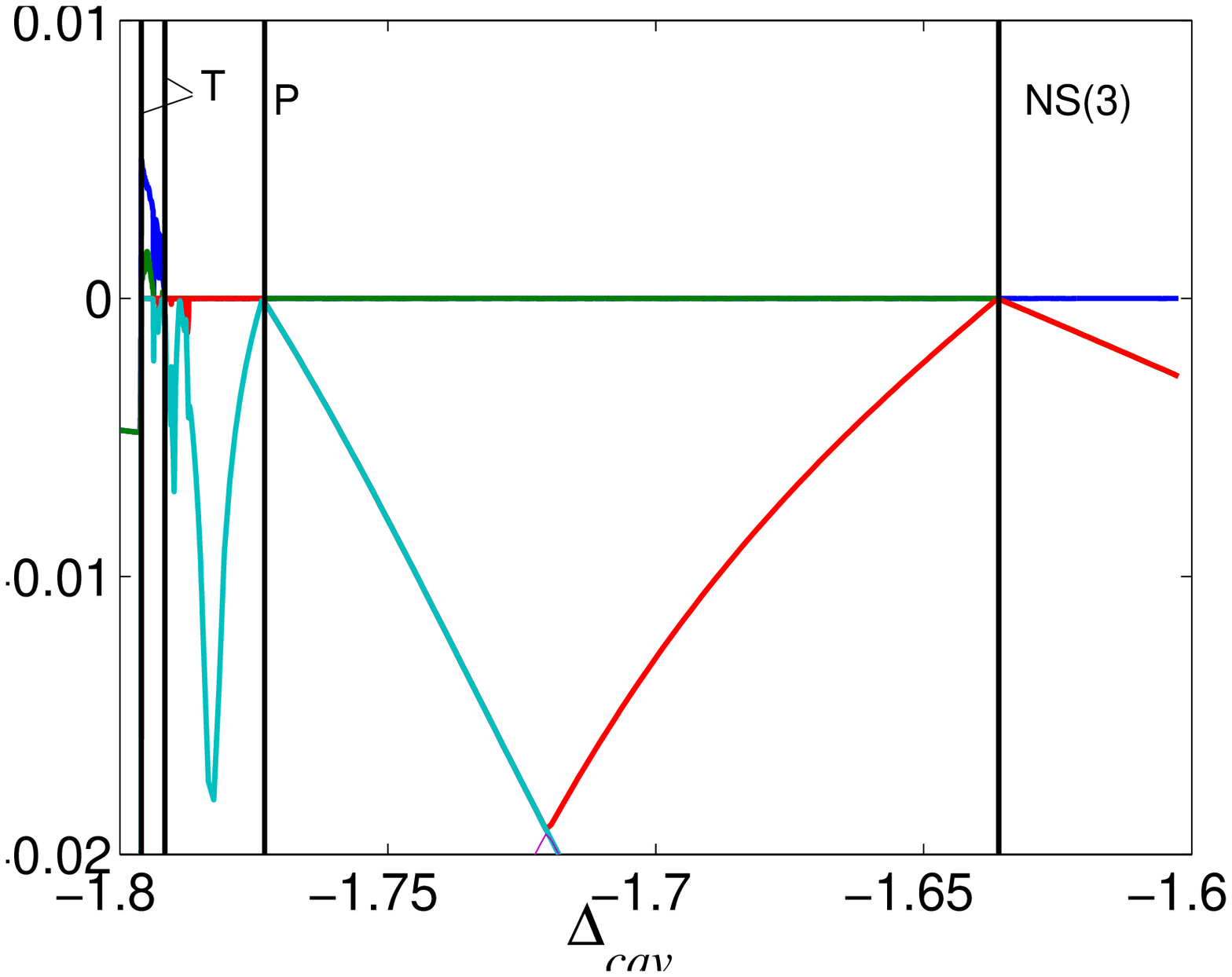}
  \includegraphics[width=.42\textwidth]{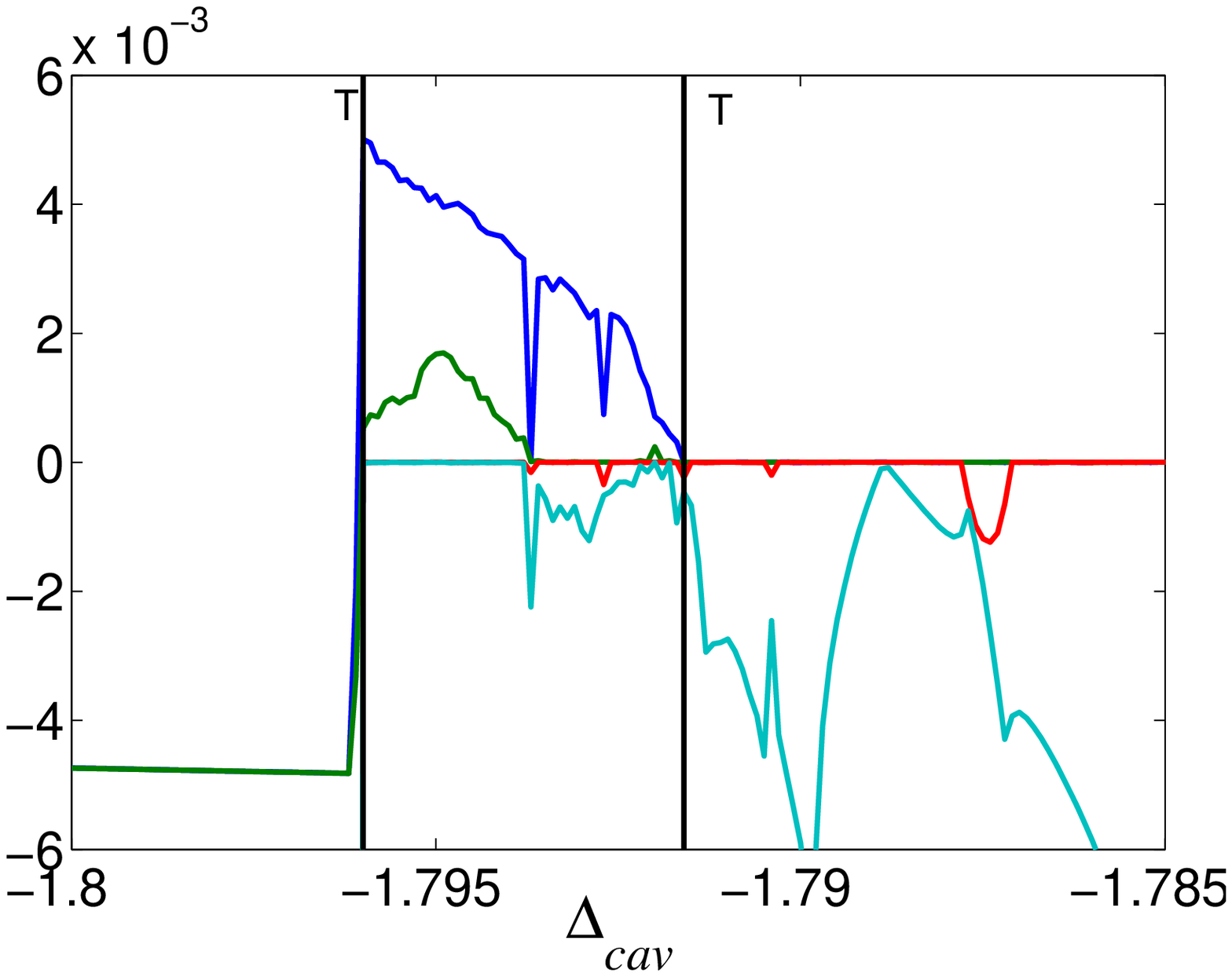}
 \caption{Lyapunov exponents computed for $\Omega_p=3.45$ close to the LPNS point
  at $(\Omega_p,\Delta_{cav})=(3.411,-1.819)$, (left) for $\Delta_{cav}\in[-1.8;-1.6]$ 
  and (right) zoomed in near the region with chaos due to heteroclinic tangles. The vertical black
  lines indicate the parameter values where a bifurcation occurs.} \label{fig:lyapunov_lpnsleft}
\end{figure}

Figure \ref{fig:lasermodel} shows three NS curves {\sf NS(1)}, {\sf NS(2)} and {\sf NS(3)} 
starting from two {\sf HH} points. On {\sf NS(3)} one of the richer situations happens. The normal 
form coefficients for the {\sf LPNS} point at $(\Omega_p,\Delta_{cav})=(3.411,-1.819)$ are 
$(s,\theta,E)=(1, -0.139, -911.248)$, so $s\theta<0$. This means that there exists a 3-torus, 
which is stable since $\theta<0$ and $E<0$. Therefore, we are in the case represented in 
Figure \ref{fig:NF_LPNS} (c), but with a stable $3$-torus. For computing the Lyapunov exponents, 
we used a code written by V. N. Govorukhin (2004). Figure \ref{fig:lyapunov_lpnsleft} (left)
shows the calculated Lyapunov exponents for $\Omega_p$ fixed at $3.45$ and $\Delta_{cav} \in [-1.8;-1.6]$. 
More detail is shown in Figure \ref{fig:lyapunov_lpnsleft} (right), where we get a clear view on 
the number of Lyapunov exponents equal to zero. For $\Delta_{cav}$ values to the right of $-1.636$, 
there is one Lyapunov exponent equal to zero, which corresponds to the stable limit cycle from region 6 
in Figure \ref{fig:NF_LPNS} (c). At $\Delta_{cav}=-1.636$, we cross {\sf NS(3)} and arrive in region 5 with a 
stable $2$-torus and therefore two Lyapunov exponents equal to zero. When crossing the $P$ curve at 
$\Delta_{cav}=-1.773$, the stable $3$-torus from region 4 arises. Remark that in some small intervals 
only two Lyapunov exponents are equal to zero, and thus not the expected three zero ones, but these 
correspond with resonances on the $3$-torus. Then, in the interval $\Delta_{cav} \in [-1.796;-1.7916]$ 
positive Lyapunov exponents appear which indicates that there is chaos. This zone corresponds with $T$. 
Afterwards, we arrive in region 3, where all Lyapunov exponents are smaller than zero.

\begin{figure}[htbp]
\centering
 \footnotesize
 \psfrag{Omegap}[][c]{$\Omega_p$}
 \psfrag{Deltacav}[][c]{$\Delta_{cav}$}
  \includegraphics[width=.7\textwidth]{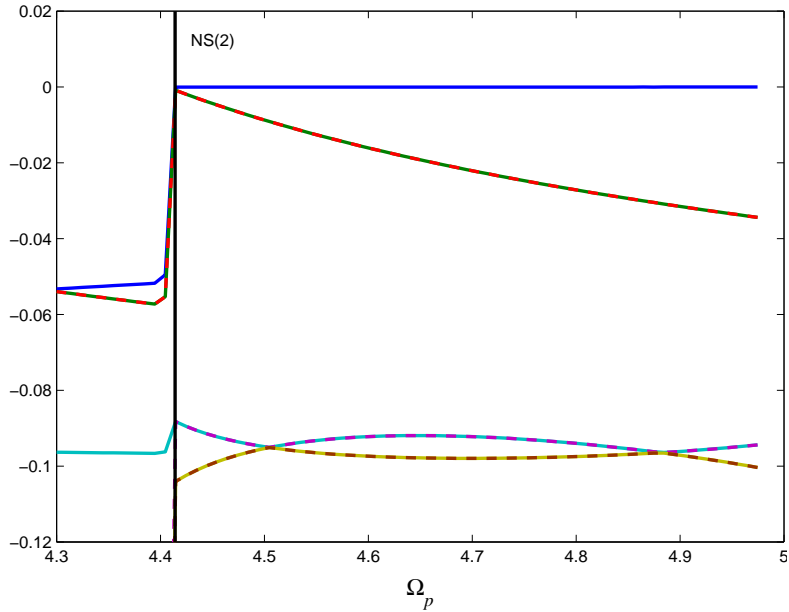}
 \caption{Lyapunov exponents computed close to the LPNS point at $(\Omega_p,\Delta_{cav})=(4.632, 1.438)$.
 The two-coloured dashed lines reveal pairs of equally large Lyapunov exponents.} 
 \label{fig:lyapunov_lpnsright}
\end{figure}

On the {\sf NS(2)} curve there is one {\sf LPNS} point for $(\Omega_p,\Delta_{cav})=(4.632, 1.438)$. 
The normal form coefficients are $(s,\theta,E)=(1,0.206,808.009)$. The product $s\theta>0$ is positive, 
so we are in a ``simple" case, where no $3$-torus is present. Since $s=1$, the torus arisen through the 
Neimark-Sacker curve exists below the {\sf NS(2)} curve. We have computed the Lyapunov exponents for a straight line 
where the beginning point $(\Omega_p,\Delta_{cav})=(4.302,0.673)$ and end point $(\Omega_p,\Delta_{cav})=
(4.984,1.984)$ lie between the curves {\sf LPC(2)} and {\sf NS(2)}, to the left and to the right of the 
{\sf LPNS point}. In Figure \ref{fig:lyapunov_lpnsright}, we plot the Lyapunov exponents for 
$\Omega_p \in [4.3,4.98]$. The stable limit cycle is situated in the upper wedge between the 
{\sf LPC(2)} and {\sf NS(2)} curves which corresponds to region $4$ in Figure  \ref{fig:NF_LPNS_1}, 
so we have one Lyapunov exponent equal to zero for $\Omega_p$-values larger than the subcritical 
{\sf NS(2)} curve. At $\Omega_p \approx 4.41$, we cross the subcritical {\sf NS(2)} curve, with to 
the left no zero Lyapunov exponents.

\subsection{Periodic predator-prey model}
As a second model we study a simple two-patch predator-prey system with periodic (seasonal) forcing. 
Simple predator-prey models lead to the `paradox of enrichment', i.e., increasing the carrying 
capacity of the prey ultimately leads to extinction of the population \cite{Rosenzweig:1971}. 
Outside the laboratory, however, stable populations are observed and not extinction. Here, spatial models have been put forward to explain this discrepancy. As the simplest spatial case, 
one may consider a two-patch predator-prey model \cite{Jansen:1995} where predator and prey can 
migrate between the two patches by diffusion. This leads to a diffusive instability of large 
oscillations and stabilizes the total population size \cite{Jansen:2001}. Here, we propose 
an extension where one of the patches experiences seasonal influences while the other can 
be seen as a wild-life refuge where human intervention minimizes seasonal influences. As a 
simplication we will only consider the case that the predators can move between the patches, 
i.e., they can cross the refuge barrier. On a proper time scale, the investigated system is defined by
\begin{equation} \label{eq:PreyPredator}
\begin{cases}
\dot x_1 = {\displaystyle r_1x_1(1-x_1)-\frac{cx_1x_2}{x_1+b_1(1+\varepsilon v_1)}}\\
\dot x_2 = {\displaystyle -x_2+\frac{cx_1x_2}{x_1+b_1(1+\varepsilon v_1)}+\gamma(y_2-x_2)}\\
\dot y_1 = {\displaystyle r_2y_1(1-y_1)-\frac{cy_1y_2}{y_1+b_2}}\\
\dot y_2 = {\displaystyle -y_2+\frac{cy_1y_2}{y_1+b_2}+\gamma(x_2-y_2)}\\
\dot v_1 = -v_2+v_1(1-v_1^2-v_2^2)\\
\dot v_2 = v_1+v_2(1-v_1^2-v_2^2).
\end{cases}
\end{equation}
The values of $x_{1}$ and $x_{2}$ denote the numbers of individuals (or densities) respectively 
of prey and predator populations living outside the refuge and $y_{1}$ and $y_{2}$ are the corresponding 
numbers or densities inside. The intrinsic growth rates $r_{i}$ and the constant attack rate $c$ are 
parameters of the model. For the predator outside the refuge, the Holling type II is 
chosen as functional response with a half saturation which varies periodically with period $2\pi$. 
To this end, the last two equations are introduced; their solutions converge to a stable limit cycle 
$v_{1}(t) =\cos(t +\phi)$ with a phase shift $\phi$ depending on the initial conditions. The terms with 
parameter $\gamma$ describe the coupling of the two patches. The fixed parameter values are 
$r_{1}=1, r_{2}=1, b_{1}=0.4, \gamma=0.1, c=2$. We will use the half saturation $b_{2}$ as a 
continuation parameter together with the amplitude of the seasonal forcing $\varepsilon$. It is not 
our aim to give a full analysis of this model, but rather analyze the codim 2 bifurcations 
relevant for this paper. We observe that a refuge can induce complex behaviour in a spatial 
population model with seasonal forcing.

\subsubsection{The PDNS points}
\begin{figure}[htbp]
\centering
 \footnotesize
% \psfrag{b2}[][c]{$b_2$}
% \psfrag{eps}[][c]{$\epsilon$}
  \includegraphics[width=.95\textwidth]{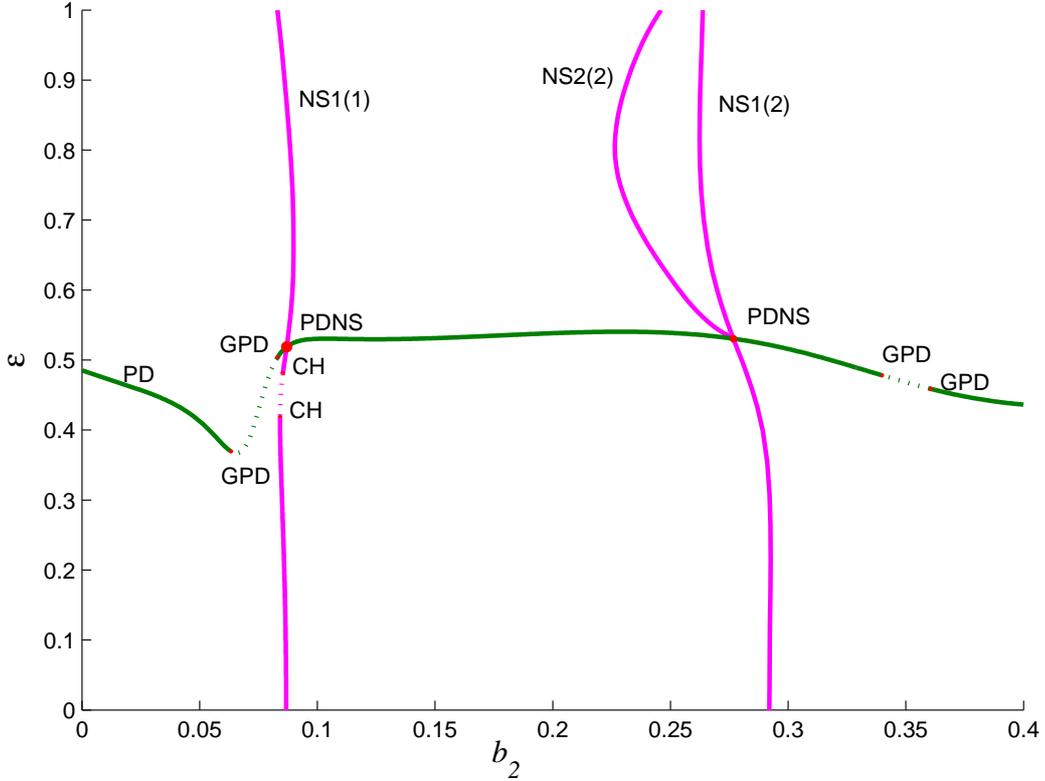}
 \caption{Bifurcation diagram of limit cycles in \eqref{eq:PreyPredator}. In green are period doubling 
 curves and in magenta Neimark-Sacker curves (of the first or of the second
iterate, respectively labeled with {\sf NS1} and {\sf NS2}). } \label{fig:periodicforced}
\end{figure}

\begin{figure}[htbp]
\centering
 \footnotesize
% \psfrag{b2}[][c]{$b_2$}
% \psfrag{eps}[][c]{$\epsilon$}
  \includegraphics[width=.95\textwidth]{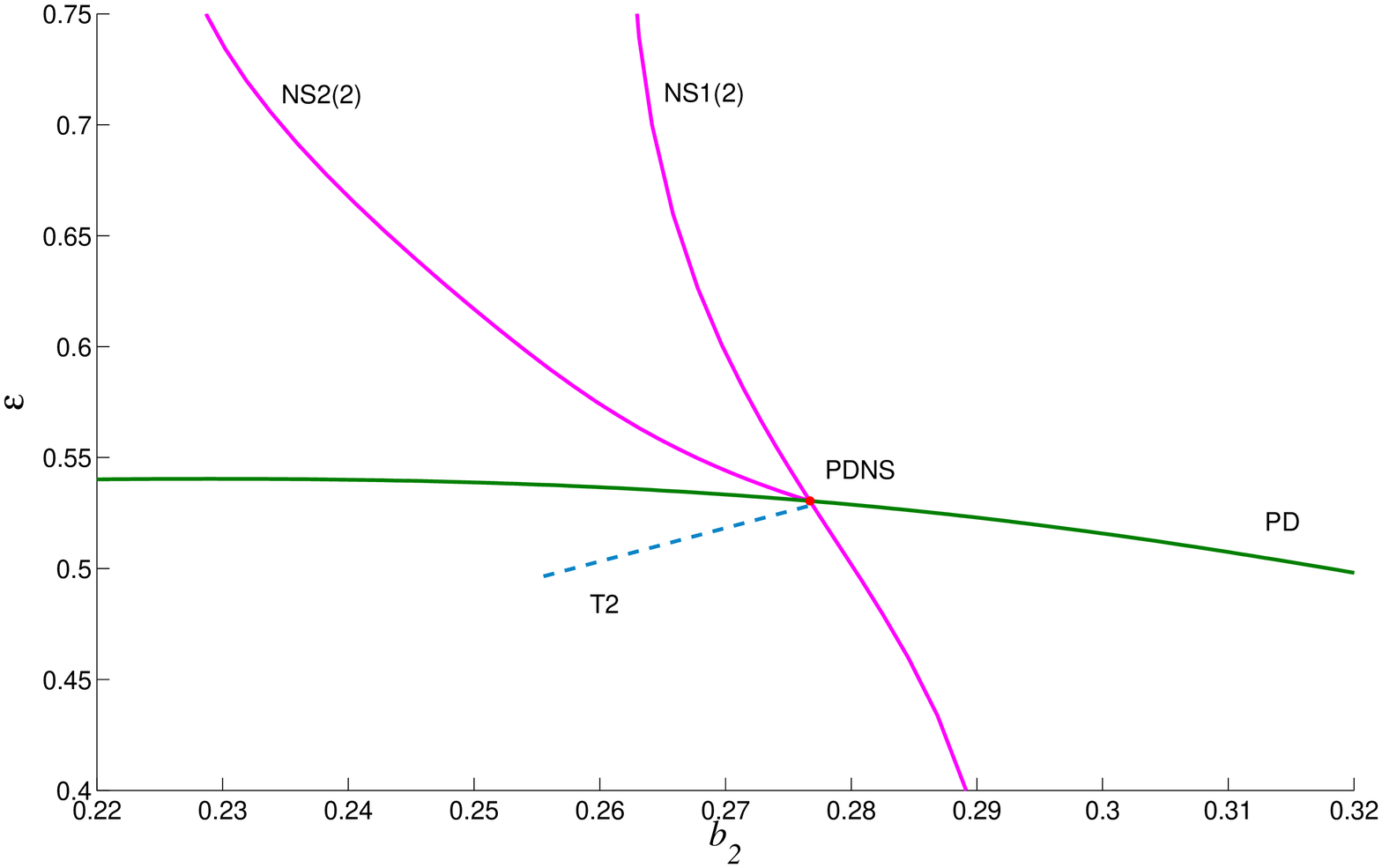}
 \caption{Zoom of the neighbourhood of the {\sf PDNS} point at $(b_2, \varepsilon)=(0.277,0.530)$ from 
 Figure \ref{fig:periodicforced}. In green are period doubling curves, in magenta Neimark-Sacker curves 
 (of the first or of the second iterate, respectively labeled with {\sf NS1(2)} and {\sf NS2(2)}), in blue is 
 the sketch of the {\sf T2} `curve'.} \label{fig:periodicforced_zoomrightpdns}
\end{figure}

\begin{figure}[htbp]
\centering
 \footnotesize
 \psfrag{b2}[][c]{$b_2$}
 \psfrag{eps}[][c]{$\varepsilon$}
  \includegraphics[width=.95\textwidth]{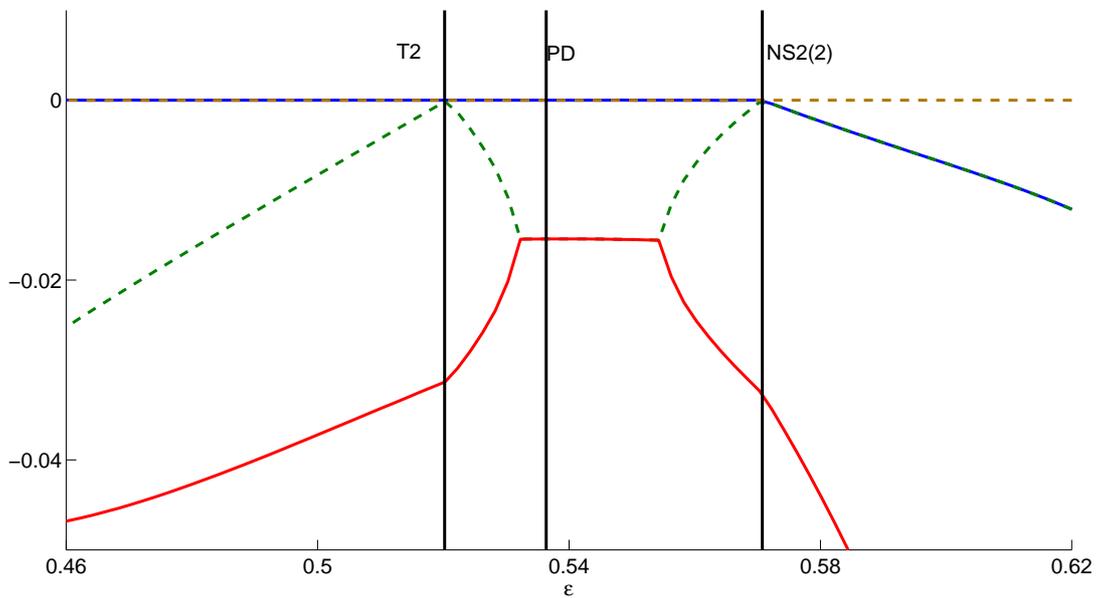}
 \caption{Lyapunov exponents computed for $b_2=0.261$, close to the PDNS point at 
 $(b_2, \varepsilon)=(0.277,0.530)$.} \label{fig:lyapunov_pdnsright}
\end{figure}

Figure \ref{fig:periodicforced} represents a bifurcation diagram for system (\ref{eq:PreyPredator})
where two {\sf PDNS} points are detected. The right {\sf PDNS} point has parameter values 
$(b_2, \varepsilon)=(0.277,0.530)$. We are in the ``simple" case of  Section \ref{Appendix:2PDNS} because 
the product of the coefficients $p_{11}=-5.01 \cdot 10^{-2}$ and $p_{22}=-0.211$ is positive. 
Since $\theta=-0.320$ and $\delta= 1.087$, Figure \ref{fig:NF_PDNSoverviewleft} indicates  
that the bifurcation diagram in a neighbourhood of the {\sf PDNS} point is as in case III in Figure 
\ref{fig:NF_PDNSsimple_1}, where $\mu_1=0$ corresponds with {\sf NS1} and $\mu_2=0$ with {\sf PD}. Curve $T_1$ 
corresponds to the Neimark-Sacker curve of the period doubled cycle {\sf NS2(2)} in 
Figure \ref{fig:periodicforced}. Therefore, we expect the period doubling `curve' {\sf T2} of the torus 
to be situated to the left of {\sf NS1(2)} and under the {\sf PD} curve. The stable limit cycles are situated 
in the lower right region of the {\sf PDNS} point. The exact location of {\sf T2} can be determined by 
computing Lyapunov exponents for fixed $b_2$ values smaller than the critical $b_2=0.277$ 
corresponding with the {\sf PDNS} point. We have plotted a sketch of this {\sf T2} curve in 
Figure \ref{fig:periodicforced_zoomrightpdns}, which represents a zoom of the neighbourhood 
of the {\sc PDNS} point and which includes a plot of {\sf NS2(2)} (curve $T_1$ in Figure 
\ref{fig:NF_PDNSsimple_1}). We have computed the Lyapunov exponents for $b_2$ fixed at $0.261$ 
and $\varepsilon \in [0.46;0.62]$, see Figure \ref{fig:lyapunov_pdnsright}. In this figure 
the black vertical lines indicate the position of the {\sf PD} and {\sf NS2(2)} curves. 
From the value of the Lyapunov exponents we derive that {\sf T2} is crossed for $\varepsilon \approx 0.52$. 
To the left of the {\sf T2} curve in Figure \ref{fig:lyapunov_pdnsright}, we have a stable torus, 
arisen through the supercritical Neimark-Sacker curve {\sf NS1(2)}, corresponding with region $2$ from 
Figure \ref{fig:NF_PDNSsimple_2}. Between the curves {\sf T2} and {\sf NS2(2)}, the $2$-torus arisen 
through {\sf T2} is attracting. These regions correspond with region $6$ (between {\sf T2} and {\sf PD}) 
and $5$ (between {\sf PD} and {\sf NS2(2)}) from Figure \ref{fig:NF_PDNSsimple_2}. When crossing the 
{\sf NS2(2)} curve, the $2$-torus disappears and the period doubled cycle becomes attracting. All this 
is in agreement with the fact that two Lyapunov exponents are equal to zero to the left of {\sf NS2(2)}, 
where afterwards only one zero Lyapunov exponent is left.\\

\begin{figure}[htbp]
\centering
 \footnotesize
 \psfrag{b2}[][c]{$b_2$}
 \psfrag{eps}[][c]{$\varepsilon$}
  \includegraphics[width=.95\textwidth]{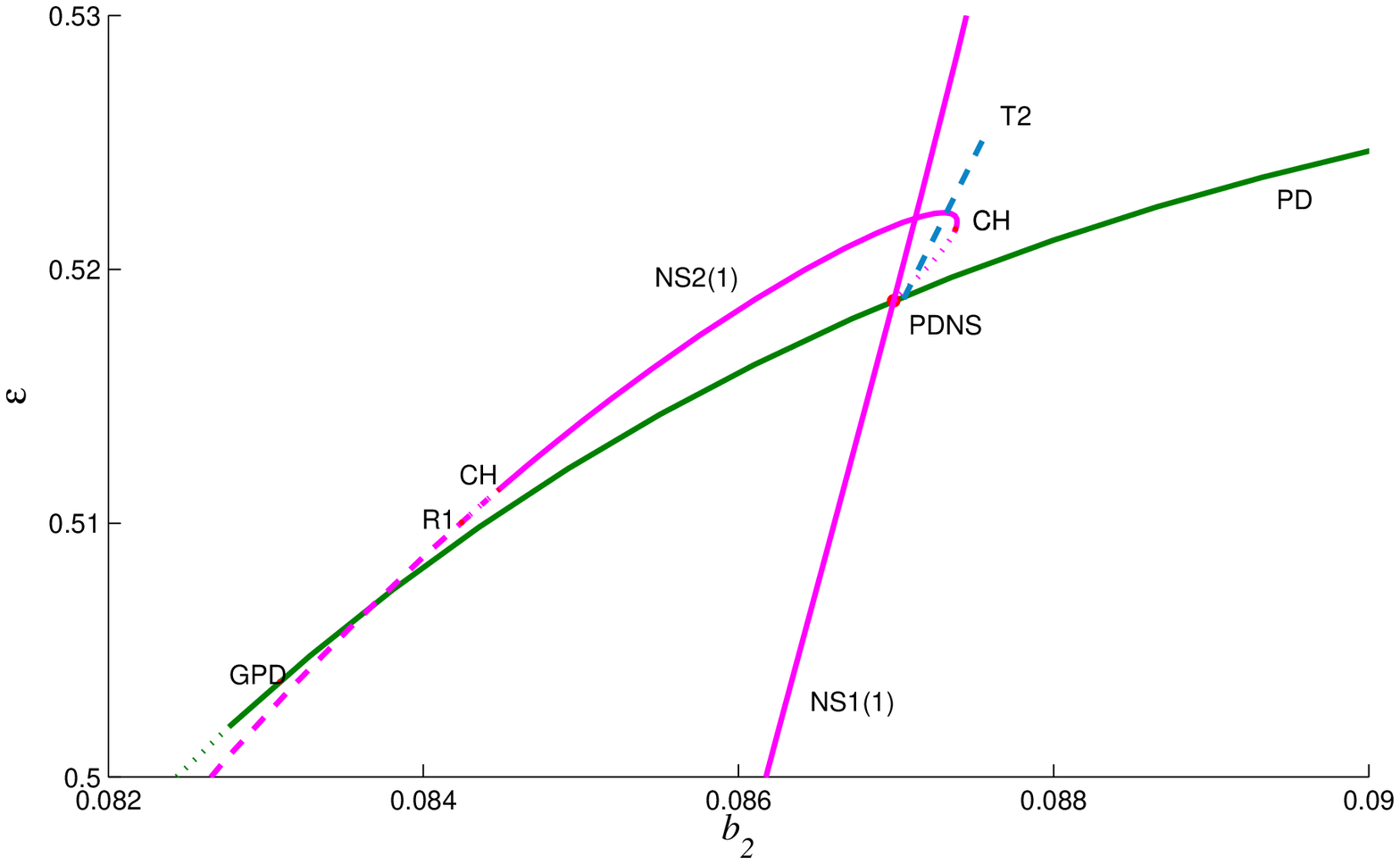}
 \caption{Zoom of the neighbourhood of the {\sf PDNS} point at 
 $(b_2, \varepsilon)=(8.699\cdot 10^{-2},0.519)$ from Figure \ref{fig:periodicforced}. 
 In green are period doubling curves, in magenta Neimark-Sacker curves (of the first 
 or of the second iterate, respectively labeled with {\sf NS1(1)} and {\sf NS2(1)}), in blue is the 
 sketch of the {\sf T2} `curve'. } \label{fig:periodicforced_zoomleftpdns}
\end{figure}

\begin{figure}[htb]
	\tiny
% 	\psfrag{b2}[][c]{$b_2$}
% 	\psfrag{eps}[][c]{$\epsilon$}
  \includegraphics[width=.95\textwidth]{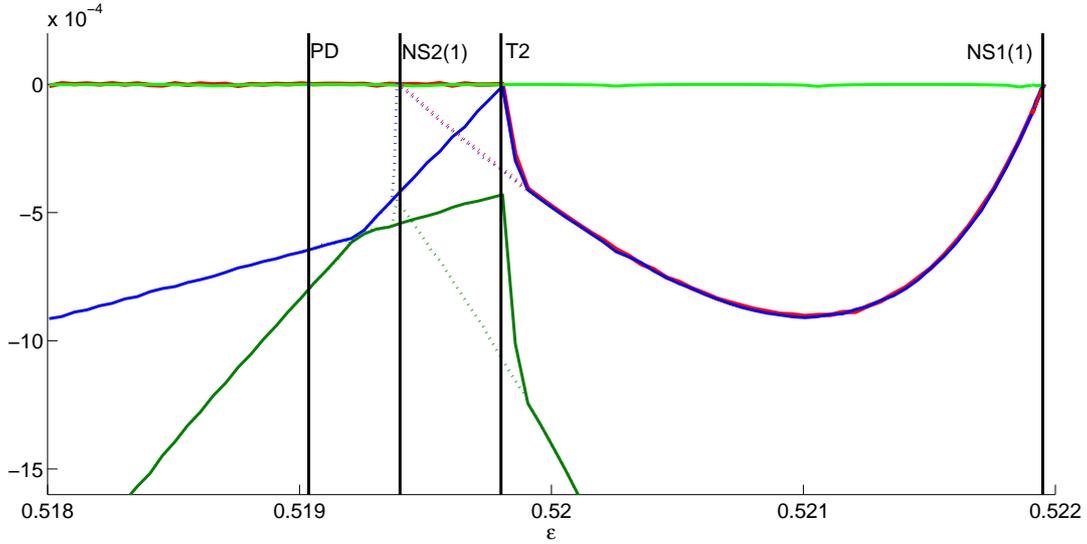}
%	\label{fig:lyapunov_pdnsleft(kip)}}
	\caption{Lyapunov exponents computed for $b_2=0.08709$, close to the {\sf PDNS} point at $(b_2, \varepsilon)=(8.699\cdot 10^{-2},0.519)$. Exponents indicated with solid lines are computed by following the attractor with increasing $\varepsilon$, dotted lines with decreasing $\varepsilon$. This highlights the bistability between {\sf NS2(1)} and {\sf T2}.}
	\label{fig:lyapunov_pdnsleft}
\end{figure}

The left {\sf PDNS} point at $(b_2, \varepsilon)=(8.699\cdot 10^{-2},0.519)$ again belongs to one of the 
``simple" situations in Section \ref{Appendix:2PDNS} ($p_{11}=-0.447,p_{22}=-1.472$). The neighbourhood 
of the bifurcation point is as in case I in Figure \ref{fig:NF_PDNSoverviewleft} 
since $(\theta,\delta)=(2.234,1.304)$. Remark that the stable limit cycles are situated 
in the lower left quadrant of the {\sf PDNS} point in Figure \ref{fig:periodicforced_zoomleftpdns}. 
The behaviour in a neighbourhood of this {\sf PDNS} point can be derived from Figure 
\ref{fig:periodicforced_zoomleftpdns}, which includes a plot of the Neimark-Sacker curve {\sf NS2(1)} 
of the period doubled cycle and also a sketch of the period doubled curve {\sf T2} of the torus, 
made on the basis of the computation of the Lyapunov exponents. We have calculated the Lyapunov 
exponents for parameter values in the upper right quadrant, close to the {\sf PDNS} point, for $b_2=0.08709$. 
The results are given in Figure \ref{fig:lyapunov_pdnsleft}. Going from the left to the right, where we 
follow the solid lines, we start with two Lyapunov exponents equal to zero which correspond with the 
stable torus from the original cycle in the regions $2, 3$ and $4$ from Figure \ref{fig:NF_PDNSsimple}. 
At the point where the second Lyapunov exponent becomes non-zero, the {\sf T2} curve is located, 
namely at $\varepsilon \approx 0.5198$. We then arrive in region $12$ from Figure \ref{fig:NF_PDNSsimple_2} 
where the $2$-torus has lost his stability and the period doubled cycle is stable. Therefore, one 
zero Lyapunov exponent remains. We scan the Lyapunov exponents for a second time where we now go 
from the right to the left and follow the dashed lines. The second Lyapunov exponent now approaches 
zero not at the {\sf T2} curve but at the {\sf NS2(1)} curve. This is explained by the bistability happening in 
region 4, where one Lyapunov exponent equal to zero indicates the stable period doubled cycle and 
two zero Lyapunov exponents indicate the stable torus. When going further, we cross region 3 and 2, 
with the stable torus of the orginal cycle. \\

Remark that since we have a periodically forced system the return time is independent of the 
distance from the limit cycle, so we could do this extra check. Indeed, for all {\sf PDNS} points, 
the $\alpha_{ijk}$ in the first equation of (\ref{eq:NF-PDNS}) are zero up to the accuracy of the 
computation. Here too, the Lyapunov exponents corroborate the prediction based on the normal form coefficients.

\subsection{Control of vibrations}
In \cite{Faver:03} a two-mass system of which the main mass is excited by a flow-induced, 
self excited force is studied. A single mass which acts as a dynamic absorber is attached to the 
main mass and, by varying the stiffness between the main mass and the absorber mass, represents a 
parametric excitation. The system is given by
\begin{equation} \label{eq:vibrationsmodel}
\begin{cases}
\dot x_1 = v_1\\
\dot x_2 = v_2\\
\dot v_1 = -k_1(v_1-v_2)-Q^2(1+\varepsilon y_1)(x_1-x_2)\\
\dot v_2 = Mk_1(v_1-v_2)+MQ^2(1+\varepsilon y_1)(x_1-x_2)-k_2v_2-x_2+\beta V^2(1-\gamma v_2^2)v_2\\
\dot y_1 = -\eta y_2+y_1(1-y_1^2-y_2^2)\\
\dot y_2 = \eta y_1+y_2(1-y_1^2-y_2^2).
\end{cases}
\end{equation}
The following parameters are fixed: $\varepsilon=0.1, k_2=0.1, \beta=0.1, V=\sqrt{2.1}, 
\gamma=4, Q=0.95, M=0.2$, $k_1$ and $\eta$ will be the continuation parameters. 

\subsubsection{The NSNS points}
\begin{figure}[htbp]
\centering
 \footnotesize
% \psfrag{k1}[][c]{$k_1$}
% \psfrag{eta}[][c]{$\eta$}
  \includegraphics[width=.7\textwidth]{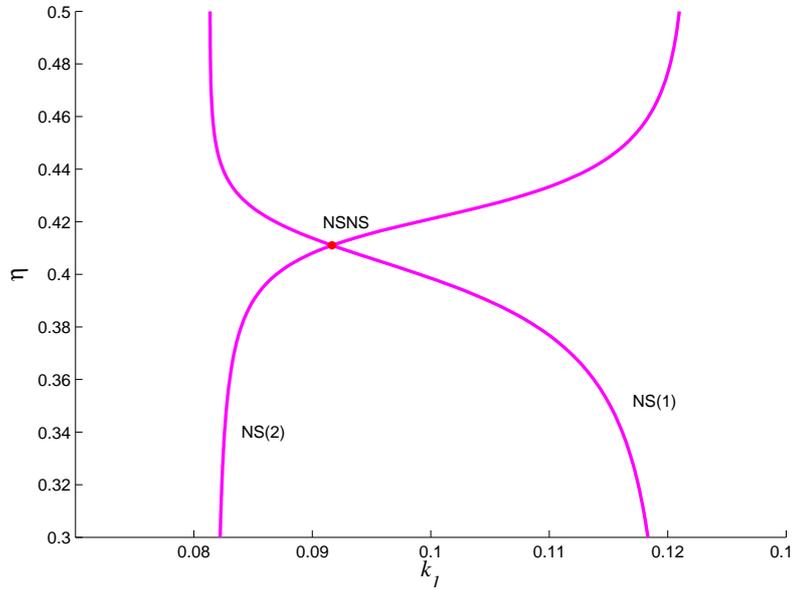}
 \caption{Partial bifurcation diagram of limit cycles in system \eqref{eq:vibrationsmodel}.} \label{fig:faver}
\end{figure}

\begin{figure}[htbp]
\centering
 \footnotesize
 \psfrag{k1}[][c]{$k_1$}
 \psfrag{eta}[][c]{$\eta$}
  \includegraphics[width=.85\textwidth]{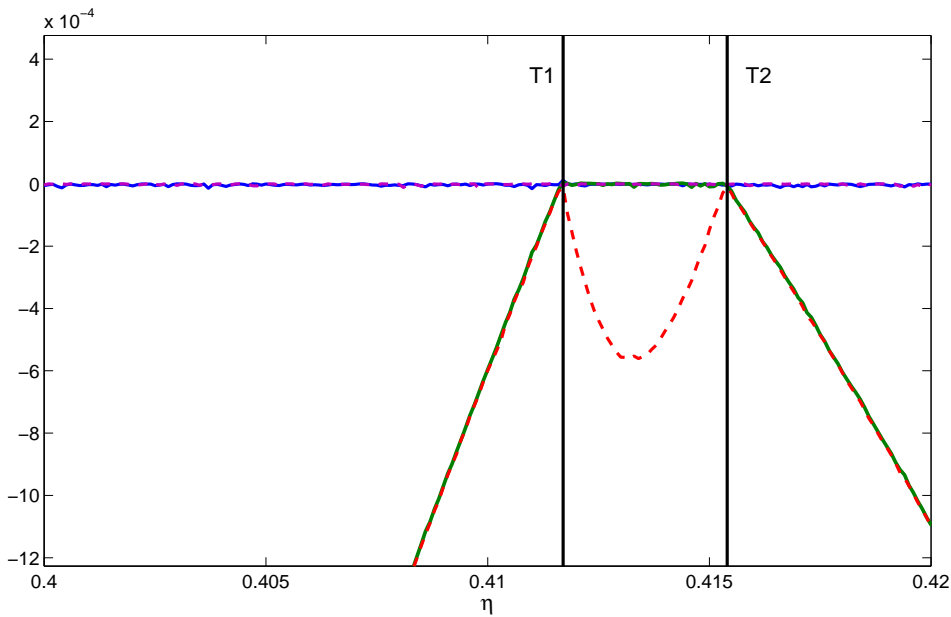}
 \caption{Lyapunov exponents computed for $k_1=0.083$.} \label{fig:lyapunov_nsns}
\end{figure}

An {\sf NSNS} point is detected for $(k_1,\eta) = (9.167\cdot 10^{-2},0.411)$, see Figure \ref{fig:faver}.
The normal form coefficients are 
$$
(p_{11},p_{22},\theta,\delta,\mbox{sign }l_1) = 
(-3.733 \cdot 10^{-3}, -6.494 \cdot 10^{-3}, 0.541, 1.203, 1).
$$ 
The positive sign of the 
product $p_{11}p_{22}$ implies that we are in a ``simple" case of Section \ref{Appendix:2NSNS}. 
Since $\delta > \theta$, the role of both coefficients has to be reversed. 
Therefore, $\theta>1, \delta<1, \theta \delta<1$ indicate that the {\sf NSNS} bifurcation is 
located in region II in Figure \ref{fig:NF_PDNSsimple_1}. As in the previous examples, we 
have computed the Lyapunov exponents to check the obtained results of the normal form 
coefficients. We have done the computations for $k_1$ fixed at $0.083$ and $\eta 
\in [0.4;0.42]$ ($\eta$ values are between the {\sf NS} curves). The results are given in 
Figure \ref{fig:lyapunov_nsns}. For $\eta$-values starting from $0.38$, we are in region 
$3$ (or $12$ due to symmetry) in Figure \ref{fig:NF_PDNSsimple_2}, where there is a 
stable 2-torus and thus two Lyapunov exponents equal to zero. A third Lyapunov exponent 
approaches zero and between $\eta\approx 0.4117$ and $\eta\approx 0.4154$ three Lyapunov 
exponents are equal to zero. This region denotes the appearance of a stable $3$-torus and 
corresponds with region $5$ from Figure \ref{fig:NF_PDNSsimple_2} II. The critical values 
of $\eta$ correspond with the curves $T_1$ and $T_2$ in Figure \ref{fig:NF_PDNSsimple_1}. 
For $\eta \geq 0.4154$, only a stable $2$-torus remains such that there are two zero 
Lyapunov exponents. Therefore, the computed Lyapunov exponents are in agreement with the 
normal form coefficients.

Also in this case all $\alpha_{ijkl}$ in the normal form (\ref{eq:NF-NSNS}) vanish since we have a 
periodically forced system.

\section{Discussion}\label{Section:Discussion}

This paper completes the development of efficient methods for the computation of the critical normal form coefficients for all 
codim 1 and 2 local bifurcations of limit cycles, started in \cite{KuDoGoDh:05,drdwgk}  and based on \cite{Io:88}. Together with 
our previous papers on the computation of the critical normal form coeffcients for codim 1 and 2 local bifurcations of equilibria in ODEs 
\cite{Ku:99} and fixed points of maps \cite{KuMe:2005,KuMe:2006},  it contributes to the development of methods, algorithms, and 
software tools for multiparameter bifurcation analysis of smooth finite-dimensional dynamical systems. 

The resulting formulas are independent of the phase space dimension and can be applied in the original basis, without preliminary linear transformations. 
As limit cycles are concerned, the formulas are directly suitable for numerical implementation using orthogonal collocation. They 
fit perfectly into a continuation context, where limit cycles and their bifurcations are computed using the BVP-approach
\cite{DoGoKu:2003}, without numerical approximation of the Poincar\'{e} map or its derivatives. Being implemented into the 
{\sc matlab} toolbox {\sc matcont} \cite{MATCONT,NewMATCONT}, the developed methods are freely available to assist an advanced two-parameter
bifurcation analysis of dynamical systems generated by ODEs and maps from various applications. 

To fully support the two-parameter bifurcation analysis of ODEs and maps, one needs special methods to switch between various 
branches of codim 1 bifurcations of fixed points and cycles rooted at codim 2 points. Such methods have been developed and 
implemented in {\sc matcont} for codim 2 equilibrium \cite{KuMeGoSa:2008} and fixed point \cite{GoGhKuMe:07} bifurcations. 
Switching at codim 2 points to the continuation of codim 1 local bifurcations of limit cycles seems to be the next natural problem to attack,
while that for codim 1 bifurcations of homoclinic and heteroclinic orbits is more difficult and probably requires new ideas. Similar remarks 
can be made about quasiperiodic bifurcations of tori.

\appendix
\section{Derivation of the critical normal forms}
\label{Appendix:1}

\subsection{Notation}\label{sectionNotation}
Let $M \in \R^{n \times n}$ be the monodromy matrix. In all codimension $2$ cases all critical multipliers, i.e., all multipliers with modulus $1$, have non-degenerate Jordan blocks. Let $M_0$ be the critical Jordan structure, i.e., the block diagonal matrix consisting of the critical Jordan blocks, starting with the block of the trivial multiplier $1$. Let $\mu_k = e^{i\theta_k} (0\leq \theta_k < 2\pi)$ be a critical multiplier with multiplicity $m_k$. The matrix $L_k \in \R^{m_k \times m_k}$ is
defined as
\[
L_k=\begin{pmatrix}
 \sigma_k & 1  & \ldots & 0 \\
 0 & \sigma_k  & \ldots & 0 \\
 \vdots  & \ddots &\ddots & 1\\
 0  & \ldots & 0 & \sigma_k
\end{pmatrix},
\]
where $\sigma_k$ is the \textit{Floquet
exponent of the multiplier} $\mu_k$, with $\sigma_k = i\theta_k/T$ in the case of a positive real multiplier or a complex multiplier $\mu_k$ and $\sigma_k = 0$ for $\mu_k=-1$. The matrix $L_0$ is the block diagonal matrix formed from
the blocks $L_k$ for which $\left|\mu_k\right|=1$, starting with the block that corresponds with multiplier $1$. The matrix $\tilde L_0$ is the
matrix $L_0$ without the first row and the first column.

\subsubsection{\tt LPNS} \label{Appendix:1_lpns} At the LPNS bifurcation the matrices described
in \ref{sectionNotation} are
\[
M_0=\begin{pmatrix} 1 & 1 & 0 & 0 \\ 0 & 1 & 0 & 0 \\ 0 & 0 & e^{i\omega T} & 0 \\ 0 & 0 & 0 & e^{-i\omega T}
\end{pmatrix}, \quad
L_0=\begin{pmatrix} 0 & 1 & 0 & 0 \\ 0 & 0 & 0 & 0 \\ 0 & 0 & i\omega & 0 \\ 0 & 0 & 0 & -i\omega
\end{pmatrix}, \quad
\tilde L_0 =\begin{pmatrix} 0 & 0 & 0 \\ 0 & i\omega & 0 \\ 0 & 0 & -i\omega
\end{pmatrix}.
\]
We are in a case in which we can apply Theorem 2 from
\cite{Io:88}. So we can define a $T$-periodic normal form
\[
\ds\dd{\tau}{t}=1+\xi_1+p(\tau,\xi), \quad \ds\dd{\xi}{\tau}=\tilde L_0
\xi + P(\tau,\xi),
\]
where $\xi=(\xi_1,\xi_2,\bar \xi_2)^{\rm T}$ with $\xi_1 \in \R, \xi_2\in \C$. The polynomials $p$ and $P$ are real, respectively complex, $T$-periodic in $\tau$ and
at least quadratic in $(\xi_1,\xi_2,\bar\xi_2)$ such that
\begin{eqnarray*}
\ds\dd{}{\tau}p(\tau,\xi)-\ds\dd{}{\xi}p(\tau,\xi)\tilde L_0^*\xi=0, \quad
\ds\dd{}{\tau}P(\tau,\xi)+\tilde L_0^* P(\tau,\xi)-\ds\dd{}{\xi}P(\tau,\xi)\tilde L_0^*\xi=0.
\end{eqnarray*}
If we write the polynomials in a Fourier expansion, namely
\begin{equation*}
p(\tau,\xi)= \sum_{l=-\infty}^{\infty} p_l(\xi) e^{i \frac{2 \pi l
\tau}{T}}, \quad P(\tau,\xi)= \sum_{l=-\infty}^{\infty} P_l(\xi) e^{i
\frac{2 \pi l \tau}{T}},
\end{equation*}
we obtain for any $l\in \mathbb Z$ the following differential equations
\begin{eqnarray*}
 &\dd{}{\xi}p_l(\xi) \tilde L_0^* \xi-i\frac{2 \pi l}{T} p_l(\xi)=0,& \\
 &\dd{}{\xi}P_l(\xi) \tilde L_0^* \xi-i\frac{2 \pi l}{T} P_l(\xi)- \tilde L_0^*
 P_l(\xi)=0.&
\end{eqnarray*}
Putting $\tilde L_0$ into the equations and writing
$$
P_l(\xi_1,\xi_2,\bar\xi_2)=(P^{(1)}_l(\xi_1,\xi_2,\bar\xi_2),P^{(2)}_l(\xi_1,\xi_2,\bar\xi_2),\bar P^{(2)}_l(\xi_1,\xi_2,\bar\xi_2))^{\rm T}
$$
we can rewrite them as a set of differential equations in variable $\xi_2$
\begin{eqnarray*}
 i\omega \xi_2\dd{}{\xi_2}p_l(\xi_1,\xi_2,\bar\xi_2) +i \frac{2 \pi l}{T} p_l(\xi_1,\xi_2,\bar \xi_2)& =& i\omega \bar \xi_2\dd{}{\bar \xi_2}p_l(\xi_1,\xi_2,\bar\xi_2), \\
 i\omega \xi_2\dd{}{\xi_2}P^{(1)}_l(\xi_1,\xi_2,\bar\xi_2) + i \frac{2 \pi l}{T} P^{(1)}_l(\xi_1,\xi_2,\bar \xi_2)& =& i\omega \bar\xi_2\dd{}{\bar\xi_2}P^{(1)}_l(\xi_1,\xi_2,\bar\xi_2), \\
 i\omega \xi_2\dd{}{\xi_2}P^{(2)}_l(\xi_1,\xi_2,\bar\xi_2) + i \frac{2 \pi l}{T} P^{(2)}_l(\xi_1,\xi_2,\bar \xi_2)& =& i\omega \bar\xi_2\dd{}{\bar\xi_2}P^{(2)}_l(\xi_1,\xi_2,\bar\xi_2)+i \omega P^{(2)}_l(\xi_1,\xi_2,\bar\xi_2),\\
 i\omega \xi_2\dd{}{\xi_2}\bar P^{(2)}_l(\xi_1,\xi_2,\bar\xi_2) + i \frac{2 \pi l}{T} \bar P^{(2)}_l(\xi_1,\xi_2,\bar \xi_2)& =& i\omega \bar\xi_2\dd{}{\bar\xi_2}\bar P^{(2)}_l(\xi_1,\xi_2,\bar\xi_2)-i \omega \bar P^{(2)}_l(\xi_1,\xi_2,\bar\xi_2).
\end{eqnarray*}
$p_l(\xi_1,\xi_2,\bar \xi_2),P^{(1)}_l(\xi_1,\xi_2,\bar \xi_2)$ and $P^{(2)}_l(\xi_1,\xi_2,\bar \xi_2)$
are polynomials, and from the equations it follows that $p_l(\xi_1,\xi_2,\bar \xi_2), P^{(1)}_l(\xi_1,\xi_2,\bar \xi_2)$ and $P^{(2)}_l(\xi_1,\xi_2,\bar \xi_2)$ are zero if $l\neq 0$. Thus the polynomials are $\tau$-independent.
We obtain
\begin{eqnarray*}
 \xi_2\dd{}{\xi_2}p_0(\xi_1,\xi_2,\bar\xi_2)& =&\bar \xi_2\dd{}{\bar \xi_2}p_0(\xi_1,\xi_2,\bar\xi_2), \\
 \xi_2\dd{}{\xi_2}P^{(1)}_0(\xi_1,\xi_2,\bar\xi_2)& =& \bar\xi_2\dd{}{\bar\xi_2}P^{(1)}_0(\xi_1,\xi_2,\bar\xi_2), \\
 \xi_2\dd{}{\xi_2}P^{(2)}_0(\xi_1,\xi_2,\bar\xi_2)& =& \bar\xi_2\dd{}{\bar\xi_2}P^{(2)}_0(\xi_1,\xi_2,\bar\xi_2)+ P^{(2)}_0(\xi_1,\xi_2,\bar\xi_2),
\end{eqnarray*}
and the complex conjugate of the last equation. From the first equation it follows that 
\begin{eqnarray*}
p_0(\xi_1,\xi_2,\bar\xi_2) = \psi_1(\xi_1)+\psi_2(\left|\xi_2\right|^2)+\psi_3(\xi_1)\psi_4(\left|\xi_2\right|^2),
\end{eqnarray*}
where $\psi_2,\psi_3$ and $\psi_4$ are at least linear in their argument and $\psi_1$ at least quadratic. Similarly, we obtain 
\begin{eqnarray*}
P^{(1)}_0(\xi_1,\xi_2,\bar\xi_2) = \phi_1(\xi_1)+\phi_2(\left|\xi_2\right|^2)+\phi_3(\xi_1)\phi_4(\left|\xi_2\right|^2),
\end{eqnarray*}
with the same conditions for $\phi$ as the ones for $\psi$. At last, from the third equation we can derive that
\begin{eqnarray*}
P^{(2)}_0(\xi_1,\xi_2,\bar\xi_2) = \xi_2\chi_1(\xi_1)+\xi_2\chi_2(\left|\xi_2\right|^2)+\xi_2\chi_3(\xi_1)\chi_4(\left|\xi_2\right|^2),
\end{eqnarray*}
where 
$\chi_1,\chi_2,\chi_3$ and $\chi_4$ are at least linear in their argument.

Assembling all the information gives us the following normal form
\[\begin{cases}
 \ds\dd{\tau}{t}=1+\xi_1+\alpha_{200} \xi_1^2+\alpha_{011}\left|\xi_2\right|^2+\alpha_{300}' \xi_1^3+\alpha_{111}' \xi_1\left|\xi_2\right|^2+\ldots, \smallskip \\
 \ds\dd{\xi_1}{\tau}= a_{200}' \xi_1^2+a_{011}'\left|\xi_2\right|^2+a_{300} \xi_1^3+a_{111} \xi_1\left|\xi_2\right|^2+\ldots, \smallskip \\
 \ds\dd{\xi_2}{\tau}= i\omega \xi_2+b_{110}' \xi_1\xi_2+b_{210} \xi_1^2 \xi_2+b_{021} \xi_2\left|\xi_2\right|^2 +\ldots.
\end{cases}
\]
We do the substitution $\xi_1 \mapsto -\xi_1$ and find the Iooss normal form (\ref{eq:NF-LPNS}), i.e.
\begin{equation*}
\begin{cases}
 \ds\dd{\tau}{t}=1-\xi_1+\alpha_{200} \xi_1^2+\alpha_{011}\left|\xi_2\right|^2+\alpha_{300} \xi_1^3+\alpha_{111} \xi_1\left|\xi_2\right|^2+\ldots, \smallskip \\
 \ds\dd{\xi_1}{\tau}= a_{200} \xi_1^2+a_{011}\left|\xi_2\right|^2+a_{300} \xi_1^3+a_{111} \xi_1\left|\xi_2\right|^2+\ldots, \smallskip \\
 \ds\dd{\xi_2}{\tau}= i\omega \xi_2+b_{110} \xi_1\xi_2+b_{210} \xi_1^2 \xi_2+b_{021} \xi_2\left|\xi_2\right|^2 +\ldots,
\end{cases}
\end{equation*}
where the dots denote $O(\left|\xi\right|^4)$ terms. Note that the time evolution can be obtained by applying the chain 
rule to this system and is then given by
\[
\begin{cases}
\ds\dd{\tau}{t} =1-\xi_1+\alpha_{200} \xi_1^2+\alpha_{011}\left|\xi_2\right|^2+\alpha_{300} \xi_1^3+\alpha_{111} \xi_1\left|\xi_2\right|^2+\ldots,\smallskip\\
\ds\dd{\xi_1}{t}  = a_{200}\xi_1^2+a_{011}\left|\xi_2\right|^2+a_{300}'\xi_1^3+a_{111}'\xi_1\left|\xi_2\right|^2+\ldots,\smallskip\\
\ds\dd{\xi_2}{t}  = i\omega \xi_2+b_{110}'\xi_1\xi_2+b_{210}'\xi_1^2\xi_2+b_{021}'\xi_2\left|\xi_2\right|^2+\ldots,
\end{cases}
\]
with $a_{300}'=-a_{200}+a_{300}$ and $a_{111}'=-a_{011}+a_{111}$, and 
with $b_{110}'=-i\omega+b_{110},b_{210}'=i\omega \alpha_{200}-b_{110}+b_{210}$ and $b_{021}'=i\omega \alpha_{011}+b_{021}$. We could use this system as our starting normal form. To draw conclusions about the bifurcation diagrams, we could then perform the time reparametrization on the center manifold to obtain an autonomous truncated ODE which approximates the Poincar\'e map as a $T$-shift (as done in \cite{drdwgk}). It could be studied by comparing the obtained ODE to the one for the Zero-Hopf bifurcation of equilibria. However, the time reparametrized ODE has exactly the same form as (\ref{eq:NF-LPNS}), so we can as well use (\ref{eq:NF-LPNS}). The approach we will follow in this paper (with (\ref{eq:NF-LPNS}) as starting normal form) is mathematically equivalent to the one of \cite{drdwgk}, but takes a shorter path.

\bigskip

\subsubsection{\tt PDNS} \label{Appendix:1_pdns} At the PDNS bifurcation the matrices described
in \ref{sectionNotation} are
\[
M_0=\begin{pmatrix} 1 & 0 & 0 & 0 \\ 0 & -1 & 0 & 0 \\ 0 & 0 & e^{i\omega T} & 0 \\ 0 & 0 & 0 & e^{-i\omega T}
\end{pmatrix}, 
L_0=\begin{pmatrix} 0 & 0 & 0 & 0 \\ 0 & 0 & 0 & 0 \\ 0 & 0 & i\omega & 0 \\ 0 & 0 & 0 & -i\omega
\end{pmatrix}, 
\tilde L_0 =\begin{pmatrix} 0 & 0 & 0 \\ 0 & i\omega & 0 \\ 0 & 0 & -i\omega
\end{pmatrix}.
\]
We are in a case in which we can apply Theorem 3 from
\cite{Io:88}. So we can define a $2T$-periodic normal form
\[
\ds\dd{\tau}{t}=1+p(\tau,\xi), \quad \ds\dd{\xi}{\tau}=\tilde L_0
\xi + P(\tau,\xi),
\]
where $\xi=(\xi_1,\xi_2,\bar \xi_2)$. The polynomials $p$ and $P$ are $2T$-periodic in $\tau$ and
at least quadratic in their argument such that
\begin{eqnarray*}\label{eq:PDNS_proprieties1}\\
\ds\dd{}{\tau}p(\tau,\xi)-\ds\dd{}{\xi}p(\tau,\xi)\tilde L_0^*\xi &=&0, \\
\ds\dd{}{\tau}P(\tau,\xi)+\tilde L_0^* P(\tau,\xi)-\ds\dd{}{\xi}P(\tau,\xi)\tilde L_0^*\xi &=&0, \\ 
p(\tau+T,\xi_1,\xi_2,\bar \xi_2) &=& p(\tau,-\xi_1,\xi_2,\bar \xi_2), \\
P^{(1)}(\tau+T,-\xi_1,\xi_2,\bar \xi_2) &=& -P^{(1)}(\tau,\xi_1,\xi_2,\bar \xi_2), \\
P^{(2)}(\tau+T,-\xi_1,\xi_2,\bar \xi_2) &=& P^{(2)}(\tau,\xi_1,\xi_2,\bar \xi_2),
\end{eqnarray*}
and the complex conjugate of the last equation.

As in the LPNS case (since the $\tilde L_0$ matrix is the same) we obtain that all
polynomials are independent from $\tau$, and thus we can rewrite the last three equations  as
\begin{eqnarray*}
 &&p(\xi_1,\xi_2,\bar \xi_2)=p(-\xi_1,\xi_2,\bar \xi_2),\\
 &&P^{(1)}(-\xi_1,\xi_2,\bar \xi_2)=-P^{(1)}(\xi_1,\xi_2,\bar \xi_2), \quad P^{(2)}(-\xi_1,\xi_2,\bar \xi_2)=P^{(2)}(\xi_1,\xi_2,\bar \xi_2),
\end{eqnarray*}
thus $p$ and $P^{(2)}$ are even in $\xi_1$ and $P^{(1)}$ is odd in $\xi_1$. Taking the results from the LPNS case into account, we obtain
\begin{eqnarray*}
&&p_0(\xi_1,\xi_2,\bar\xi_2) = \psi_1(\xi_1^2)+\psi_2(\left|\xi_2\right|^2)+\psi_3(\xi_1^2)\psi_4(\left|\xi_2\right|^2),\\
&&P^{(1)}_0(\xi_1,\xi_2,\bar\xi_2) = \xi_1\phi_1(\xi_1^2)+\xi_1\phi_2(\left|\xi_2\right|^2)+\xi_1\phi_3(\xi_1^2)\phi_4(\left|\xi_2\right|^2),\\
&&P^{(2)}_0(\xi_1,\xi_2,\bar\xi_2) = \xi_2\chi_1(\xi_1^2)+\xi_2\chi_2(\left|\xi_2\right|^2)+\xi_2\chi_3(\xi_1^2)\chi_4(\left|\xi_2\right|^2),
\end{eqnarray*}
with all functions at least linear in their argument. 

Assembling all the information gives us the Iooss normal form (\ref{eq:NF-PDNS}), i.e.
\begin{equation*}
\begin{cases}
 \ds\dd{\tau}{t}=1+\alpha_{200} \xi_1^2+\alpha_{011}\left|\xi_2\right|^2+\alpha_{400} \xi_1^4+\alpha_{022}\left|\xi_2\right|^4+\alpha_{211} \xi_1^2\left|\xi_2\right|^2+\ldots, \smallskip \\
 \ds\dd{\xi_1}{\tau}= a_{300} \xi_1^3+a_{111}\xi_1\left|\xi_2\right|^2+a_{500} \xi_1^5+a_{122} \xi_1\left|\xi_2\right|^4+a_{311} \xi_1^3\left|\xi_2\right|^2+\ldots, \smallskip \\
 \ds\dd{\xi_2}{\tau}= i\omega \xi_2+b_{210} \xi_1^2\xi_2+b_{021} \xi_2\left|\xi_2\right|^2+b_{410} \xi_1^4\xi_2+b_{221} \xi_1^2\xi_2\left|\xi_2\right|^2+b_{032} \xi_2\left|\xi_2\right|^4+\ldots,
\end{cases}
\label{normalformpdns}
\end{equation*}
where the dots denote $O(\left|\xi\right|^6)$ terms. Note that the time evolution can be obtained by applying the chain rule to this
system and is given by
\begin{equation*}
\begin{cases}
 \dd{\tau}{t}=1+\alpha_{200} \xi_1^2+\alpha_{011}\left|\xi_2\right|^2+\alpha_{400} \xi_1^4+\alpha_{022}\left|\xi_2\right|^4+\alpha_{211} \xi_1^2\left|\xi_2\right|^2 + \ldots, \smallskip \\
 \dd{\xi_1}{t}= a_{300}\xi_1^3+a_{111}\xi_1\left|\xi_2\right|^2+a_{500}'\xi_1^5+a_{311}'\xi_1^3\left|\xi_2\right|^2+a_{122}'\xi_1\left|\xi_2\right|^4+ \ldots, \smallskip \\
 \dd{\xi_2}{t} = i\omega \xi_2+b_{210}'\xi_1^2\xi_2+b_{021}'\xi_2\left|\xi_2\right|^2+b_{410}'\xi_1^4\xi_2+b_{221}'\xi_1^2\xi_2\left|\xi_2\right|^2+b_{032}'\xi_2\left|\xi_2\right|^4+\ldots, \smallskip\\
\end{cases}
\end{equation*}
with $a_{500}'=a_{300}\alpha_{200}+a_{500}, a_{311}'=a_{300}\alpha_{011}+a_{111}\alpha_{200}+a_{311}$, $a_{122}'=a_{111}\alpha_{011}+a_{122}$, $b_{210}'=i\omega \alpha_{200}+b_{210},b_{021}'=i\omega \alpha_{011}+b_{021}, b_{410}'=i\omega \alpha_{400}+b_{210}\alpha_{200}+b_{410}, b_{221}'=i\omega \alpha_{211}+b_{210}\alpha_{011}+b_{021}\alpha_{200}+b_{221}$ and $b_{032}'=i\omega \alpha_{022}+b_{021}\alpha_{011}+b_{032}$.\\

\subsubsection{NSNS} \label{Appendix:1_nsns} 
At the NSNS bifurcation the matrices described
in \ref{sectionNotation} are
\[
M_0=\begin{pmatrix} 1 & 0 & 0 & 0 & 0 \\ 0 & e^{i\omega_1 T} & 0 & 0 & 0 \\ 0 & 0 & e^{-i\omega_1 T} & 0 & 0 \\ 0 & 0 & 0 & e^{i\omega_2 T} & 0\\0 & 0 & 0 & 0 & e^{-i\omega_2 T}
\end{pmatrix},\]
\[
L_0=\begin{pmatrix}  0 & 0 & 0 & 0 & 0 \\ 0 & i\omega_1 & 0 & 0 & 0 \\ 0 & 0 & -i\omega_1 & 0 & 0 \\ 0 & 0 & 0 & i\omega_2& 0\\0 & 0 & 0 & 0 & -i\omega_2
\end{pmatrix}, 
\tilde L_0 =\begin{pmatrix} i\omega_1& 0 & 0 & 0 \\ 0 & -i\omega_1& 0 & 0 \\ 0 & 0 & i\omega_2& 0\\0 &0&0& -i\omega_2
\end{pmatrix}.
\]
We are in a case in which we can apply Theorem 1 from
\cite{Io:88}. So we can define a $T$-periodic normal form
\[
\ds\dd{\tau}{t}=1+p(\tau,\xi), \quad \ds\dd{\xi}{\tau}=\tilde L_0
\xi + P(\tau,\xi),
\]
where $\xi=(\xi_1,\bar \xi_1,\xi_2,\bar \xi_2)$. The polynomials $p$ and $P$ are $T$-periodic in $\tau$ and
at least quadratic in $(\xi_1,\bar \xi_1,\xi_2,\bar\xi_2)$ such that
\begin{eqnarray*}
\ds\dd{}{\tau}p(\tau,\xi)-\ds\dd{}{\xi}p(\tau,\xi)\tilde L_0^*\xi=0, \quad
\ds\dd{}{\tau}P(\tau,\xi)+\tilde L_0^* P(\tau,\xi)-\ds\dd{}{\xi}P(\tau,\xi)\tilde L_0^*\xi=0.
\end{eqnarray*}
Writing down the polynomials in a Fourier expansion results in the following equations
\begin{eqnarray*}
 &&i\omega_1 \xi_1\dd{}{\xi_1}p_l(\xi_1,\bar \xi_1,\xi_2,\bar\xi_2)+i\omega_2 \xi_2\dd{}{\xi_2}p_l(\xi_1,\bar \xi_1,\xi_2,\bar\xi_2) +i \frac{2 \pi l}{T} p_l(\xi_1,\bar \xi_1,\xi_2,\bar \xi_2)\\
 &&= i \omega_1 \bar \xi_1\dd{}{\bar \xi_1}p_l(\xi_1,\bar \xi_1,\xi_2,\bar\xi_2)+i \omega_2 \bar \xi_2\dd{}{\bar \xi_2}p_l(\xi_1,\bar \xi_1,\xi_2,\bar\xi_2), \\
 &&i\omega_1 \xi_1\dd{}{\xi_1}P^{(1)}_l(\xi_1,\bar \xi_1,\xi_2,\bar\xi_2)+i\omega_2 \xi_2\dd{}{\xi_2}P^{(1)}_l(\xi_1,\bar \xi_1,\xi_2,\bar\xi_2) +i \frac{2 \pi l}{T} P^{(1)}_l(\xi_1,\bar \xi_1,\xi_2,\bar \xi_2)\\
 && = i \omega_1 \bar \xi_1\dd{}{\bar \xi_1}P^{(1)}_l(\xi_1,\bar \xi_1,\xi_2,\bar\xi_2)+i \omega_2 \bar \xi_2\dd{}{\bar \xi_2}P^{(1)}_l(\xi_1,\bar \xi_1,\xi_2,\bar\xi_2)+i\omega_1P^{(1)}_l(\xi_1,\bar \xi_1,\xi_2,\bar \xi_2), \\
 &&i\omega_1 \xi_1\dd{}{\xi_1}\bar P^{(1)}_l(\xi_1,\bar \xi_1,\xi_2,\bar\xi_2)+i\omega_2 \xi_2\dd{}{\xi_2}\bar P^{(1)}_l(\xi_1,\bar \xi_1,\xi_2,\bar\xi_2)+i\omega_1 \bar P^{(1)}_l(\xi_1,\bar \xi_1,\xi_2,\bar \xi_2) \\
 && = i \omega_1 \bar \xi_1\dd{}{\bar \xi_1}\bar P^{(1)}_l(\xi_1,\bar \xi_1,\xi_2,\bar\xi_2)+i \omega_2 \bar \xi_2\dd{}{\bar \xi_2}\bar P^{(1)}_l(\xi_1,\bar \xi_1,\xi_2,\bar\xi_2)-i \frac{2 \pi l}{T} \bar P^{(1)}_l(\xi_1,\bar \xi_1,\xi_2,\bar \xi_2), \\
  &&i\omega_1 \xi_1\dd{}{\xi_1}P^{(2)}_l(\xi_1,\bar \xi_1,\xi_2,\bar\xi_2)+i\omega_2 \xi_2\dd{}{\xi_2}P^{(2)}_l(\xi_1,\bar \xi_1,\xi_2,\bar\xi_2) +i \frac{2 \pi l}{T} P^{(2)}_l(\xi_1,\bar \xi_1,\xi_2,\bar \xi_2) \\
  &&= i \omega_1 \bar \xi_1\dd{}{\bar \xi_1}P^{(2)}_l(\xi_1,\bar \xi_1,\xi_2,\bar\xi_2)+i \omega_2 \bar \xi_2\dd{}{\bar \xi_2}P^{(2)}_l(\xi_1,\bar \xi_1,\xi_2,\bar\xi_2)+i\omega_2P^{(2)}_l(\xi_1,\bar \xi_1,\xi_2,\bar \xi_2)\\
 &&i\omega_1 \xi_1\dd{}{\xi_1}\bar P^{(2)}_l(\xi_1,\bar \xi_1,\xi_2,\bar\xi_2)+i\omega_2 \xi_2\dd{}{\xi_2}\bar P^{(2)}_l(\xi_1,\bar \xi_1,\xi_2,\bar\xi_2)+i\omega_1 \bar P^{(2)}_l(\xi_1,\bar \xi_1,\xi_2,\bar \xi_2) \\
 && = i \omega_1 \bar \xi_1\dd{}{\bar \xi_1}\bar P^{(2)}_l(\xi_1,\bar \xi_1,\xi_2,\bar\xi_2)+i \omega_2 \bar \xi_2\dd{}{\bar \xi_2}\bar P^{(2)}_l(\xi_1,\bar \xi_1,\xi_2,\bar\xi_2)-i \frac{2 \pi l}{T} \bar P^{(2)}_l(\xi_1,\bar \xi_1,\xi_2,\bar \xi_2).
\end{eqnarray*}
$p_l(\xi_1,\bar \xi_1,\xi_2,\bar \xi_2),P^{(1)}_l(\xi_1,\bar \xi_1,\xi_2,\bar \xi_2)$ and $P^{(2)}_l(\xi_1,\bar \xi_1,\xi_2,\bar \xi_2)$
are polynomials, and from the equations follows that $l$ has to be equal to zero, thus the polynomials are $\tau$-independent.
We obtain
\begin{eqnarray*}
 \xi_1\dd{}{\xi_1}p_0(\xi_1,\bar \xi_1,\xi_2,\bar\xi_2)& =&\bar \xi_1\dd{}{\bar \xi_1}p_0(\xi_1,\bar \xi_1,\xi_2,\bar\xi_2), \\
 \xi_2\dd{}{\xi_2}p_0(\xi_1,\bar \xi_1,\xi_2,\bar\xi_2)& =&\bar \xi_2\dd{}{\bar \xi_2}p_0(\xi_1,\bar \xi_1,\xi_2,\bar\xi_2), \\
 \xi_1\dd{}{\xi_1}P^{(1)}_0(\xi_1,\bar \xi_1,\xi_2,\bar\xi_2)& =& \bar\xi_1\dd{}{\bar\xi_1}P^{(1)}_0(\xi_1,\bar \xi_1,\xi_2,\bar\xi_2) + P^{(1)}_0(\xi_1,\bar \xi_1,\xi_2,\bar\xi_2), \\
 \xi_2\dd{}{\xi_2}P^{(1)}_0(\xi_1,\bar \xi_1,\xi_2,\bar\xi_2)& =& \bar\xi_2\dd{}{\bar\xi_2}P^{(1)}_0(\xi_1,\bar \xi_1,\xi_2,\bar\xi_2), \\
 \xi_1\dd{}{\xi_1}P^{(2)}_0(\xi_1,\bar \xi_1,\xi_2,\bar\xi_2)& =& \bar\xi_1\dd{}{\bar\xi_1}P^{(2)}_0(\xi_1,\bar \xi_1,\xi_2,\bar\xi_2), \\
 \xi_2\dd{}{\xi_2}P^{(2)}_0(\xi_1,\bar \xi_1,\xi_2,\bar\xi_2)& =& \bar\xi_2\dd{}{\bar\xi_2}P^{(2)}_0(\xi_1,\bar \xi_1,\xi_2,\bar\xi_2)+ P^{(2)}_0(\xi_1,\bar \xi_1,\xi_2,\bar\xi_2). 
\end{eqnarray*}
From the first two equations follows that 
\begin{eqnarray*}
p_0(\xi_1,\bar \xi_1,\xi_2,\bar\xi_2) = \psi_1(\left|\xi_1\right|^2)+\psi_2(\left|\xi_2\right|^2)+\psi_3(\left|\xi_1\right|^2)\psi_4(\left|\xi_2\right|^2).
\end{eqnarray*}
From the third and fourth equation, we obtain 
\begin{eqnarray*}
P^{(1)}_0(\xi_1,\bar \xi_1,\xi_2,\bar\xi_2) = \xi_1\phi_1(\left|\xi_1\right|^2)+\xi_1\phi_2(\left|\xi_2\right|^2)+\xi_1\phi_3(\left|\xi_1\right|^2)\phi_4(\left|\xi_2\right|^2),
\end{eqnarray*}
and analogously
\begin{eqnarray*}
P^{(2)}_0(\xi_1,\bar \xi_1,\xi_2,\bar\xi_2) = \xi_2\chi_1(\left|\xi_2\right|^2)+\xi_2\chi_2(\left|\xi_1\right|^2)+\xi_2\chi_3(\left|\xi_1\right|^2)\chi_4(\left|\xi_2\right|^2),
\end{eqnarray*}
where all functions are at least linear in their argument. 

Assembling all the information gives us the  Iooss normal form (\ref{eq:NF-NSNS}), i.e.
\begin{equation*}
\begin{cases} \ds\dd{\tau}{t}=1+\alpha_{1100}\left|\xi_1\right|^2+\alpha_{0011}\left|\xi_2\right|^2+\alpha_{2200}\left|\xi_1\right|^4+\alpha_{0022} \left|\xi_2\right|^4+\alpha_{1111} \left|\xi_1\right|^2\left|\xi_2\right|^2+\ldots, \smallskip \\
 \ds\dd{\xi_1}{\tau}= i\omega_1 \xi_1 + a_{2100} \xi_1\left|\xi_1\right|^2+a_{1011}\xi_1\left|\xi_2\right|^2+a_{3200}\xi_1\left|\xi_1\right|^4+a_{1022}\xi_1\left|\xi_2\right|^4+a_{2111}\xi_1\left|\xi_1\right|^2\left|\xi_2\right|^2+\ldots, \smallskip \\
 \ds\dd{\xi_2}{\tau}= i\omega_2 \xi_2 + b_{0021} \xi_2\left|\xi_2\right|^2+b_{1110}\xi_2\left|\xi_1\right|^2+b_{0032}\xi_2\left|\xi_2\right|^4+b_{2210}\xi_2\left|\xi_1\right|^4 +b_{1121}\xi_2\left|\xi_1\right|^2\left|\xi_2\right|^2+\ldots,
\end{cases}
\label{normalformnsns}
\end{equation*}
where the dots denote $O(\left|\xi\right|^6)$ terms. Note that the time evolution can be obtained by applying the chain 
rule to this system and is of the form
\begin{equation*}
\begin{cases} \dd{\tau}{t}=1+\alpha_{1100}\left|\xi_1\right|^2+\alpha_{0011}\left|\xi_2\right|^2+\alpha_{2200}\left|\xi_1\right|^4+\alpha_{0022} \left|\xi_2\right|^4+\alpha_{1111} \left|\xi_1\right|^2\left|\xi_2\right|^2+\ldots, \smallskip \\
 \dd{\xi_1}{t}=i\omega_1\xi_1 +a_{2100}'\xi_1\left|\xi_1\right|^2+a_{1011}'\xi_1\left|\xi_2\right|^2+a_{3200}'\xi_1\left|\xi_1\right|^4+a_{1022}'\xi_1\left|\xi_2\right|^4+a_{2111}'\xi_1\left|\xi_1\right|^2\left|\xi_2\right|^2+\ldots, \smallskip \\
 \dd{\xi_2}{t} =i\omega_2\xi_2 +b_{0021}'\xi_2\left|\xi_2\right|^2+b_{1110}'\xi_2\left|\xi_1\right|^2+b_{0032}'\xi_2\left|\xi_2\right|^4+b_{2210}'\xi_2\left|\xi_1\right|^4+b_{1121}'\xi_2\left|\xi_1\right|^2\left|\xi_2\right|^2+\ldots, \smallskip\\
\end{cases}
\end{equation*}
where the coefficients with primes are functions of the original coefficients.\\

\section{Bifurcations of the amplitude system for Hopf-Hopf bifurcation in the ``difficult" case}
\label{Appendix:6}

Here we derive quadratic approximations of the Hopf and heteroclinic bifurcation 
curves for the Hopf-Hopf truncated amplitude system (\ref{pdnsamplitude}), that can be written in the rescaled form
\begin{equation*}\label{dh:amp}
\left(\begin{array}{c} \dot{x}\\ \dot{y}\end{array}\right)=
\left(\begin{array}{c} x (\mu_{1}+x-\theta y+\Theta y^{2})\\ y (\mu_{2}+\delta x- y+\Delta x^{2})\end{array}\right) . \end{equation*}

The main results are
\begin{eqnarray*}
\mu_{2,Hopf}&=&-\frac{\delta-1}{\theta-1}\mu_{1}-\frac{(\delta-1)\Theta+(\theta-1)\Delta}{(\theta-1)^{3}}\mu_{1}^{2} 
+ O(\mu_1^3),\\
\mu_{2,Het}&=&-\frac{\delta-1}{\theta-1}\mu_{1}+\frac{\theta\Theta(\delta-1)^{3}+\delta\Delta(\theta-1)^{3}}{(\theta-1)^{3}(2\delta\theta-\delta-\theta)} \mu_{1}^{2} + O(\mu_1^3),\\
l_{1} &=& -C\left[\theta\left(\delta(\delta-1)\Theta+\theta(\theta-1)\Delta\right)\right],~~~ C>0.
\end{eqnarray*}

For the Hopf bifurcation curve we impose the conditions 
$\dot{x}=0,\dot{y}=0$ and $\frac{\partial \dot{x}}{\partial x} + \frac{\partial \dot{y}}{\partial y}=0$. 
Solving a series expansion yields the result for the curve. Next, the first Lyapunov coefficient $l_{1}$ 
is computed using the invariant formula (5.39) from \cite{Ku:2004} that leads to (\ref{l_1}).

For the heteroclinic curve we proceed as follows \cite{Chow:1994}. We assume $\delta,\theta<0$ and $\delta\theta-1>0$ and 
we transform variables to obtain a system that is a perturbation of a Hamiltonian system. This enables 
us to formulate a Melnikov function. Setting this function to zero yields an equation from which we extract 
the quadratic approximation to the heteroclinic curve. We introduce the transformation 
$(t,x,y,\mu_{1},\mu_{2})\rightarrow(\varepsilon x^{p-1}y^{q-1}t, \varepsilon x, \varepsilon y,-\varepsilon,c_{1} \varepsilon +c_{2}\varepsilon^{2})$ where 
$$
c_{1}=\frac{\delta-1}{\theta-1}, p=\frac{1-\delta}{\delta\theta-1}, q=\frac{1-\theta}{\delta\theta-1}.
$$
Then we obtain 
\begin{equation*}\label{dh:ham}
\left(\begin{array}{c} \dot{x}\\ \dot{y}\end{array}\right)=
x^{p-1}y^{q-1}\left(\begin{array}{c} x (-1+x-\theta y)\\ y (c_{1}+\delta x- y)\end{array}\right) +\varepsilon x^{p-1}y^{q-1} \left(\begin{array}{c} \Theta xy^{2}\\ \Delta c_{2}y+yx^{2}\end{array}\right), \end{equation*}
which for $\epsilon=0$ is a Hamilton system with Hamiltonian
$$
H(x,y)=\frac{1}{q}x^{p}y^{q}\left(-1+x+\frac{\theta-1}{\delta-1}y\right).
$$
Define $g_1=\Theta x^{p}y^{q+1}$ and $g_2=x^{p-1}y^{q}\left(c_{2}+\Delta x^{2}\right)$. 
The Melnikov function along $H(x,y)=h$ is given by the following integral
\begin{eqnarray*}
M(h) &=& \int_{H=h} g_{1}dy-g_{2}dx \\
&=& \int_{H=h} -\Theta x^{p}y^{q+1}dy - x^{p-1}y^{q}\left(c_{2}+\Delta x^2\right) dx\\
&=& \int_{H=h} -\left(x^{p-1}y^{q+2}\frac{p\Theta}{q+2} +x^{p-1}y^{q}\left(c_{2}+\Delta x^2\right)\right) dx,
\end{eqnarray*}
where we used Green's Theorem to convert the $dy$ term to $dx$. Now along the nontrivial critical 
curve $H(x,y)=0$ we have $y=\frac{\delta-1}{\theta-1}(1-x)$ so that 
\begin{eqnarray*}
M(0)&=&-\left(\frac{\delta-1}{\theta-1}\right)^{q}\int_{0}^{1} 
x^{p-1}(1-x)^{q} c_{2}+x^{p-1}y^{q+2}\left(\frac{\Theta p(\delta-1)^2}{(q+2)(\theta-1)^2} +\Delta\right)dx\\
&\sim &c_{2} I_{p-1,q}+\Theta I_{p-1,q+2}\left(\frac{\Theta p(\delta-1)^2}{(q+2)(\theta-1)^2} +\Delta\right),
\end{eqnarray*}
where we defined 
$$I_{a,b}=\int_{0}^{1}(1-x)^{a}x^{b}dx=\frac{\Gamma(1+a)\Gamma(1+b)}{\Gamma(2+a+b)}.$$
Solving $M(0)=0$ and substituting $p,q$ we obtain 
$$c_{2}=\frac{1}{(\theta-1)^{3}(2\delta\theta-\delta-\theta)}
	\left(\theta\Theta(\delta-1)^{3}+\delta\Delta(\theta-1)^{3}\right).$$

As a final check we consider the difference between the Hopf and heteroclinic curves
\begin{equation*}
\mu_{2,HET}-\mu_{2,HOPF}=\frac{(\delta\theta-1) l_{1}}{C\theta(\theta-1)^{3}(2\delta\theta-\delta-\theta)} \mu_1^2 + O(\mu_1^3).
\end{equation*}
We see that the quadratic approximations of the curves coincide when the Hopf bifurcation is degenerate.

\section{Higher order coefficients}
\label{Appendix:5}
In this appendix we list the third order normal form coefficients for the LPNS bifurcation 
and the fourth and fifth order coefficients for PDNS  and NSNS, which are necessary 
to determine the stability of the tori (if they exist). Remark that we have not listed the 
coefficients (or necessary coefficients of the expansion of the critical center manifold) 
which can be obtained by complex conjugacy or the 
similar expressions for $\omega_2$ instead of $\omega_1$ in the case of { NSNS}.

\subsection{Third order coefficients for LPNS} \label{Appendix:5_lpns}
The normal form coefficients in (\ref{eq:NF-LPNS}):
\begin{equation*} \label{eq:a300_LPNS}
 a_{300}=\frac{1}{6}\int_0^{T} \langle \varphi^*,C(\tau;v_1,v_1,v_1)+3 B(\tau;v_1,h_{200})+3\dot h_{200} -6a_{200} h_{200} -6\alpha_{200}A(\tau) v_1 \rangle\; d \tau+a_{200}
\end{equation*}
\begin{equation*} \label{eq:b210_LPNS}
 b_{210}=\frac{1}{2}\int_0^{T} \langle v_2^*,C(\tau;v_1,v_1,v_2) + B(\tau;v_2,h_{200}) + 2B(\tau;v_1,h_{110})+2\dot h_{110}-2\alpha_{200}A(\tau) v_2 \rangle\; d \tau+b_{110}
\end{equation*}
\begin{equation*} \label{eq:b021_LPNS}
 b_{021}=\frac{1}{2}\int_0^{T} \langle v_2^*,C(\tau;v_2,v_2,\bar v_2) +B(\tau;\bar v_2,h_{020}) + 2B(\tau;v_2,h_{011})  -2\alpha_{011} A(\tau)v_2 \rangle \; d\tau
\end{equation*}
\begin{eqnarray*} \label{eq:a111_LPNS}\\
 a_{111}&=&\int_0^{T} \langle \varphi^*,C(\tau;v_1,v_2,\bar v_2)+2\Re(B(\tau;v_2,h_{101}))+B(\tau;v_1,h_{011}) +\dot h_{011}-\alpha_{011}A(\tau) v_1\nonumber\\
 &&-2\Re(b_{110})h_{011}-a_{011}h_{200} \rangle \; d\tau +a_{011}\nonumber
\end{eqnarray*}

\subsection{Fourth and fifth order coefficients for PDNS}
The normal form coefficients in (\ref{eq:NF-PDNS}):

\begin{eqnarray*}
\alpha_{400}&=& \frac{1}{24} \int_0^{T} \langle \varphi^*,D(\tau;v_1,v_1,v_1,v_1)+6C(\tau;v_1,v_1,h_{200})+3B(\tau;h_{200},h_{200})\\
&&+4B(\tau;v_1,h_{300})-12 \alpha_{200} \dot h_{200}\rangle \; d\tau = 0
\end{eqnarray*}

\begin{eqnarray*}
\alpha_{211}&= &\frac{1}{2} \int_0^{T} \langle \varphi^*,D(\tau;v_1,v_1,v_2,\bar v_2)+C(\tau;v_1,v_1,h_{011})+C(\tau;v_2,\bar v_2,h_{200})+4\Re(C(\tau;v_1,v_2,h_{101}))\nonumber\\
 &&+2\Re(B(\tau;v_2,h_{201}))+B(\tau;h_{200},h_{011})+2B(\tau;h_{101},h_{110})+2B(\tau;v_1,h_{111})\nonumber\\
 &&- \alpha_{011} \dot h_{200}- 2\alpha_{200} \dot h_{011}\rangle \; d\tau\nonumber \label{eq:alpha111_PDNS}
\end{eqnarray*}

\begin{eqnarray*}
\alpha_{022}&= &\frac{1}{4} \int_0^{T} \langle \varphi^*,D(\tau;v_2,v_2,\bar v_2,\bar v_2)+4C(\tau;v_2,\bar v_2,h_{011})+2\Re(C(\tau;v_2,v_2,h_{002}))\nonumber\\
 &&+B(\tau;h_{020},h_{002})+2B(\tau;h_{011},h_{011})+4\Re(B(\tau;v_2,h_{012}))\nonumber\\
 &&-4\alpha_{011} \dot h_{011}\rangle \; d\tau\nonumber \\ \label{eq:alpha022_PDNS}
\end{eqnarray*}

Fourth order coefficients of the expansion of the critical center manifold can be computed by solving the
following BVPs on $[0,T]$:

\begin{equation*}  \label{eq:BVP_h400_PDNS}
\left\{\begin{array}{rcl}
 \dot h_{400}-A(\tau) h_{400}-D(\tau;v_1,v_1,v_1,v_1)-6C(\tau;v_1,v_1,h_{200})-3B(\tau;h_{200},h_{200})&&\\
 -4B(\tau;v_1,h_{300})+12 \alpha_{200} \dot h_{200}+24\alpha_{400}\dot u_0+24a_{300}h_{200} & = & 0\\
h_{400}(T)- h_{400}(0) & = & 0\\
\int_0^{T} \langle \varphi^*,h_{400}\rangle \; d\tau & = & 0
\end{array}
\right.
\end{equation*}

\begin{equation*}  \label{eq:BVP_h040_PDNS}
\left\{\begin{array}{rcl}
 \dot h_{040}-A(\tau) h_{040}+4i\omega h_{040}-D(\tau;v_2,v_2,v_2,v_2)-6 C(\tau;v_2,v_2,h_{020})&&\\
 -4B(\tau;v_2,h_{030}) -3B(\tau;h_{020},h_{020}) & = & 0\\
h_{040}(T)- h_{040}(0) & = & 0
\end{array}
\right.
\end{equation*}

\begin{equation*}  \label{eq:BVP_h310_PDNS}
\left\{\begin{array}{rcl}
 \dot h_{310}-A(\tau) h_{310}+i\omega h_{310}-D(\tau;v_1,v_1,v_1,v_2) -3C(\tau;v_1,v_1,h_{110})&&\\
 -3C(\tau;v_1,v_2,h_{200})-B(\tau;v_2,h_{300}) -3B(\tau;v_1,h_{210})-3B(\tau;h_{200},h_{110})&&\\
 +6 \alpha_{200} \dot h_{110}+6a_{300}h_{110}+6b_{210}h_{110}+6i\omega \alpha_{200}h_{110}& = & 0 \\
h_{310}(T)+ h_{310}(0) & = & 0
\end{array}
\right.
\end{equation*}

\begin{equation*}  \label{eq:BVP_h130_PDNS}
\left\{\begin{array}{rcl}
\dot h_{130}-A(\tau)h_{130}+3i\omega h_{130}-D(\tau;v_1,v_2,v_2,v_2)-3C(\tau;v_2,v_2,h_{110})&&\\
-3C(\tau;v_1,v_2,h_{020})-B(\tau;v_1,h_{030}) -3B(\tau;h_{020},h_{110})-3B(\tau;v_2,h_{120})& = & 0 \\
h_{130}(T)+ h_{130}(0) & = & 0
\end{array}
\right.
\end{equation*}

\begin{equation*}  \label{eq:BVP_h031_PDNS}
\left\{\begin{array}{rcl}
\dot h_{031}-A(\tau)h_{031}+2i\omega h_{031} -D(\tau;v_2,v_2,v_2,\bar v_2)-3C(\tau;v_2,v_2,h_{011})&&\\
-3C(\tau;v_2,\bar v_2,h_{020})-B(\tau;\bar v_2,h_{030})-3B(\tau;h_{020},h_{011})-3B(\tau;v_2,h_{021})&&\\
+3\alpha_{011}\dot h_{020}+6b_{021}h_{020}+6i\omega\alpha_{011}h_{020}& = & 0\\
h_{031}(T)- h_{031}(0) & = & 0
\end{array}
\right.
\end{equation*}

\begin{eqnarray*}  \label{eq:BVP_h211_PDNS}
\left\{\begin{array}{rcl}
 \dot h_{211}-A(\tau) h_{211}-D(\tau;v_1,v_1,v_2,\bar v_2)-C(\tau;v_1,v_1,h_{011})-C(\tau;v_2,\bar v_2,h_{200})&&\\ -4\Re(C(\tau;v_1,v_2,h_{101}))-2\Re(B(\tau;v_2,h_{201}))-B(\tau;h_{200},h_{011})&&\\
 -2B(\tau;h_{101},h_{110})-2B(\tau;v_1,h_{111})+ \alpha_{011} \dot h_{200}+ 2\alpha_{200} \dot h_{011}&&\\
 +2\alpha_{211}\dot u_0+2a_{111}h_{200}+4\Re(b_{210})h_{011} & = & 0\\
 h_{211}(T)- h_{211}(0) & = & 0\\
\int_0^{T} \langle \varphi^*,h_{211}\rangle \; d\tau & = & 0
\end{array}
\right.
\end{eqnarray*}

\begin{equation*}  \label{eq:BVP_h121_PDNS}
\left\{\begin{array}{rcl}
\dot h_{121}-A(\tau)h_{121}+i\omega h_{121} -D(\tau;v_1,v_2,v_2,\bar v_2)-C(\tau;v_1,\bar v_2,h_{020})&&\\
-2C(\tau;v_1,v_2,h_{011})-C(\tau;v_2,v_2,h_{101})-2C(\tau;v_2,\bar v_2,h_{110})-B(\tau;v_1,h_{021})\nonumber\\
 -B(\tau;h_{020},h_{101})-2B(\tau;h_{011},h_{110})-2B(\tau;v_2,h_{111})-B(\tau;\bar v_2,h_{120})&&\nonumber\\
 +2\alpha_{011}\dot h_{110}+2b_{021}h_{110}+2a_{111}h_{110}+2i\omega\alpha_{011}h_{110}& = & 0 \\
h_{121}(T)+ h_{121}(0) & = & 0
\end{array}
\right.
\end{equation*}

\begin{equation*}  \label{eq:BVP_h220_PDNS}
\left\{\begin{array}{rcl}
\dot h_{220}-A(\tau)h_{220}+2i\omega h_{220}-D(\tau;v_1,v_1,v_2,v_2)-C(\tau;v_2,v_2,h_{200})&&\\
-4C(\tau;v_1,v_2,h_{110})-C(\tau;v_1,v_1,h_{020})-B(\tau;h_{200},h_{020})-2B(\tau;v_2,h_{210})&&\\
-2B(\tau;h_{110},h_{110})-2B(\tau;v_1,h_{120})+2\alpha_{200}\dot h_{020}+4b_{210}h_{020}+4i\omega\alpha_{200}h_{020}& = & 0 \\
h_{220}(T)- h_{220}(0) & = & 0
\end{array}
\right.
\end{equation*}

\begin{equation*}  \label{eq:BVP_h022_PDNS}
\left\{\begin{array}{rcl}
 \dot h_{022}-A(\tau) h_{022}-D(\tau;v_2,v_2,\bar v_2,\bar v_2)-4C(\tau;v_2,\bar v_2,h_{011})-2\Re(C(\tau;v_2,v_2,h_{002}))&&\\
 -B(\tau;h_{020},h_{002})-2B(\tau;h_{011},h_{011})-4\Re(B(\tau;v_2,h_{012}))\nonumber&&\\
 +4\alpha_{011} \dot h_{011}+4\alpha_{022}\dot u_0+8\Re(b_{021})h_{011} & = & 0\\
 h_{022}(T) - h_{022}(0) &=& 0\\
\int_0^{T} \langle \varphi^*,h_{022}\rangle \; d\tau & = & 0
\end{array}
\right.
\end{equation*}

Fifth order normal form coefficients in (\ref{eq:NF-PDNS}):

\begin{eqnarray*} 
 a_{500}&=&\frac{1}{120}\int_0^{T} \langle v_1^*, E(\tau;v_1,v_1,v_1,v_1,v_1)+10D(\tau;v_1,v_1,v_1,h_{200})+10C(\tau;v_1,v_1,h_{300})\\
 &&+15C(\tau;v_1,h_{200},h_{200})+10B(\tau;h_{200},h_{300})+5B(\tau;v_1,h_{400})-20\alpha_{200}\dot h_{300}\\
 &&-120\alpha_{400}A(\tau) v_1\rangle\; d \tau-\alpha_{200}a_{300}\nonumber \\  \label{eq:a500_PDNS}
\end{eqnarray*}

\begin{eqnarray*}
 b_{410}&=&\frac{1}{24}\int_0^{T} \langle v_2^*, E(\tau;v_1,v_1,v_1,v_1,v_2)+6D(\tau;v_1,v_1,v_2,h_{200})+4D(\tau;v_1,v_1,v_1,h_{110})\\ &&+4C(\tau;v_1,v_2,h_{300})+6C(\tau;v_1,v_1,h_{210})+3C(\tau;v_2,h_{200},h_{200})+12C(\tau;v_1,h_{200},h_{110})\\
 &&+4B(\tau;v_1,h_{310})+4B(\tau;h_{110},h_{300})+6B(\tau;h_{200},h_{210})+B(\tau;v_2,h_{400})\\
 &&-24 \alpha_{400}A(\tau) v_2 -12\alpha_{200}\dot h_{210}\rangle\; d \tau-\alpha_{200}b_{210} \nonumber \\  \label{eq:b410_PDNS}
\end{eqnarray*}

\begin{eqnarray*} 
 a_{311}&=&\frac{1}{6}\int_0^{T} \langle v_1^*,E(\tau;v_1,v_1,v_1,v_2,\bar v_2)+D(\tau;v_1,v_1,v_1,h_{011})+6\Re(D(\tau;v_1,v_1,v_2,h_{101}))\\
 &&+3D(\tau;v_1,v_2,\bar v_2,h_{200})+3C(\tau;v_1,h_{200},h_{011})+6\Re(C(\tau;v_2,h_{200},h_{101}))\nonumber\\
 &&+6\Re(C(\tau;v_1,v_2,h_{201}))+C(\tau;v_2,\bar v_2,h_{300})\nonumber\\
 &&+3C(\tau;v_1,v_1,h_{111})+6C(\tau;v_1,h_{101},h_{110})+3B(\tau;h_{200},h_{111})\nonumber\\
 &&+6\Re(B(\tau;h_{201},h_{110}))+2\Re(B(\tau;v_2,h_{301}))+B(\tau;h_{011},h_{300})\nonumber\\
 &&+3B(\tau;h_{211},v_1)-6 \alpha_{211}A(\tau) v_1\nonumber\\
  &&-\alpha_{011}\dot h_{300} -6\alpha_{200}\dot h_{111}\rangle\; d \tau-\alpha_{200}a_{111}-\alpha_{011}a_{300} \nonumber \\ \label{eq:a311_PDNS}
\end{eqnarray*}

\begin{eqnarray*} 
 b_{221}&=&\frac{1}{4}\int_0^{T} \langle v_2^*,E(\tau;v_1,v_1,v_2,v_2,\bar v_2)+D(\tau;v_2,v_2,\bar v_2,h_{200})\\
 &&+2D(\tau;v_1,v_2,v_2,h_{101})+2D(\tau;v_1,v_1, v_2,h_{011})\nonumber\\
 &&+D(\tau;v_1,v_1,\bar v_2,h_{020})+4D(\tau;v_1,v_2,\bar v_2,h_{110})+2C(\tau;\bar v_2,h_{110},h_{110})+C(\tau;v_1,v_1,h_{021})\nonumber\\
 &&+C(\tau;v_2,v_2,h_{201})+C(\tau;\bar v_2,h_{200},h_{020})+2C(\tau;v_2,\bar v_2,h_{210})\nonumber\\
 &&+2C(\tau;v_1,\bar v_2,h_{120})+2C(\tau;v_1,h_{020},h_{101})+4C(\tau;v_1,v_2,h_{111})+4C(\tau;v_2,h_{101},h_{110})\nonumber\\
 &&+2C(\tau;v_2,h_{200},h_{011})+4C(\tau;v_1,h_{110},h_{011})+B(\tau;\bar v_2,h_{220})\nonumber\\
 &&+2B(\tau;v_1,h_{121})+2B(\tau;h_{120},h_{101})+4B(\tau;h_{110},h_{111})+2B(\tau;h_{210},h_{011})\nonumber\\
 &&+2B(\tau;v_2,h_{211})+B(\tau;h_{200},h_{021})+B(\tau;h_{201},h_{020})\nonumber\\
 &&-2\alpha_{011}\dot h_{210}-4\alpha_{211}A(\tau) v_2-2\alpha_{200}\dot h_{021}\rangle\; d \tau-\alpha_{200}b_{021}-\alpha_{011}b_{210} \nonumber \\ \label{eq:b221_PDNS}
\end{eqnarray*}

\begin{eqnarray*} 
 a_{122}&=&\frac{1}{4}\int_0^{T} \langle v_1^*,E(\tau;v_1,v_ 2,v_2,\bar v_2,\bar v_2)+4\Re(D(\tau;v_2,v_2,\bar v_2,h_{101}))+2\Re(D(\tau;v_1,v_2, v_2,h_{002}))\nonumber\\
 &&+4D(\tau;v_1,v_2,\bar v_2,h_{011})+2C(\tau;v_1,h_{011},h_{011})\nonumber\\
 &&+8\Re(C(\tau;v_2,h_{011},h_{101}))+4\Re(C(\tau;v_2,h_{002},h_{110}))\nonumber\\
 &&+C(\tau;v_1,h_{020},h_{002})+4\Re(C(\tau;v_1,v_2,h_{012}))+4C(\tau;v_2,\bar v_2,h_{111})\nonumber\\
 &&+2\Re(C(\tau;v_2,v_2,h_{102}))+4B(\tau;h_{011},h_{111})\nonumber\\
 &&+2\Re(B(\tau;h_{020},h_{102}))+4\Re(B(\tau;v_2,h_{112}))+B(\tau;v_1,h_{022})\nonumber\\
 &&+4\Re(B(\tau;h_{110},h_{012}))-4\alpha_{011}\dot h_{111}\nonumber\\
 &&-4\alpha_{022}A(\tau) v_1\rangle\; d \tau-\alpha_{011}a_{111}\nonumber \\  \label{eq:a122_PDNS}
\end{eqnarray*}

\begin{eqnarray*} 
 b_{032}&=&\frac{1}{12}\int_0^{T} \langle v_2^*,E(\tau;v_2,v_2,v_2,\bar v_2,\bar v_2)+D(\tau;v_2,v_2, v_2,h_{002})+3D(\tau; v_2,\bar v_2,\bar v_2,h_{020})\\
 &&+6D(\tau;v_2,v_2,\bar v_2,h_{011})+6C(\tau;\bar v_2,h_{020},h_{011})+6C(\tau;v_2,\bar v_2,h_{021})+C(\tau;\bar v_2,\bar v_2,h_{030})\nonumber\\
 &&+3C(\tau;v_2,h_{002},h_{020})+3C(\tau;v_2,v_2,h_{012})+6C(\tau;v_2,h_{011},h_{011})\nonumber\\
 &&+3B(\tau;h_{020},h_{021})+6B(\tau;h_{021},h_{011})+3B(\tau;v_2,h_{022})+2B(\tau;\bar v_2,h_{031})\nonumber\\
 &&+B(\tau;h_{002},h_{030})-6\alpha_{011}\dot h_{021}-12\alpha_{022}A(\tau) v_2\rangle\; d \tau -\alpha_{011}b_{021}\nonumber \\\label{eq:b032_PDNS}
\end{eqnarray*}

\subsection{Fourth and fifth order coefficients for NSNS}

Fourth order normal form coefficients for (\ref{eq:NF-NSNS}):

\begin{eqnarray*} \label{eq:alpha2200_NSNS}
 \alpha_{2200}&=&\frac{1}{4} \int_0^{T} \langle \varphi^*,D(\tau;v_1,v_1,\bar v_1,\bar v_1)+2\Re(C(\tau;v_1,v_1,h_{0200}))\nonumber\\
 &&+4C(\tau;v_1,\bar v_1,h_{1100})+B(\tau;h_{2000},h_{0200})+2B(\tau;h_{1100},h_{1100})+4\Re(B(\tau;v_1,h_{1200}))\\
 &&-4\alpha_{1100}\dot h_{1100}\rangle\; d \tau\nonumber\\
\end{eqnarray*}

\begin{eqnarray*} \label{eq:alpha1111_NSNS}
 \alpha_{1111}&=&\int_0^{T} \langle \varphi^*,D(\tau;v_1,\bar v_1,v_2,\bar v_2)+C(\tau;v_1,\bar v_1,h_{0011})+2\Re(C(\tau;v_1,v_2,h_{0101}))\nonumber\\
 &&+2\Re(C(\tau;v_1,\bar v_2,h_{0110}))\nonumber \\
 &&+C(\tau;v_2,\bar v_2,h_{1100}) +2\Re(B(\tau;v_1,h_{0111}))+B(\tau;h_{0110},h_{1001})+B(\tau;h_{0101},h_{1010})\nonumber \\
 &&+B(\tau;h_{0011},h_{1100})+2\Re(B(\tau;v_2,h_{1101}))\nonumber\\
 &&-\alpha_{0011}\dot h_{1100}-\alpha_{1100}\dot h_{0011}\rangle\; d \tau\nonumber\\
\end{eqnarray*}

Fourth order coefficients of the expansion of the critical center manifold can be computed by solving the
following BVPs on $[0,T]$:

\begin{equation*}  \label{eq:BVP_h4000_NSNS}
\left\{\begin{array}{rcl}
\dot h_{4000}-A(\tau) h_{4000}+4i\omega_1 h_{4000}-D(\tau;v_1,v_1,v_1,v_1)-6C(\tau;v_1,v_1,h_{2000})&&\nonumber\\
 -3B(\tau;h_{2000},h_{2000})-4B(\tau;v_1,h_{3000})& = & 0 \\
h_{4000}(T)- h_{4000}(0) & = & 0
\end{array}
\right.
\end{equation*}

\begin{equation*}  \label{eq:BVP_h3100_NSNS}
\left\{\begin{array}{rcl}
\dot h_{3100}-A(\tau) h_{3100}+2i\omega_1 h_{3100}-D(\tau;v_1,v_1,v_1,\bar v_1)-3C(\tau;v_1,v_1,h_{1100})&&\\
-3C(\tau;v_1,\bar v_1,h_{2000})-B(\tau;\bar v_1,h_{3000})-3B(\tau;v_1,h_{2100})-3B(\tau;h_{2000},h_{1100})&&\\
+3\alpha_{1100}\dot h_{2000}+6a_{2100}h_{2000}+i\omega_1\alpha_{1100}h_{2000}& = & 0 \\
h_{3100}(T)- h_{3100}(0) & = & 0
\end{array}
\right.
\end{equation*}

\begin{equation*}  \label{eq:BVP_h3010_NSNS}
\left\{\begin{array}{rcl}
\dot h_{3010}-A(\tau) h_{3010}+3i\omega_1 h_{3010}+i\omega_2 h_{3010}-D(\tau;v_1,v_1,v_1,v_2)-3C(\tau;v_1,v_1,h_{1010})&&\nonumber\\
-3C(\tau;v_1,v_2,h_{2000}) -B(\tau;v_2,h_{3000})-3B(\tau;v_1,h_{2010})-3B(\tau;h_{2000},h_{1010})& = & 0 \\
h_{3010}(T)- h_{3010}(0) & = & 0
\end{array}
\right.
\end{equation*}

\begin{equation*}  \label{eq:BVP_h3001_NSNS}
\left\{\begin{array}{rcl}
\dot h_{3001}-A(\tau) h_{3001}+3i\omega_1 h_{3001}-i\omega_2 h_{3001}-D(\tau;v_1,v_1,v_1,\bar v_2)-3C(\tau;v_1,v_1,h_{1001})&&\nonumber\\
-3C(\tau;v_1,\bar v_2,h_{2000}) -B(\tau;\bar v_2,h_{3000})-3B(\tau;v_1,h_{2001})-3B(\tau;h_{2000},h_{1001})& = & 0 \\
h_{3001}(T)- h_{3001}(0) & = & 0
\end{array}
\right.
\end{equation*}

\begin{equation*}  \label{eq:BVP_h2200_NSNS}
\left\{\begin{array}{rcl}
\dot h_{2200}-A(\tau) h_{2200}-D(\tau;v_1,v_1,\bar v_1,\bar v_1)-2\Re(C(\tau;v_1,v_1,h_{0200}))&&\nonumber\\
 -4C(\tau;v_1,\bar v_1,h_{1100})-B(\tau;h_{2000},h_{0200})-2B(\tau;h_{1100},h_{1100})-4\Re(B(\tau;v_1,h_{1200}))\nonumber\\
 +8\Re(a_{2100})h_{1100}+4\alpha_{2200}\dot u_0+4\alpha_{1100}\dot h_{1100}& = & 0 \\
h_{2200}(T)- h_{2200}(0) & = & 0\\
\int_0^{T} \langle \varphi^*,h_{2200}\rangle \; d\tau & = & 0
\end{array}
\right.
\end{equation*}

\begin{equation*}  \label{eq:BVP_h2020_NSNS}
\left\{\begin{array}{rcl}
 \dot h_{2020}-A(\tau) h_{2020}+2i\omega_1h_{2020}+2i\omega_2h_{2020}-D(\tau;v_1,v_1,v_2,v_2)-C(\tau;v_1,v_1,h_{0020})&&\nonumber\\
 -C(\tau;v_2,v_2,h_{2000}) -4C(\tau;v_1,v_2,h_{1010})-B(\tau;h_{2000},h_{0020})-2B(\tau; v_2,h_{2010})&&\nonumber\\
 -2B(\tau;h_{1010},h_{1010})-2B(\tau;v_1,h_{1020})& = & 0\\
h_{2020}(T)- h_{2020}(0) & = & 0
\end{array}
\right.
\end{equation*}

\begin{equation*}  \label{eq:BVP_h2002_NSNS}
\left\{\begin{array}{rcl}
 \dot h_{2002}-A(\tau) h_{2002}+2i\omega_1h_{2002}-2i\omega_2h_{2002}-D(\tau;v_1,v_1,\bar v_2,\bar v_2)-C(\tau;\bar v_2,\bar v_2,h_{2000})&&\nonumber\\
 -4C(\tau;v_1,\bar v_2,h_{1001}) -C(\tau;v_1,v_1,h_{0002})-2B(\tau;\bar v_2,h_{2001})-B(\tau; h_{2000},h_{0002})&&\nonumber\\
 -2B(\tau;h_{1001},h_{1001})-2B(\tau;v_1,h_{1002})& = & 0 \\
h_{2002}(T)- h_{2002}(0) & = & 0\\
\end{array}
\right.
\end{equation*}

\begin{equation*}  \label{eq:BVP_h2110_NSNS}
\left\{\begin{array}{rcl}
 \dot h_{2110}-A(\tau) h_{2110}+i\omega_1h_{2110}+i\omega_2h_{2110}-D(\tau;v_1,v_1,\bar v_1,v_2)&&\\
 -C(\tau;v_1, v_1,h_{0110})-2C(\tau;v_1,\bar v_1,h_{1010})-C(\tau;\bar v_1,v_2,h_{2000})-2C(\tau;v_1,v_2,h_{1100})&&\\
 -B(\tau;\bar v_1,h_{2010})-2B(\tau;h_{1010},h_{1100})-B(\tau;v_2,h_{2100})-B(\tau;h_{2000},h_{0110})&&\\
 -2B(\tau;v_1,h_{1110}) +2a_{2100}h_{1010}+2b_{1110}h_{1010}+2\alpha_{1100}\dot h_{1010}+2i(\omega_1+\omega_2)\alpha_{1100}h_{1010}& = & 0 \\
h_{2110}(T)- h_{2110}(0) & = & 0\\
\end{array}
\right.
\end{equation*}

\begin{equation*}  \label{eq:BVP_h2101_NSNS}
\left\{\begin{array}{rcl}
 \dot h_{2101}-A(\tau) h_{2101}+i\omega_1h_{2101}-i\omega_2h_{2101}-D(\tau;v_1,v_1,\bar v_1,\bar v_2)-C(\tau;v_1, v_1,h_{0101})&&\nonumber\\
 -2C(\tau;v_1,\bar v_1,h_{1001})-C(\tau;\bar v_1,\bar v_2,h_{2000})-2C(\tau;v_1,\bar v_2,h_{1100})-2B(\tau;h_{1001},h_{1100})&&\nonumber\\
 -2B(\tau;v_1,h_{1101})-B(\tau;\bar v_2,h_{2100})-B(\tau;\bar v_1,h_{2001})-B(\tau;h_{2000},h_{0101})&&\\
 +2a_{2100}h_{1001}+2b_{1101}h_{1001}+2\alpha_{1100}\dot h_{1001}+2i(\omega_1-\omega_2)\alpha_{1100}h_{1001}& = & 0 \\
h_{2101}(T)- h_{2101}(0) & = & 0\\
\end{array}
\right.
\end{equation*}

\begin{equation*}  \label{eq:BVP_h2011_NSNS}
\left\{\begin{array}{rcl}
 \dot h_{2011}-A(\tau) h_{2011}+2i\omega_1h_{2011}-D(\tau;v_1,v_1,v_2,\bar v_2)-C(\tau;v_1, v_1,h_{0011})&&\\
 -2C(\tau;v_1,v_2,h_{1001})-C(\tau;v_2,\bar v_2,h_{2000})-2C(\tau;v_1,\bar v_2,h_{1010})-B(\tau;\bar v_2,h_{2010})&&\\
 -B(\tau;v_2,h_{2001})-B(\tau;h_{2000} ,h_{0011})-2B(\tau;h_{1001},h_{1010})-2B(\tau;v_1,h_{1011})&&\\
 +2a_{1011}h_{2000}+\alpha_{0011}\dot h_{2000}+2i\omega_1\alpha_{0011}h_{2000}& = & 0 \\
h_{2011}(T)- h_{2011}(0) & = & 0\\
\end{array}
\right.
\end{equation*}

\begin{equation*}  \label{eq:BVP_h1111_NSNS}
\left\{\begin{array}{rcl}
\dot h_{1111}-A(\tau) h_{1111}-D(\tau;v_1,\bar v_1,v_2,\bar v_2)-C(\tau;v_1,\bar v_1,h_{0011})&&\\
-2\Re(C(\tau;v_1,v_2,h_{0101}))-2\Re(C(\tau;v_1,\bar v_2,h_{0110}))\nonumber &&\\
 -C(\tau;v_2,\bar v_2,h_{1100}) -2\Re(B(\tau;v_1,h_{0111}))-B(\tau;h_{0110},h_{1001})&&\\
 -B(\tau;h_{0101},h_{1010})-B(\tau;h_{0011},h_{1100})-2\Re(B(\tau;v_2,h_{1101}))\nonumber&&\\
 +2\Re(a_{0111})h_{1100}+\alpha_{0011}\dot h_{1100}+2\Re(b_{1101})h_{0011}+\alpha_{1111}\dot u_0+\alpha_{1100}\dot h_{0011} & = & 0\\
 h_{1111}(T)- h_{1111}(0) & = & 0\\
\int_0^{T} \langle \varphi^*,h_{1111}\rangle \; d\tau & = & 0
\end{array}
\right.
\end{equation*}

Fifth order normal form coefficients for (\ref{eq:NF-NSNS}):

\begin{eqnarray*} \label{eq:a3200_NSNS}
 a_{3200}&=&\frac{1}{12} \int_0^{T} \langle v_1^*,E(\tau;v_1,v_1,v_1,\bar v_1,\bar v_1)+D(\tau;v_1,v_1,v_1,h_{0200})+3D(\tau;v_1,\bar v_1,\bar v_1,h_{2000})\nonumber\\
 &&+6D(\tau;v_1,v_1,\bar v_1,h_{1100})+6C(\tau;v_1,h_{1100},h_{1100})+3C(\tau;v_1,v_1,h_{1200})\nonumber\\
 &&+C(\tau;\bar v_1,\bar v_1,h_{3000})+6C(\tau;v_1,\bar v_1,h_{2100})+6C(\tau;\bar v_1,h_{2000},h_{1100})+3C(\tau;v_1,h_{0200},h_{2000})\nonumber\\
 &&+B(\tau;h_{0200},h_{3000})+2B(\tau;\bar v_1,h_{3100})+3B(\tau;v_1,h_{2200})+6B(\tau;h_{2100},h_{1100})\nonumber\\
 &&+3B(\tau;h_{2000},h_{1200})-6\alpha_{1100}\dot h_{2100}\nonumber\\
 &&-12\alpha_{2200}A(\tau) v_1 \rangle\; d \tau-\alpha_{1100}a_{2100}\nonumber\\
\end{eqnarray*}

\begin{eqnarray*} \label{eq:b0032_NSNS}
 b_{0032}&=&\frac{1}{12} \int_0^{T} \langle v_2^*,E(\tau;v_2,v_2,v_2,\bar v_2,\bar v_2)+D(\tau;v_2,v_2,v_2,h_{0002})+3D(\tau;v_2,\bar v_2,\bar v_2,h_{0020})\nonumber\\
 &&+6D(\tau;v_2,v_2,\bar v_2,h_{0011})+6C(\tau;v_2,h_{0011},h_{0011})+3C(\tau;v_2,v_2,h_{0012})\nonumber\\
 &&+C(\tau;\bar v_2,\bar v_2,h_{0030})+6C(\tau;v_2,\bar v_2,h_{0021})+6C(\tau;\bar v_2,h_{0020},h_{0011})+3C(\tau;v_2,h_{0002},h_{0020})\nonumber\\
 &&+B(\tau;h_{0002},h_{0030})+2B(\tau;\bar v_2,h_{0031})+3B(\tau;v_2,h_{0022})+6B(\tau;h_{0021},h_{0011})\nonumber\\
 &&+3B(\tau;h_{0020},h_{0012})-6\alpha_{0011}\dot h_{0021}\nonumber\\
 &&-12\alpha_{0022}A(\tau) v_2 \rangle\; d \tau-\alpha_{0011}b_{0021}\nonumber\\
\end{eqnarray*}

\begin{eqnarray*} \label{eq:a1022_NSNS}
 a_{1022}&=&\frac{1}{4} \int_0^{T} \langle v_1^*,E(\tau;v_1,v_2,v_2,\bar v_2,\bar v_2)+D(\tau;v_1,\bar v_2,\bar v_2,h_{0020})+D(\tau;v_1,v_2,v_2,h_{0002})\nonumber \\
 &&+2D(\tau;v_2,v_2,\bar v_2,h_{1001})+2D(\tau;v_2,\bar v_2,\bar v_2,h_{1010})\nonumber\\
 &&+4D(\tau;v_1,v_2,\bar v_2,h_{0011})+2C(\tau;v_1,\bar v_2,h_{0021})+C(\tau;v_1,h_{0020},h_{0002})\nonumber\\
 &&+2C(\tau;v_1,v_2,h_{0012})+2C(\tau;\bar v_2,h_{1001},h_{0020})+4C(\tau;v_2,h_{1001},h_{0011})+2C(\tau;v_2,h_{1010},h_{0002})\nonumber\\
 &&+C(\tau;v_2,v_2,h_{1002})+4C(\tau;\bar v_2,h_{1010},h_{0011})+4C(\tau;v_2,\bar v_2,h_{1011})+C(\tau;\bar v_2,\bar v_2,h_{1020})\nonumber\\
 &&+2C(\tau;v_1,h_{0011},h_{0011})\nonumber\\
 &&+B(\tau;v_1,h_{0022})+2B(\tau;h_{0021},h_{1001})+B(\tau;h_{0020},h_{1002})+2B(\tau;h_{0012},h_{1010})\nonumber\\
 &&+4B(\tau;h_{0011},h_{1011})+2B(\tau;v_2,h_{1012})+B(\tau;h_{0002},h_{1020})+2B(\tau;\bar v_2,h_{1021})\nonumber\\
 &&-4\alpha_{0011}\dot h_{1011}-4\alpha_{0022}A(\tau)v_1\rangle\; d \tau -\alpha_{0011}a_{1011}\nonumber\\
\end{eqnarray*}

\begin{eqnarray*} \label{eq:b2210_NSNS}
 b_{2210}&=&\frac{1}{4} \int_0^{T} \langle v_2^*,E(\tau;v_1,v_1,\bar v_1,\bar v_1,v_2)+D(\tau;\bar v_1,\bar v_1, v_2,h_{2000})+D(\tau;v_1,v_1,v_2,h_{0200})\nonumber \\
 &&+2D(\tau;v_1,v_1,\bar v_1,h_{0110})+2D(\tau;v_1,\bar v_1,\bar v_1,h_{1010})\nonumber\\
 &&+4D(\tau;v_1,\bar v_1,v_2,h_{1100})+2C(\tau;\bar v_1,v_2,h_{2100})+C(\tau;v_2,h_{2000},h_{0200})\nonumber\\
 &&+2C(\tau;v_1,v_2,h_{1200})+2C(\tau;\bar v_1,h_{0110},h_{2000})+4C(\tau;v_1,h_{0110},h_{1100})+2C(\tau;v_1,h_{1010},h_{0200})\nonumber\\
 &&+C(\tau;v_1,v_1,h_{0210})+4C(\tau;\bar v_1,h_{1010},h_{1100})+4C(\tau;v_1,\bar v_1,h_{1110})+C(\tau;\bar v_1,\bar v_1,h_{2010})\nonumber\\
 &&+2C(\tau;v_2,h_{1100},h_{1100})\nonumber\\
 &&+B(\tau;v_2,h_{2200})+2B(\tau;h_{2100},h_{0110})+B(\tau;h_{2000},h_{0210})+2B(\tau;h_{1200},h_{1010})\nonumber\\
 &&+4B(\tau;h_{1100},h_{1110})+2B(\tau;v_1,h_{1210})+B(\tau;h_{0200},h_{2010})+2B(\tau;\bar v_1,h_{2110})\nonumber\\
 &&-4\alpha_{1100}\dot h_{1110}-4\alpha_{2200}A(\tau)v_2\rangle\; d \tau -\alpha_{1100}b_{1110}\nonumber\\
\end{eqnarray*}

\begin{eqnarray*} \label{eq:a2111_NSNS}
 a_{2111}&=&\frac{1}{2} \int_0^{T} \langle v_1^*,E(\tau;v_1,v_1,\bar v_1,v_2,\bar v_2)+D(\tau;v_1,v_1, v_2,h_{0101})+D(\tau;v_1,v_1,\bar v_1,h_{0011})\nonumber\\
 &&+D(\tau;v_1,v_1,\bar v_2,h_{0110})+2D(\tau;v_1,\bar v_1,\bar v_2,h_{1010})+2D(\tau;v_1,v_2,\bar v_2,h_{1100})\nonumber\\
 &&+D(\tau;\bar v_1,v_2,\bar v_2,h_{2000})+2D(\tau;v_1,\bar v_1,v_2,h_{1001})+2C(\tau;v_1,h_{1001} ,h_{0110})\nonumber\\ 
 &&+C(\tau;v_2,\bar v_2,h_{2100})+2C(\tau;v_2,h_{1001},h_{1100})+2C(\tau;v_1,h_{1100},h_{0011})\nonumber\\ 
 &&+2C(\tau;v_1,h_{1010},h_{0101})+C(\tau;\bar v_1,\bar v_2,h_{2010})+C(\tau;\bar v_1,v_2,h_{2001})\nonumber\\ 
 &&+2C(\tau;v_1,\bar v_2,h_{1110})+2C(\tau;\bar v_2,h_{1010},h_{1100})+2C(\tau;v_1,v_2,h_{1101})\nonumber\\
 &&+C(\tau;v_2,h_{2000},h_{0101})+2C(\tau;\bar v_1,h_{1001},h_{1010})+C(\tau;\bar v_1,h_{2000},h_{0011})\nonumber\\ 
 &&+C(\tau;\bar v_2,h_{2000},h_{0110})+2C(\tau;v_1,\bar v_1,h_{1011})+C(\tau;v_1,v_1,h_{0111})\nonumber\\
 &&+B(\tau;v_2,h_{2101})+B(\tau;h_{2100},h_{0011})+2B(\tau;v_1,h_{1111})+B(\tau;\bar v_1,h_{2011})\nonumber\\
 &&+B(\tau;h_{0101},h_{2010})+B(\tau;h_{2001},h_{0110})+B(\tau;h_{2000},h_{0111})+B(\tau;\bar v_2,h_{2110})\nonumber\\
 &&+2B(\tau;h_{1011},h_{1100})+2B(\tau;h_{1010},h_{1101})+2B(\tau;h_{1001},h_{1110})\nonumber\\
 &&-2\alpha_{1111}A(\tau) v_1-2\alpha_{1100}\dot h_{1011}-\alpha_{0011}\dot h_{2100}\rangle\; d \tau-\alpha_{0011}a_{2100}-\alpha_{1100}a_{1011}\nonumber\\
\end{eqnarray*}

\begin{eqnarray*} \label{eq:b1121_NSNS}
 b_{1121}&=&\frac{1}{2} \int_0^{T} \langle v_2^*,E(\tau;v_1,\bar v_1,v_2,v_2,\bar v_2)+D(\tau;v_1,v_2, v_2,h_{0101})+D(\tau;v_2,v_2,\bar v_2,h_{1100})\nonumber\\
 &&+D(\tau;\bar v_1,v_2,v_2,h_{1001})+2D(\tau;\bar v_1,v_2,\bar v_2,h_{1010})+2D(\tau;v_1,\bar v_1,v_2,h_{0011})\nonumber\\
 &&+D(\tau;v_1,\bar v_1,\bar v_2,h_{0020})+2D(\tau;v_1,v_2,\bar v_2,h_{0110})+2C(\tau;v_2,h_{0110} ,h_{1001})\nonumber\\ 
 &&+C(\tau;v_1,\bar v_1,h_{0021})+2C(\tau;v_1,h_{0110},h_{0011})+2C(\tau;v_2,h_{0011},h_{1100})\nonumber\\ 
 &&+2C(\tau;v_2,h_{1010},h_{0101})+C(\tau;\bar v_1,\bar v_2,h_{1020})+C(\tau;v_1,\bar v_2,h_{0120})\nonumber\\ 
 &&+2C(\tau;\bar v_1,v_2,h_{1011})+2C(\tau;\bar v_1,h_{1010},h_{0011})+2C(\tau;v_1,v_2,h_{0111})\nonumber\\
 &&+C(\tau;v_1,h_{0020},h_{0101})+2C(\tau;\bar v_2,h_{0110},h_{1010})+C(\tau;\bar v_2,h_{0020},h_{1100})\nonumber\\ 
 &&+C(\tau;\bar v_1,h_{0020},h_{1001})+2C(\tau;v_2,\bar v_2,h_{1110})+C(\tau;v_2,v_2,h_{1101})\nonumber\\
 &&+B(\tau;v_1,h_{0121})+B(\tau;h_{0021},h_{1100})+2B(\tau;v_2,h_{1111})+B(\tau;\bar v_2,h_{1120})\nonumber\\
 &&+B(\tau;h_{0101},h_{1020})+B(\tau;h_{0120},h_{1001})+B(\tau;h_{0020},h_{1101})+B(\tau;\bar v_1,h_{1021})\nonumber\\
 &&+2B(\tau;h_{1110},h_{0011})+2B(\tau;h_{1010},h_{0111})+2B(\tau;h_{0110},h_{1011})\nonumber\\
 &&-2\alpha_{1111}A(\tau) v_2-2\alpha_{0011}\dot h_{1110}-\alpha_{1100}\dot h_{0021}\rangle\; d \tau - \alpha_{1100}b_{0021}-\alpha_{0011}b_{1110}\nonumber\\
\end{eqnarray*}

\bibliographystyle{amsplain}
\bibliography{LC2NF}

\providecommand{\bysame}{\leavevmode\hbox to3em{\hrulefill}\thinspace}
\providecommand{\MR}{\relax\ifhmode\unskip\space\fi MR }
% \MRhref is called by the amsart/book/proc definition of \MR.
\providecommand{\MRhref}[2]{%
  \href{http://www.ams.org/mathscinet-getitem?mr=#1}{#2}
}
\providecommand{\href}[2]{#2}
\begin{thebibliography}{10}

\bibitem{capd}
\emph{{CAPD}: {C}omputer {A}ssisted {P}roofs in {D}ynamics}, {\tt
  http://capd.ii.uj.edu.pl}.

\bibitem{tides}
\emph{{TIDES}: A {T}aylor {I}ntegrator for {D}ifferential {E}quation{S}}, {\tt
  http://gme.unizar.es/software/tides}.

\bibitem{Ar:83}
{\mbox{V}}.~I. Arnol'd, \emph{Geometrical {M}ethods in the {T}heory of
  {O}rdinary {D}ifferential {E}quations}, Springer-Verlag, New York, 1983.

\bibitem{BVPbook:95}
Uri~M. Ascher, Robert M.~M. Mattheij, and Robert~D. Russell, \emph{Numerical
  solution of boundary value problems for ordinary differential equations},
  Classics in Applied Mathematics, vol.~13, Society for Industrial and Applied
  Mathematics (SIAM), Philadelphia, PA, 1995, Corrected reprint of the 1988
  original.

\bibitem{Chow:1994}
S.-N. Chow, C.~Li, and D.~Wang, \emph{Normal forms and bifurcations of planar
  vector fields}, Cambridge University Press, Cambridge, 1994.

\bibitem{BoSw:73}
C.~{De Boor} and B.~Swartz, \emph{Collocation at {G}aussian points}, SIAM J.
  Numer. Anal. \textbf{10} (1973), no.~4, 582--606.

\bibitem{drdwgk}
V.~{De Witte}, F.~{Della Rossa}, W.~Govaerts, and {\mbox{Yu}}.A. Kuznetsov,
  \emph{Numerical periodic normalization for codim $2$ bifurcations of limit
  cycles - computational formulas, numerical implementation, and examples},
  Submitted to SIADS.

\bibitem{MATCONT}
A.~Dhooge, W.~Govaerts, and {\mbox{Yu}}.~A. Kuznetsov, \emph{{MATCONT}: {A}
  {MATLAB} package for numerical bifurcation analysis of {ODE}s}, ACM Trans.
  Math. Software \textbf{29} (2003), no.~2, 141--164.

\bibitem{NewMATCONT}
A.~Dhooge, W.~Govaerts, Yu.~A. Kuznetsov, H.~G.~E. Meijer, and B.~Sautois,
  \emph{New features of the software {M}at{C}ont for bifurcation analysis of
  dynamical systems}, Math. Comput. Model. Dyn. Syst. \textbf{14} (2008),
  no.~2, 147--175.

\bibitem{sac2003}
A.~Dhooge, W.~Govaerts, {\mbox{Yu}}.~A. Kuznetsov, W.~Mestrom, and A.~M. Riet,
  \emph{{CL\_matcont}: {A} continuation toolbox in {Matlab}}, Symposium on
  Applied Computing (Melbourne, Florida), ACM, 2003, pp.~161--166.

\bibitem{AUTO97}
E.~J. Doedel, A.~R. Champneys, T.~F. Fairgrieve, {\mbox{Yu}}.~A. Kuznetsov,
  B.~Sandstede, and X.~J. Wang, \emph{{\sc auto97}: Continuation and
  bifurcation software for ordinary differential equations (with {H}om{C}ont)},
  1997.

\bibitem{DoGoKu:2003}
E.~J. Doedel, W.~Govaerts, and {\mbox{Yu}}.~A. Kuznetsov, \emph{Computation of
  periodic solution bifurcations in {ODEs} using bordered systems}, SIAM J.
  Numer. Anal. \textbf{41} (2003), no.~2, 401--435.

\bibitem{Faver:03}
S.~Fatimah and F.~Verhulst, \emph{Suppressing flow-induced vibration by
  parametric excitation}, Nonlinear Dynam. \textbf{31} (2003), 275--297.

\bibitem{GoGhKuMe:07}
W.~Govaerts, R.~{Khoshsiar Ghaziani}, {\mbox{Yu}}.~A. Kuznetsov, and H.~G.~E.
  Meijer, \emph{Numerical methods for two-parameter local bifurcation analysis
  of maps}, SIAM J. Sci. Comput. \textbf{29} (2007), no.~6, 2644--2667.

\bibitem{GuHo:83}
J.~Guckenheimer and {\mbox{Ph}}.~Holmes, \emph{Nonlinear {O}scillations,
  {D}ynamical {S}ystems and {B}ifurcations of {V}ector {F}ields},
  Springer-Verlag, New York, 1983.

\bibitem{Guckenheimer:2000}
J.~Guckenheimer and B.~Meloon, \emph{Computing periodic orbits and their
  bifurcations with automatic differentiation}, SIAM J. Sci. Comput.
  \textbf{22} (2000), no.~3, 951--985.

\bibitem{Io:79}
G.~Iooss, \emph{Bifurcation of {M}aps and {A}pplications}, North-Holland
  Mathematics Studies, vol.~36, North-Holland Pub. Co., Amsterdam, 1979.

\bibitem{Io:88}
\bysame, \emph{Global characterization of the normal form for a vector field
  near a closed orbit}, J. Differential Equations \textbf{76} (1988), no.~1,
  47--76.

\bibitem{IoAd:92}
G.~Iooss and M.~Adelmeyer, \emph{Topics in {B}ifurcation {T}heory and
  {A}pplications}, Advanced Series in Nonlinear Dynamics, vol.~3, World Sci.
  Pub. Co. Inc., River Edge, New York, 1992.

\bibitem{Jansen:1995}
V.A.A. Jansen, \emph{Regulation of predator-prey systems through spatial
  interactions: A possible solution to the paradox of enrichment}, Oikos
  \textbf{74} (1995), no.~3, 384--390.

\bibitem{Jansen:2001}
\bysame, \emph{The dynamics of two diffusively coupled predator-prey
  populations}, Theoretical Population Biology \textbf{59} (2001), no.~2, 119
  -- 131.

\bibitem{Ku:99}
{\mbox{Yu}}.~A. Kuznetsov, \emph{Numerical normalization techniques for all
  codim {$2$} bifurcations of equilibria in {ODE}s}, SIAM J. Numer. Anal.
  \textbf{36} (1999), no.~4, 1104--1124.

\bibitem{Ku:2004}
{\mbox{Yu}}.~A. Kuznetsov, \emph{Elements of {A}pplied {B}ifurcation {T}heory},
  Springer-Verlag, New York, 2004, 3rd~ed.

\bibitem{KuDoGoDh:05}
{\mbox{Yu}}.~A. Kuznetsov, W.~Govaerts, E.~J. Doedel, and A.~Dhooge,
  \emph{Numerical periodic normalization for codim 1 bifurcations of limit
  cycles}, SIAM J. Numer. Anal. \textbf{43} (2005), no.~4, 1407--1435.

\bibitem{CONTENT}
{\mbox{Yu}}.~A. Kuznetsov and V.~V. Levitin, \emph{{CONTENT}: {A} multiplatform
  environment for analyzing dynamical systems, {D}ynamical {S}ystems
  {L}aboratory, {CWI}, {A}msterdam}, 1997.

\bibitem{KuMe:2005}
{\mbox{Yu}}.~A. Kuznetsov and H.~G.~E. Meijer, \emph{Numerical normal forms for
  codim 2 bifurcations of fixed points with at most two critical eigenvalues},
  SIAM J. Sci. Comput. \textbf{26} (2005), no.~6, 1932--1954.

\bibitem{KuMe:2006}
\bysame, \emph{Remarks on interacting \mbox{Neimark-Sacker} bifurcations}, J.
  Difference Equ. Appl. \textbf{12} (2006), no.~10, 1009--1035.

\bibitem{KuMeGoSa:2008}
{\mbox{Yu}}.~A. Kuznetsov, H.~G.~E. Meijer, W.~Govaerts, and B.~Sautois,
  \emph{Switching to nonhyperbolic cycles from codim 2 bifurcations of
  equilibria in {ODE}s}, Phys. D \textbf{237} (2008), no.~23, 3061--3068.

\bibitem{Los:1989}
J.E. Los, \emph{Nonnormally hyperbolic invariant curves for maps in r 3 and
  doubling bifurcation}, Nonlinearity \textbf{2} (1989), no.~1, 149.

\bibitem{Rosenzweig:1971}
M.L. Rosenzweig, \emph{Paradox of enrichment: Destabilization of exploitation
  ecosystems in ecological time}, Science \textbf{171} (1971), no.~3969,
  385--387.

\bibitem{Simo}
C.~Sim\'o, \emph{Analytical and numerical computation of invariant manifolds},
  Modern Methods in Celestial Mechanics (C.~Benest and C.~Froeschl\'e, eds.),
  Editions Fronti\'{e}res, 1990, pp.~285--330.

\bibitem{ViBrSi:2011}
R.~Vitolo, H.~Broer, and C.~Sim{\'o}, \emph{Quasi-periodic bifurcations of
  invariant circles in low-dimensional dissipative dynamical systems}, Regul.
  Chaotic Dyn. \textbf{16} (2011), no.~1-2, 154--184.

\bibitem{BSV:2010}
R.~Vitolo, H.W. Broer, and C.~Sim\'o, \emph{Routes to chaos in the
  {H}opf-saddle-node bifurcation for fixed points of {3D}-diffeomorphisms},
  Nonlinearity \textbf{23} (2010), no.~8, 1919--1948.

\bibitem{WieChow:06}
S.~Wieczorek and W.~W. Chow, \emph{Self-induced chaos in a single-mode
  inversionless laser}, Phys. Rev. Lett. \textbf{97} (2006).

\end{thebibliography}

\end{document}